\documentclass{gtart_h}  

\def\ifplaintex{\expandafter\ifx\csname documentclass\endcsname\relax}

\def\ifplaintex{\expandafter\ifx\csname documentclass\endcsname\relax}

\ifplaintex 
\else
\fi

\expandafter\ifx\csname epsfbox\endcsname\relax\input epsf\fi

\def\gt{{\mathsurround=0pt\it $\cal G\mskip-2mu$eometry \&\ 
$\cal T\!\!$opology}}        %  journal title in recommended style

\def\gtp{{\mathsurround=0pt\it $\cal G\mskip-2mu$eometry \&\ 
$\cal T\!\!$opology $\cal P\!$ublications}}  % GT publications

\def\lognumber#1{\def\thelognumber{#1}}
\def\volumenumber#1{\def\thevolumenumber{#1}}
\def\papernumber#1{\def\thepapernumber{#1}}
\def\volumeyear#1{\def\thevolumeyear{#1}}

\def\pagenumbers#1#2{\def\startpage{#1}\def\finishpage{#2}}
\def\published#1{\def\publishdate{#1}}
\def\proposed#1{\def\theproposer{#1}}
\def\seconded#1{\def\theseconders{#1}}
\def\received#1{\def\receiveddate{#1}}
\def\revised#1{\def\reviseddate{#1}}
\def\accepted#1{\def\accepteddate{#1}}

\long\def\asciiabstract#1{\long\def\theasciiabstract{#1}}

\let\\\par\let\thelognumber\relax
\let\thevolumenumber\relax\let\thepapernumber\relax
\let\thevolumeyear\relax\let\thesamplenumber\relax\let\startpage\relax
\let\finishpage\relax\let\publishdate\relax\let\receiveddate\relax
\let\reviseddate\relax\let\accepteddate\relax\let\theasciititle\relax
\let\theasciiauthors\relax
\let\theasciiabstract\relax
\let\theasciiemail\relax\let\theshortauthors\relax\let\theshorttitle\relax

\long\def\maketitlep{   % start of definition of \maketitlep

\count0=\startpage

\gt\hfill      %   Journal title (top left) 
\hbox to 77pt{\vbox to 0pt{\vglue -15pt\epsfbox{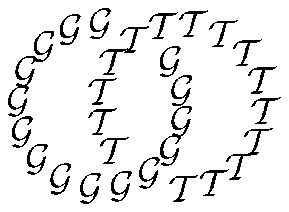}\vss}\hss}
\break
{\small\ifx\thesamplenumber\relax % sample?  
Volume \else Sample
\fi\thevolumenumber\ (\thevolumeyear)
\startpage--\finishpage\nl
Published: \publishdate}
\vglue 0.5truein plus 0.4fil minus 0.1truein

{\parskip=0pt\leftskip 0pt plus 1fil\def\\{\par\smallskip}{\ifplaintex\large
\else\Large\fi\bf\thetitle}\par\medskip}   

\vglue 0pt plus 0.1fil 

{\parskip=0pt\leftskip 0pt plus 1fil\def\\{\par}{\sc\theauthors}
\par\medskip}

\vglue 0pt plus 0.1fil 

{\small\parskip=0pt\let\newline\\
{\leftskip 0pt plus 1fil\def\\{\par}{\sl\theaddress}\par}
\expandafter\ifx\theemail\relax    % email address?
\relax\else\vglue 5pt plus 0.02fil minus 2pt\def\\{\stdspace{\rm 
and}\stdspace} 
\cl{Email:\stdspace\tt\theemail}\fi
\ifx\theurl\relax                  % URL given?
\relax\else\vglue 5pt plus 0.02fil minus 2pt\def\\{\stdspace{\rm 
and}\stdspace}
\cl{URL:\stdspace\tt\theurl}\fi\par}

\vglue 7pt plus 0.3fil minus 3pt

{\bf Abstract}
\vglue 5pt plus 0.1fil minus 2pt

\theabstract

\vglue 7pt plus 0.3fil minus 3pt

{\bf AMS Classification numbers}\quad Primary:\quad \theprimaryclass

Secondary:\quad \thesecondaryclass

\vglue 5pt plus 0.3fil minus 2pt

{\bf Keywords:}\quad \thekeywords

\vglue 10pt plus 0.5fil minus 5pt

{\small  Proposed: \theproposer\hfill Received: \receiveddate\nl
Seconded: \theseconders\hfill 
\ifx\reviseddate\relax                         % paper revised?
Accepted: \accepteddate                        % no
\else
Revised: \reviseddate                          % yes
\fi}
\eject
}       %  end of definition of \maketitlep

\font\phead=cmsl9 scaled 950
\font=cmsl9 scaled 1050
\font\pnum=cmbx10 scaled 913
\font\lnum=cmbx10 
\font\pfoot=cmsl9 scaled 950
\font=cmsl9 scaled 1050
\ifplaintex
\headline{\vbox to 0pt{\vskip -4.5mm\line{\small\phead\ifnum
\count0=\startpage ISSN 1364-0380 (on line)
1465-3060 (printed) \hfill {\pnum\folio}\else\ifodd\count0\def\\{ }% 
\ifx\theshorttitle\relax\thetitle\else\theshorttitle\fi\hfill{\pnum\folio}
\else\def\\{ and }{\pnum\folio}\hfill\ifx\theshortauthors\relax\theauthors
\else\theshortauthors\fi\fi\fi}\vss}}
\footline{\vbox to 0pt{\vglue 0mm\line{\small\pfoot\ifnum\count0=\startpage
\copyright\ \gtp\hfill\else
\gt, Volume \thevolumenumber\ (\thevolumeyear)\hfill\fi}\vss
}}
\else
\makeatletter
\def\@oddhead{{\small\lhead\ifnum\count0=\startpage ISSN 1364-0380 (on line)
1465-3060 (printed) \hfill {\lnum\number\count0}\else\ifodd\count0
\def\\{ }\ifx\theshorttitle\relax \thetitle \else\theshorttitle\fi\hfill
{\lnum\number\count0}\else\def\\{ and }{\lnum\number\count0}
\hfill\ifx\theshortauthors\relax 
\theauthors\else\theshortauthors\fi\fi\fi}}\def\@evenhead{\@oddhead}
\def\@oddfoot{\small\lfoot\ifnum\count0=\startpage\copyright\ \gtp\hfill\else
\gt, Volume \thevolumenumber\ (\thevolumeyear)\hfill\fi}
\def\@evenfoot{\@oddfoot}
\makeatother
\fi

\newwrite\gtoutfile
\long\gdef\makeheadfile{  %%% start of definition of \makeheadfile
{\def\\{, }\def\s{ }
\immediate\openout\gtoutfile head.xxx
\immediate\write\gtoutfile{Proxy-for: \ifx\theasciiauthors\relax
\theauthors\else\theasciiauthors\fi\s<\ifx\theasciiemail\relax\theemail\else\theasciiemail\fi>}
\immediate\write\gtoutfile{\noexpand\\}
\immediate\write\gtoutfile{Authors: \ifx\theasciiauthors\relax
\theauthors\else\theasciiauthors\fi}
{\def\\{ }\immediate\write\gtoutfile{Title: \ifx\theasciititle\relax
\thetitle\else\theasciititle\fi}}
\immediate\write\gtoutfile{Subj-class: GT or SG or MG etc}
\immediate\write\gtoutfile{MSC-class: \theprimaryclass\ifx\thesecondaryclass\relax\else, \thesecondaryclass\fi}
\immediate\write\gtoutfile{Journal-ref: Geom. Topol. \thevolumenumber
(\thevolumeyear) \startpage-\finishpage}
\immediate\write\gtoutfile{Comments: Published by Geometry and Topology at}
\immediate\write\gtoutfile{\s\s http://www.maths.warwick.ac.uk/gt/GTVol\thevolumenumber/paper\thepapernumber.abs.html}
\immediate\write\gtoutfile{\noexpand\\}
\immediate\write\gtoutfile{}
\ifx\theasciiabstract\relax
\immediate\write\gtoutfile{\theabstract}\else
\immediate\write\gtoutfile{\theasciiabstract}\fi
\immediate\write\gtoutfile{}
\immediate\write\gtoutfile{\noexpand\\}
\immediate\write\gtoutfile{}
\immediate\closeout\gtoutfile}}  %%% end of definition of \makeheadfile

\def\maketitlepage{\maketitlep\makeheadfile}
\let\maketitle\maketitlepage

\lognumber{518}
\received{10 November 2004}
\volumenumber{9}\papernumber{28}\volumeyear{2005}
\pagenumbers{1221}{1252}   
\revised{3 December 2004}
\published{24 July 2005}
\accepted{4 July 2005}
\proposed{Yasha Eliashberg}
\seconded{Robion Kirby, Joan Birman}

\usepackage{amsmath,amssymb,graphicx,labelfig} 
\usepackage{verbatim} % for comments

\newtheorem{theorem}{Theorem}[section]

\newtheorem{lemma}[theorem]{Lemma}

\newtheorem{proposition}[theorem]{Proposition}

\newtheorem*{conjecture}{Conjecture}
\newtheorem*{dualitytheorem}{Duality Theorem}
\newtheorem*{assumption}{Assumption}

\theoremstyle{definition}
\newtheorem{definition}[theorem]{Definition}
\newtheorem{remark}[theorem]{Remark}

\newtheorem{problem}[theorem]{Problem}

\newcommand{\im}{\operatorname{im}}
\newcommand{\ba}{\bar a}
\newcommand{\bb}{\bar b}
\newcommand{\bc}{\bar c}
\newcommand{\br}{\ensuremath{\mathbb R}}
\newcommand{\bz}{\ensuremath{\mathbb Z}}

\newcommand{\e}{\varepsilon}
\newcommand{\aug}{\operatorname{Aug}}
\newcommand{\ca}{\ensuremath{\mathcal{A}}}
\newcommand{\cc}{\ensuremath{\mathcal{C}}}
\newcommand{\cv}{C}

\newcommand{\df}{\partial}
\newcommand{\dfb}{\df_*}
\newcommand{\dfe}{\df_\e}
\newcommand{\Pb}{P_*}

\newcommand{\Pe}{P_\e}
\newcommand{\pb}{\pi_*}
\newcommand{\pe}{\pi_\e}
\newcommand{\rank}{\operatorname{rk}}
\newcommand{\ab}[1]{\langle#1\rangle}
\newcommand{\ket}[1]{|#1\rangle}
\newcommand{\braket}[2]{\langle#1|#2\rangle}
\newcommand{\ch}{ch}
\newcommand{\tb}{tb}
\def\S{Section~}
\newcommand{\fig}[2]{\includegraphics[scale=#2]{#1.eps}}

\newcommand{\foot}{\setcounter{footnote}{1}\footnote}\newcommand{\up}[1]{\raisebox{5ex}{$#1$}}
\newcommand{\tabletop}{\centering
\begin{tabular}{c r c c c r}
\hline\hline \\
knot & & reduced & & reduced &  \\ [-2ex]
type & \raisebox{1ex}{$(\tb,r)$}  & polynomials & \raisebox{1.5ex}{front \
and \ mirror}& polynomials & \raisebox{1ex}{$(\tb,r)$} \\ [1ex]
\hline \\ }

\begin{document}

\title{The nonuniqueness of Chekanov polynomials\\of Legendrian knots}
\author{Paul Melvin\\Sumana Shrestha}

\begin{abstract}
Examples are given of prime Legendrian knots in the standard contact
3--space that have arbitrarily many distinct Chekanov polynomials,
refuting a conjecture of Lenny Ng.  These are constructed using a new
``Legendrian tangle replacement" technique.  This technique is then
used to show that the phenomenon of multiple Chekanov polynomials is
in fact quite common.  Finally, building on unpublished work of Yufa
and Branson, a tabulation is given of Legendrian fronts, along with
their Chekanov polynomials, representing maximal Thurston--Bennequin
Legendrian knots for each knot type of nine or fewer crossings.  These
knots are paired so that the front for the mirror of any knot is
obtained in a standard way by rotating the front for the knot.
\end{abstract}

\asciiabstract{%
Examples are given of prime Legendrian knots in the standard contact
3-space that have arbitrarily many distinct Chekanov polynomials,
refuting a conjecture of Lenny Ng.  These are constructed using a new
`Legendrian tangle replacement' technique.  This technique is then
used to show that the phenomenon of multiple Chekanov polynomials is
in fact quite common.  Finally, building on unpublished work of Yufa
and Branson, a tabulation is given of Legendrian fronts, along with
their Chekanov polynomials, representing maximal Thurston-Bennequin
Legendrian knots for each knot type of nine or fewer crossings.  These
knots are paired so that the front for the mirror of any knot is
obtained in a standard way by rotating the front for the knot.}

\address{Department of Mathematics, Bryn Mawr 
College\\Bryn Mawr, PA 19010, USA}
\email{pmelvin@brynmawr.edu}

\maketitlepage

%%% SECTION 1 : INTRODUCTION %%%

\section{Introduction}
\label{sec1}

A smooth knot $K$ in $\br^3$ is {\it Legendrian} if it is everywhere
tangent to the two-plane distribution $\ker(\alpha)$ of the standard
contact one-form $\alpha=dz-ydx$.  It is tacitly assumed that the
projection of $K$ onto the $xz$--plane, its {\it front}, is generic,
that is all self-intersections of the front are transverse double
points.  Two Legendrian knots are {\it Legendrian isotopic} if they
are smoothly isotopic through (not necessarily generic) Legendrian
knots.  Much of the recent work in Legendrian knot theory has been
motivated by the problem of classifying Legendrian knots up to
Legendrian isotopy within a fixed knot type \cite{en}.

The two ``classical" Legendrian isotopy invariants of a Legendrian
knot $K$ are its {\it
Thurston--Bennequin number} $\tb(K)$ and its (absolute) {\it rotation
number} $r(K)$.
They can be computed by the formulas
$$
\tb(K) = w-c/2 \qquad\qquad  r(K) = |u-d|/2
$$
where $w$ is the writhe of the front of $K$ (the difference of the
number of positive and
negative crossings) and $c$ is the number of cusps in the front, of
which $u$ are
traversed upward and $d$ are traversed downward with respect to any chosen
orientation on $K$.   These two invariants serve to classify the
Legendrian knots in
some knot types, including the unknot \cite{ef}, the figure eight and
all torus knots \cite{eh1}.  Such knot types are called {\it
Legendrian simple}.

There are many knot types that are not Legendrian simple.  The first
examples, discovered independently by Chekanov \cite{c1} and
Eliashberg--Hofer \cite{eh}, were the positive twist knots (other than
the trefoil and figure eight).  The inequivalence of suitable Legendrian
representatives of these knots (with identical classical invariants)
was detected using {\it contact homology} \cite{h,egh}, a version
of Floer homology adapted to contact manifolds.  More recently, convex
surface techniques have been used to detect other examples among the
iterated torus knots \cite{eh3}.

Chekanov's approach to contact homology is purely combinatorial, and
yields an
easily computable family of integer Laurent polynomial invariants $\Pe(t)$
indexed by the {\it augmentations} $\e$ of (the {\it differential
graded algebra}
of) $K$; see \S\ref{sec2} for the definitions.  These will be called the {\it
Chekanov
polynomials} of $K$, and the number $\ch(K)$ of distinct such
polynomials will be
called the {\it Chekanov number} of $K$.

It has been observed by Chekanov \cite[\S2.2]{c2}, Ng
\cite[\S3.16]{n2} and \cite{n3}, and others that computational
evidence (at the time) supported the following conjecture.

\begin{conjecture}[Uniqueness of Chekanov polynomials]
The Che\-ka\-nov number $\ch(K)$ of any Legendrian knot $K$ is at most $1$.
\end{conjecture}

In this paper it is shown that this conjecture fails.  The first
counterexample we discovered, representing the mirror of the knot type
$8_{21}$, has Chekanov number $2$.  Knots with arbitrarily large Chekanov
numbers can then be constructed using connected sums.  With a little
more work, {\sl prime} knots of this kind can be found.

\begin{theorem}
\label{thm1.1}
For every integer $k\ge0$, there is a prime Legendrian knot with
Chekanov number $k$.
\end{theorem}

The proof is given in \S\ref{sec3}, and a conceptual  explanation of why these
examples arise is given in \S\ref{sec4} using a {\it Legendrian tangle
replacement} method.  It should be noted that these knots also
provide  counterexamples to the more speculative conjecture (cf
\cite[\S2.2]{c2}) that the {\sl  ranks} $P(1)$ are equal for all
Chekanov polynomials $P$ of a given Legendrian knot $K$; for the
knots in question, these ranks are in fact all distinct.  In
contrast, it is an easy exercise to show that the {\sl Euler
characteristic}
$
P(-1)=\tb(K)
$
for any Chekanov polynomial $P$ of
$K$.

Our computations throughout are simplified by a beautiful recent
result of Josh Sabloff \cite{s3}:

\begin{dualitytheorem}
Any Chekanov polynomial of a Legendrian knot can be written in the form
$$
P(t)=t+p(t)+p(t^{-1})
$$
for some (honest) polynomial $p(t)$ with positive integer coefficients.
\end{dualitytheorem}

Such a polynomial $p$ will be called a {\it reduced} Chekanov
polynomial of the knot.
Using this result and elementary observations about connected sums of
Legendrian knots, we give (in \S\ref{sec3}) a characterization of reduced Chekanov
polynomials.

\begin{theorem}
\label{thm1.2}
Every polynomial with positive integer coefficients arises as
the reduced Chekanov polynomial of some Legendrian knot.
\end{theorem}

It is an interesting problem to determine when $\ch(K)$ is nonzero.
Of course this can
be solved algorithmically, since the augmentations of $K$ are the
solutions of a system
of polynomial equations over $\bz_2$ that can be read off from the
front of $K$ (see \S\ref{sec2}).  Unfortunately, the known algorithms for
solving this system are all of
exponential complexity.  It would be useful to find more computable
criteria.

Fuchs and Ishkhanov \cite{fi}, and independently Sabloff \cite{s2},
have found a useful
geometric criterion that is equivalent to the existence of an
augmentation, namely the
existence of a {\it ruling} in the sense of Chekanov and Pushkar
\cite{c2}\cite{cp}.   The
obvious algorithms for finding a ruling, however, are still of
exponential complexity.

Another necessary (but far from sufficient) condition for the existence
of an
augmentation is that $K$ be {\it nondestabilizable}, ie not
Legendrian isotopic to a front with a kink $\fig{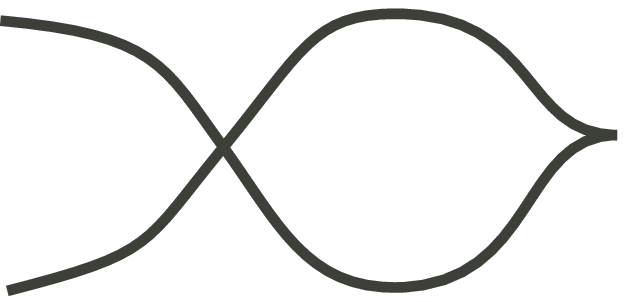}{.08}$.  An
important class of such knots are the ones of
maximal Thurston--Bennequin number in their knot types, which we
simply call {\it
maximal} knots.  (Note that there exist nondestabilizable knots that
are not maximal \cite{eh3}.)  All the knots tabulated in \S\ref{sec5} are
maximal.  Less that half of them, however,
have augmentations.  Moreover, there is no known algorithm for
deciding in the first
place if a knot is maximal (or nondestabilizable).

A less obvious (but readily computable) necessary condition is the
vanishing of the rotation number $r(K)$.  This was surmised from our
calculations in
\S\ref{sec5}, and confirmed by Sabloff as an easy consequence of the existence
of a ruling \cite{s2}.

\begin{proposition}
\label{prop1.3}
If $r(K)\ne0$, then $\ch(K)=0$.
\end{proposition}

Most of the knots with rotation number zero tabulated below have
augmentations,
although presumably this becomes less prevalent with increasing
crossing number.  Among these is the knot of type $7_4$, of Chekanov
number $1$.  In contrast, Ng \cite[\S4.2]{n2} has given examples of
knots in this knot type, with the same classical invariants, which do not
have augmentations.  Thus the Chekanov number can be used to distinguish
knots of the same type with the same classical invariants.

Note that in general, $\tb(K)$ and $r(K)$ are of opposite parity;
this can be proved in variety of ways, for example by induction on the
number of crossings
in the front of $K$.  It follows from the Proposition that only knots
with odd
Thurston--Bennequin invariant can have augmentations.

It remains an open problem to find complete criteria, computable in
polynomial time, for the existence of an augmentation.

\medskip
\textbf{Acknowledgments}\qua  We are grateful to Josh Sabloff for
introducing us to this
subject through his beautiful lectures in the TriCo Contact Seminar,
and for providing
us with an early version of his book on Legendrian knots.  We would
also like to thank
Lisa Traynor for stimulating discussions and graphics advice.
The first author is partially supported by
National Science Foundation grant FRG-0244460.

%%% SECTION 2 : CHEKANOV POLYOMIALS %%%

\section{Chekanov polynomials}
\label{sec2}

Let $K$ be a Legendrian knot with rotation number zero.   For
computational convenience we assume following Ng \cite{n1} that the
front of $K$ is {\it simple}, ie  all the right cusps have the
same $x$--coordinate.  This can be arranged by a Legendrian isotopy of
$K$.  Any smooth path in the front joining a left and right cusp will
be called a {\it spanning arc} of $K$ (so if there are $c$ cusps, then the
front consists of $c$ spanning arcs).  Any connected component of the
complement of the front in the $xz$--plane will be called a  {\it
region} in the front, and a disk formed from the closure of a union
of regions will be called a {\it disk with corners}.  The corners
refer to the nonsmooth points in the boundary that occur at
right cusps or crossings.

Now recall the definition of Chekanov's {\it differential graded
algebra} $(\ca,\df)$ for $K$.
The underlying algebra $\ca$ is the free noncommutative algebra over
$\bz_2$
generated by the crossings $c_1,\dots,c_n$ and right cusps
$c_{n+1},\dots,c_{n+r}$ in
the front.  Thus the elements of $\ca$ are finite sums of words in
the $c_i$, where the
empty word, denoted by $1$, is the identity.  The full set
$\{c_1,\dots,c_{n+r}\}$ of generators will be denoted by $\cc$.

The grading on $\ca$ is defined by assigning an integer degree
$|c|$ to each  $c\in\cc$, and then extending to higher order terms by
the rule $|ab|=|a|+|b|$.  To define $|c|$, choose a {\it  Maslov
potential} $\mu$ on the front of $K$.  By definition, $\mu$ assigns a
real number to  each spanning arc in the front in such a way that
the upper
arc at any cusp is assigned one more  than the lower arc; such an
assignment can be made consistently since $r(K)=0$.  Now set $|c|=1$
if $c$ is a cusp, and
$$
|c| = \mu(\alpha)-\mu(\beta)
$$
if $c$ is a crossing, where $\alpha$ is the upper (smaller slope) arc
at the crossing and $\beta$ is the lower arc.  We will denote the
subset of $\cc$ of generators of degree $k$ by $\cc_k$, and its
cardinality by $n_k$.

Finally the differentials $\df c$ for $c\in\cc$ are defined using
suitable embedded disks, and then extended by the Leibnitz rule and
linearity, setting $\df1=0$, to all of $\ca$.  (If the front is not
simple, one must consider immersed disks.)  In particular, a disk $D$
with corners is called an {\it admissible disk for $c$} if its
right-most corner is at $c$, its left-most point is at a left cusp,
and the remaining corners, called the {\it negative} corners of $D$
for reasons explained in \cite{s1}, are {\it convex} (meaning that
they occur at crossings that lie on the boundary of {\sl exactly} one
of the regions in $D$).  There is an associated monomial $\df
D\in\ca$ obtained by reading off the labels at the negative corners
of $D$, proceeding from $c$ counterclockwise around the boundary.
Now define
$$
\df c = \begin{cases}
\sum\df D & \text{if $c$ is a crossing}\\
1+\sum\df D & \text{if $c$ is a right cusp}
\end{cases}
$$
where the sums are over all admissible disks $D$ for $c$.  It is not
difficult to show that
$\df$ lowers degree by 1, and that $\df^2=0$.   The graded homology
$H_*(\ca,\df)$, which coincides with the contact homology of the
standard contact $\br^3$ relative to $K$ \cite{ens}, is invariant
under Legendrian isotopy of $K$.

The algebra $(\ca,\df)$ is infinite dimensional (over $\bz_2$) and
generally cumbersome
to deal with as a whole, but useful information can be
extracted from its finite dimensional quotients.  In particular,
consider projections $\pb\colon\ca\to C$ where $C$ is the
finite dimensional space spanned by $\cc$.  Under suitable conditions
the induced endomorphism $\dfb=\pb\df$ of $C$ is a differential (that
is $\dfb^2=0$) and the associated Poincar\'e polynomial $\Pb(t) =
\sum_k \dim\left(H_k(C,\dfb)\right)t^k$ is a Legendrian isotopy
invariant.   The Chekanov polynomials $\Pe$ arise in this way from
projections $\pe\colon\ca\to\cv$ associated with
augmentations $\e$ of $(\ca,\df)$.  Here are the details.

\begin{definition}
\label{def2.1}
An {\it augmentation} of $(\ca,\df)$ (also
called an augmentation of $K$) is an algebra map
$$
\e\colon\ca\to\bz_2
$$
that vanishes on elements of nonzero degree and satisfies $\e\df=0$.
The generators $c\in\cc$ with $\e(c)=1$, which are all crossings of
degree zero, will be called the {\it augmented crossings} of $\e$.

The associated  {\it projection $\pe$} is defined on monomials
$m=c_{i_1}\ldots c_{i_k}$  (which form a  basis for $\ca$ over
$\bz_2$) by extracting the linear term in
$$(c_{i_1} + \e(c_{i_1})) \ldots (c_{i_k} + \e(c_{i_k}))$$
where $\pe(1)=0$ by convention.
Thus
$$
\pe(m)=c_{i_1}+\cdots +c_{i_k}
$$
if all the factors
$c_{i_j}$ are augmented.  Such a monomial will be called {\it pure}.
If all but one of the factors,  say $c$, are augmented, then
$\pe(m)=c$.  In this case we say that $m$ is {\it full}.  In all
other cases $\pe(m)=0$.

The projection $\pe$ induces a graded
differential $\dfe=\pe\df$ on $C$, or more precisely a chain
complex
$$
\cdots \longrightarrow C_{k+1}
\overset{\df_{k+1}}{\longrightarrow} C_k
\overset{\df_k}{\longrightarrow} C_{k-1} \longrightarrow \cdots
$$
where $C_k$ is the vector space spanned by the set $\cc_k$ of
generators of degree $k$, and $\df_k = \dfe|C_k$.  Explicitly, the
differential $\dfe b$ for any $b\in\cc$ is the sum of all the factors
in all the pure monomials in $\df b$ added to the unaugmented factors
in the full monomials in $\df b$.

The {\it Chekanov polynomial}
$\Pe$ is the Poincar\'e polynomial of this complex,
$$
\Pe(t) = \sum
\dim(H_k)\,t^k
$$
where $H_k =\ker\df_k / \im\df_{k+1}$.
\end{definition}

\begin{remark}
\label{rem2.2}
An augmentation $\e$  is uniquely determined by its set of augmented
crossings, which by  abuse of notation will also be denoted by $\e$.
The condition  $\e\df=0$ simply asserts that (the set) $\e$ is the support
of a solution to the system $\{\df  c=0:c\in\cc\}$ of (commutative)
polynomial equations over $\bz_2$.   Thus the set $\aug(K)$ of all
augmentations can be viewed as the set of all subsets of $\cc_0$ which
are supports of solutions to this system.

This condition can be interpreted in simple geometric terms as follows:
Start with an arbitrary subset $\e$ of $\cc_0$.  In analogy with the
terminology above, we say that an admissible disk for a generator
$c\in\cc$ is a {\it pure} disk (with respect to $\e$) if all of its
negative corners are in $\e$, and {\it full} if all but one are in
$\e$.  Let $\ket c\in\bz_2$ denote the mod 2 reduction of either the
number of pure disks for $c$, or one more than that number, according
to whether $c$ is a crossing or a right cusp.   Then, noting that
$\ket c=0$ if $|c|\ne1$, we have
$$
\e\in\aug(K) \iff \ket c=0 \text{ for
all } c\in\cc_1.
$$
\indent In the same vein, the differential
$\dfe$ for an augmentation $\e$ can be defined by
$$
\dfe c =\textstyle\sum_b \braket bc b
$$
where $\braket bc$ is the number of
admissible disks for $c$ which have a corner at $b$ and for which all
remaining corners are augmented.  (These disks are all pure if $b$ is
augmented, and all full otherwise.)
\end{remark}

\begin{remark}
\label{rem2.3}
The coefficient of $t^k$ in $\Pe(t)$ can be expressed as
$$
n_k-r_k-r_{k+1}
$$
where $n_k=\dim C_k = |\cc_k|$ and $r_k=\rank\df_k$.  It follows that
$\Pe$ is determined by the ranks $r_k$ for all $k>0$.  Indeed these
give the coefficients for all $k>0$, and thence for $k<0$ by the
Duality Theorem (see \S\ref{sec1}).  The constant coefficient is then computed
using the
identity
$
\Pe(-1) = \textstyle\sum_k(-1)^k n_k.
$
Note that  this last sum also computes $\tb(K)$, which is odd since
$r(K)=0$ is
even.

By the same argument $\Pe$ is determined by the ranks $r_k$
for all $k<0$.  Thus if $|c|\ge-1$ for all $c\in\cc$, as for the
examples in the next section, then the reduced polynomial $p_\e$ is
linear and is determined by the single rank $r_0$:
$$
p_e(t) =
at+b
$$
where $a=n_{-1}-r_0$ and $b=a+(\tb(K)+1)/2$.  In particular,
if there are no crossings of negative degree then $p_\e$ is
constant.
\end{remark}

%%% SECTION 3 : EXAMPLES %%%

\section{Examples}
\label{sec3}

\begin{figure}[ht!]\small
%\centerline{\fig{31m}{.5}}
\centerline{
%\ShowGrid
\SetLabels
  \T(0.5*0.0) $0$ \\
    (0.3*0.6) $1$ \\
    (0.7*0.6) $1$ \\
  \T(0.5*0.9) $2$ \\
\L\T(1.0*0.7) $r$ \\
\L\B(1.0*0.2) $s$ \\
\R\T(0.25*0.4) $a$ \\
  \T(0.5*0.4) $b$ \\
  \T(0.75*0.4) $c$ \\
\endSetLabels
\AffixLabels{\includegraphics[scale=0.5]{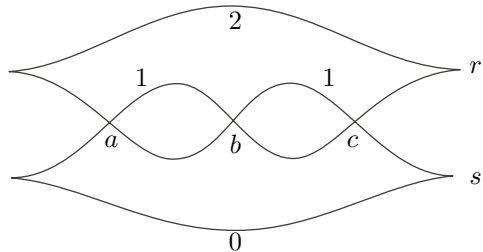}}
}
\caption{Front for the knot $K_1$
representing the right-handed trefoil}
\label{fig1}
\end{figure}

First a familiar computation.  Let $K_1$ be the
Legendrian trefoil whose front is shown in Figure~\ref{fig1}.

%%% FIG 1 %%%

%%%%%%%%%

The graded generating sets are $\cc_1=\{r,s\}$ (the
right cusps) and $\cc_0=\{a,b,c\}$ (computed using the Maslov
potential indicated by the numbers above each arc in the diagram).
The differentials of $a$, $b$ and $c$ vanish since there are no
admissible disks for these crossings, while
$$
\df r = 1+[abc]
\quad\text{and}\quad \df s = 1+[cba]
$$
where by definition
$[xyz]=x+xyz+z$.

The augmentations of $K_1$ correspond to solutions to the system $\df
r = \df s = 0$, or equivalently to the single equation $[abc]=1$ (since
$[abc]=[cba]$ after abelianizing).  By inspection
$$\aug(K_1) = \{abc,ab,a,bc,c\}$$
where for notational economy we write  $abc, ab, \ldots$ for the sets
$\{a,b,c\}, \{a,b\}, \ldots$  (These five subsets will be referred to below
as the {\it admissible} subsets of the ordered set or ``triple" $abc$.)
By direct computation, or an appeal to the duality theorem as in
Remark~\ref{rem2.3} above, all of these augmentations yield the same
Chekanov polynomial $t+2$, or equivalently the same reduced polynomial
$1$.  In particular, $K_1$ has Chekanov number $1$.

\begin{figure}[ht!]\small
\centerline{
%\ShowGrid
\SetLabels
\L\T(0.5*0.09) $0$ \\
\L\B(0.5*0.18) $1$ \\
\L\T(0.5*0.38) $2$ \\
\L\B(0.5*0.51) $1$ \\
\L\T(0.5*0.68) $2$ \\
\L\T(0.5*0.92) $3$ \\
\B(0.15*0.2) $\bar a$ \\
\B\L(0.3*0.2) $\bar b$ \\
\R\T(0.8*0.28) $\bar c$ \\
\L\T(0.1*0.62) $a$ \\
\T(0.3*0.62) $b$ \\
\L\B(0.8*0.6) $c$ \\
\T(0.4*0.5) $p$ \\
\T\R(0.7*0.5) $q$ \\
\T\L(1.0*0.7) $r$ \\
\B\L(1.0*0.4) $s$ \\
\L(1.0*0.2) $t$ \\
\endSetLabels
\AffixLabels{\includegraphics[scale=0.5]{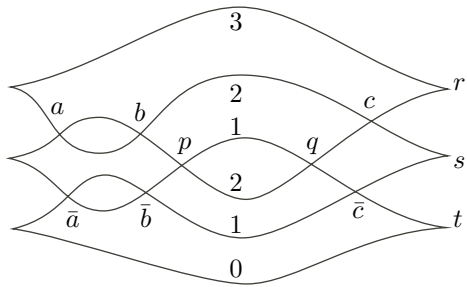}}
}
\caption{Front for the knot $K_2$
representing $\bar8_{21}$}
\label{fig2}
\end{figure}

Next  consider the Legendrian knot $K_2$ given by the front in
Figure~\ref{fig2},
which is topologically the mirror image of the knot $8_{21}$ (also
denoted $\bar 8_{21}$).

%%% FIG 2 %%%

%%%%%%%%%

Here  $\cc_1=\{p,r,s,t\}$, $\cc_0=\{a,b,c,\ba,\bb,\bc\}$ and
$\cc_{-1}=\{q\}$, and the nonzero differentials are
$$\begin{array}{lll}
\df c =
a[\ba\bb]q & \df \bc = q[ba]\ba & \df p = [ba][\ba\bb] \\
\df r =
1+[abc]+apq & \df s = [a\ba]+a[\ba\bb\bc]+[cba]\ba & \df t =
1+[\bc\bb\ba]+qp\ba
\end{array}$$
where by definition $[xy]=1+xy$.  Note that two of the three appearances
of $a\ba$ in $\df s$ can be cancelled since we are working mod 2.
The differential is written in this way in order to facilitate the
calculation of augmentations, which amounts to solving four equations
over $\bz_2$:  the vanishing of $\df p,\ \df r,\ \df s$ and $\df t$.

The solutions to these equations can be described simply by employing
the following terminology, which will be used throughout the rest of
the paper.

\begin{definition}
\label{def3.1}
The {\it admissible} subsets of a triple $xyz$ are $xyz,\ xy,\ x,\ yz$
and $z$.  The first two will be referred to as {\it ample}, since they
contain both $x$ and $y$, and the remaining three will be called {\it
sparse}.  Note that the admissible subsets correspond to solutions to
$[xyz]=1$, and the additional constraint $[xy]=1$ selects the sparse ones.
\end{definition}

Now consider the two level triples $abc$, $\ba\bb\bc$ of degree zero
crossings of $K_2$.  The subsets of these triples that are augmented
for a given augmentation $\e$ will be called the {\it levels} of $\e$.
Now the vanishing of $\df r$ and $\df t$ simply says that both levels
of $\e$ must be admissible, and $\df p = 0$ then forces at least one of
these to be ample.   The last equation $\df s = 0$ imposes no further
restrictions as it is a consequence of the other three equations.
It follows that $K_2$ has sixteen augmentations, four of which
$$
abc\ba\bb\bc,\ abc\ba\bb,\ ab\ba\bb\bc,\ ab\ba\bb
$$
have no sparse levels, while the remaining twelve
$$
abc\ba,\ abc\bb\bc,\ abc\bc,\ ab\ba,\ ab\bb\bc,\ ab\bc,\ a\ba\bb\bc,\
a\ba\bb,\ bc\ba\bb\bc,\ bc\ba\bb, \ c\ba\bb\bc, \ c\ba\bb
$$
have exactly one.  It is now straightforward (although tedious without
a computer) to show that the reduced Chekanov polynomial is $t+2$ for
the augmentations in the first group, and $1$ for those in the second.
Alternatively this can be proved using the Duality Theorem (see the last
remark in \S\ref{sec2}):  It suffices to show that an augmentation is in the first
group if and only if the rank $r_0$ of the differential $C_0\to C_{-1}$
vanishes, or equivalently $a[\ba\bb] = [ba]\ba = 0$.  But $a = b = \ba =
\bb = 1$ is clearly the only solution to these equations.  Thus $K_2$
has Chekanov number $2$.

\begin{figure}[ht!]\small
\centerline{
%\ShowGrid
\SetLabels
(0.05*0.45) $K$ \\
(0.207*0.45) $\#$ \\
(0.36*0.45) $K'$ \\
\R(0.5*0.45) $=$ \\
(0.6*0.45) $K$ \\
(0.75*0.45) $c$ \\
(0.95*0.45) $K'$ \\
\endSetLabels
\AffixLabels{\includegraphics[scale=0.6]{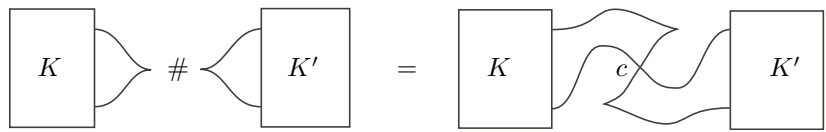}}
}
\caption{Legendrian connected sum $K\# K'$}
\label{fig3}
\end{figure}

Knots with arbitrarily large Chekanov number can now be
constructed by taking connected sums of copies of $K_2$.  To see
this, recall that the connected sum $K\# K'$ of Legendrian knots can
be formed in a variety of ways, all equivalent up to Legendrian
isotopy \cite{eh2}.  A convenient one for our purposes is shown in
Figure~\ref{fig3}.

%%% FIG 3 %%%

%%%%%%%%%

With this choice, the front of $K\# K'$
has one extra degree zero crossing $c$.  Its differential graded
algebra is generated by $\cc\cup\cc'\cup\{c\}$, with augmentations
$$
\e\#\e' = \e\cup\e'\cup\{c\}
$$
where $\e$ and $\e'$ range over
all the augmentations of $K$ and $K'$.  Fixing such an augmentation
yields the chain complex
$$
C\# C' = C\oplus C' \oplus
\bz_2
$$
where $C$ and $C'$ are the chain complexes for $\e$ and
$\e'$, and the last factor is generated by $c$.  The $d$th
differential for $C\# C'$ is the block sum of the corresponding
differentials for $C$ and $C'$, with an extra column of zeros when
$d=0$, and an extra nonzero row when $d=1$.  Its rank is therefore
the sum of the corresponding ranks $r_d$ and $r'_d$, except possibly
when $d=1$.  In fact in this case the rank is $r_1+r'_1+1$, as can be
seen using the duality theorem.  It follows that
$$
P_{\e\#\e'}(t) =
P_\e(t)+P_{\e'}(t) - t,
$$
or equivalently, the reduced Chekanov
polynomials add
$
p_{\e\#\e'} = p_{\e} + p_{\e'}.
$
This proves the
following (well-known) result, cf \cite[\S12]{c1}.

\begin{lemma}
\label{lem3.2}
The reduced Chekanov polynomials of a Legendrian connected sum
are exactly the sums of the reduced Chekanov polynomials of its
factors. \qed
\end{lemma}

Thus the set of reduced Chekanov polynomials
of the connected sum $nK_2$ of $n$ copies of $K_2$ is $\{kt+n+k :
k=0,\dots,n\}$.  In particular  $\ch(nK_2) = n+1$.

\begin{figure}\small
\centerline{
%\ShowGrid
\SetLabels
\L\B(1.0*0.5) $d{+}1$ crossings \\
\endSetLabels
\AffixLabels{\includegraphics[scale=0.6]{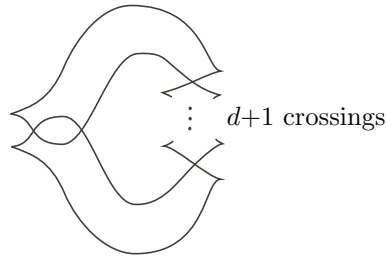}}
}
\caption{Twist knot $T_d$}
\label{fig4}
\end{figure}

Connected sums can also be used to construct Legendrian knots with
prescribed
Chekanov polynomials.  For example, noting that the twist knot
$T_{d}$ in Figure~\ref{fig4} has $t^d$ as its unique reduced Chekanov
polynomial (as is easily verified), the knot $\#_d (a_dT_d)$ is seen
using the lemma to have (unique) reduced polynomial $\sum_d a_dt^d$.
This proves Theorem~\ref{thm1.2} in the Introduction.

%%% FIG 4 %%%

%%%%%%%%%

One must work a bit harder to produce {\sl prime} knots
with large Chekanov numbers or with prescribed Chekanov polynomials.
We show how to do the former, as claimed in Theorem~\ref{thm1.1}.

\medskip
{\bf Proof of Theorem \ref{thm1.1}}\qua
Let $K_0$ be the left-handed trefoil in the tables below (representing
$3_1$).  This
knot has Chekanov number zero, as can be seen by direct calculation
or by appealing to the fact that its rotation number is nonzero (see
Proposition~\ref{prop1.3}).  It is also prime, even in the Legendrian
sense (since it has maximal Thurston--Bennequin number).

\begin{figure}\small
\centerline{
%\ShowGrid
\SetLabels
(0.08*0.70) $a_1$ \\
(0.30*0.66) $b_1$ \\
(0.78*0.66) $c_1$ \\
(0.08*0.48) $a_2$ \\
(0.31*0.41) $b_2$ \\
(0.79*0.43) $c_2$ \\
(0.08*0.23) $a_3$ \\
(0.33*0.17) $b_3$ \\
(0.79*0.19) $c_3$ \\
\L(1.00*0.80) $s_0$ \\
\L(1.00*0.58) $s_1$ \\
\L(1.00*0.36) $s_2$ \\
\L(1.00*0.15) $s_3$ \\
(0.42*0.56) $p_1$ \\
(0.66*0.56) $q_1$ \\
(0.43*0.27) $p_2$ \\
(0.67*0.28) $q_2$ \\
\endSetLabels
\AffixLabels{\includegraphics[scale=0.5]{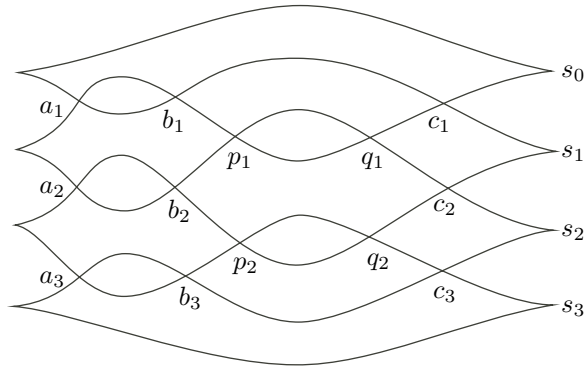}}
}
\caption{Front for the knot $K_3$}
\label{fig5}
\end{figure}

It is well-known that the knots $K_1$ and $K_2$ above, with Chekanov
numbers $1$ and $2$ respectively, are prime.  There is a natural way
to view them as the first two knots in a sequence $K_n$, noting that
each can be described as the ``plat" closure of a positive braid (see
eg \cite{bi}) where the plats correspond to cusps in the front:
$K_1$ is the closure of the $4$--braid $\sigma_2^3$, and $K_2$ is the
closure of the $6$--braid $(\sigma_2\sigma_4)^2 (\sigma_3)^2
(\sigma_2\sigma_4)$.  For $n\ge3$, let $K_n$ be the closure of the
$(2n+2)$--braid $e_n^2 o_n^2 e_n$, where $e_n = \sigma_2\sigma_4
\ldots \sigma_{2n}$ and $o_n = \sigma_3\sigma_5 \cdots
\sigma_{2n-1}$.  The case $n=3$ is shown in Figure~\ref{fig5}.

%%% FIG 5 %%%

%%%%%%%%%

The labels in Figure~\ref{fig5} suggest a general procedure for labeling $K_n$.
The degree zero crossings $a_i,b_i,c_i$  occur in descending level triples
(one for each $i=1,\dots,n$) while the inner crossings $p_j,q_j$, of
degrees $\pm1$ respectively, occur in descending level pairs (one for each
$j=1,\dots,n-1$).  The right cusps are labeled $s_j$ for $j=0,\dots,n$
from top to bottom.  Setting variables with subscripts outside their
prescribed ranges equal to $1$, the nonzero differentials are
$$
\begin{aligned}
\df p_j &= [b_ja_j]  [a_{j+1}b_{j+1}]
\\
\df c_i &= q_{i-1} [b_{i-1}a_{i-1}] a_i  + a_i [a_{i+1}b_{i+1}]  q_i
+ \cdots
\\
\df s_j &= [a_ja_{j+1}] + a_j [a_{j+1}b_{j+1}c_{j+1}] +
[c_jb_ja_j]a_{j+1} + \cdots
\end{aligned}
$$
where as above $[xy]=1+xy$ and $[xyz]=x+xyz+z$.  Here all monomials with
more than one factor of nonzero degree are suppressed, since they never
contribute to the differential associated with an augmentation.

Arguing as above it is seen that the augmentations of $K_n$ correspond
to subsets of the set of degree zero crossings whose levels (ie
augmented subsets of the level triples $a_ib_ic_i$) are admissible, and
whose sparse levels are nonadjacent.  Indeed, the equation $\df s_0 = 1 +
[a_1b_1c_1]+\cdots = 0$ establishes the admissibility of the top level.
Substituting this into $\df s_1 = 0$ gives $[a_1a_2] + a_1[a_2b_2c_2] +
a_2 = 0,$ which combined with $\df p_1 = 0$ establishes the admissibility
of the second level and forces it to be ample if the top level is sparse.
Similarly $\df s_2 = \df p_2 = 0$ then gives the admissibility of the
third level, and forces it to be ample if the level above is sparse.
Proceeding in this way gives the desired conclusion.

By the last remark in \S\ref{sec2}, the reduced Chekanov polynomial for $\e \in
\aug(K_n)$ is
$$
p_\e(t) = (n-1-r_0)t+(n-r_0)
$$
where $r_0$ is the rank  of the associated differential
$$\df_0 \colon C_0 \to C_{-1} =  \ab{q_1,\dots,q_{n-1}}.$$
It is shown below that $r_0$
can range from $0$ to $\lfloor 2n/3 \rfloor$, depending on the positions
of the sparse levels of $\e$, and so
$
\ch(K_n) =  \lfloor 2n/3 \rfloor + 1.
$
In particular
$$
\ch(K_{\lfloor 3k/2 \rfloor -1}) = k
$$
for any given integer $k$.

To compute $r_0$, observe that it is the rank of the
$(n-1){\times}n$--matrix $Q$ representing the restriction of $\df_0$
to the
subspace $\ab{c_1,\dots,c_n}$.  Now focus on the $i$th level for some
$i=1,\dots,n$.   The formula for $\df c_i$ above shows that $\df_0 c_i  =
q_i$ if the level above is sparse, $q_{i-1}$ if the level below is sparse,
$q_i+q_{i-1}$ if both are sparse, and $0$ otherwise.   It follows that
the nonzero rows of $Q$ are all standard basis vectors for $\br^n$, and
any particular one $e_j$ arises if and only if either level adjacent to
the $j$th one is sparse.  Thus $r_0$ is just a count of the number of
levels adjacent to the sparse levels, and an inductive argument shows
that this number assumes all possible values between $0$ and $\lfloor
2n/3 \rfloor$, as claimed.

It remains to show that the knots $K_n$ are prime.   One way to do this
is to use the notion of {\it prime tangles} introduced by Kirby and
Lickorish in \cite{kl}.  The reader is referred to Lickorish \cite{l},
where prime tangles were first used for primality testing in knot theory,
for the relevant definitions.  The main result of that paper is that
any knot that can be written as a sum of two prime tangles is prime,
and so it suffices to show that $K_n$ is of that form for $n>2$.

\begin{figure}[ht!]\small
\centerline{
%\ShowGrid
\SetLabels
(0.10*0.43) $S$ \\
(0.28*0.43) $T$ \\
(0.73*0.43) $S$ \\
(0.93*0.43) $T$ \\
\L(0.12*-0.15) (a) sum \\
\L(0.70*-0.15) (b) partial sum \\
\endSetLabels
\AffixLabels{\includegraphics[scale=0.5]{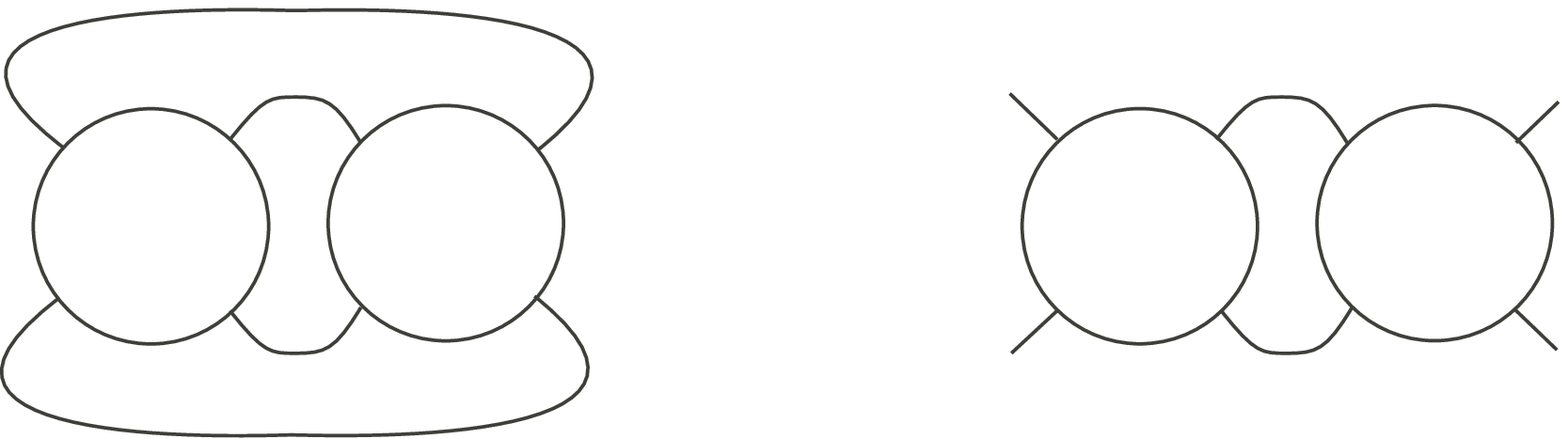}}
}
\vskip 5pt
\caption{Tangle sums}
\label{fig6}
\end{figure}

It should be noted that in \cite{l}, a {\it sum} of tangles $S$ and
$T$ refers to any link obtained by gluing the tangles together by a
homeomorphism of their boundaries, as in Figure~\ref{fig6}(a), rather
than gluing them together along only part of their boundaries, as in
Figure~\ref{fig6}(b).  The latter produces another tangle called a {\it
partial sum} (or {\it tangle sum}) of $S$ and $T$.  To distinguish the
two, write $S\oplus T$ for a sum and $S+T$ for a partial sum.

%%% FIG 6 %%%
%%%%%%%%%

\begin{figure}\small
\centerline{
%\ShowGrid
\SetLabels
(0.14*0.80) $P$ \\
(0.90*0.40) $T_n$ \\
\endSetLabels
\AffixLabels{\includegraphics[scale=0.4]{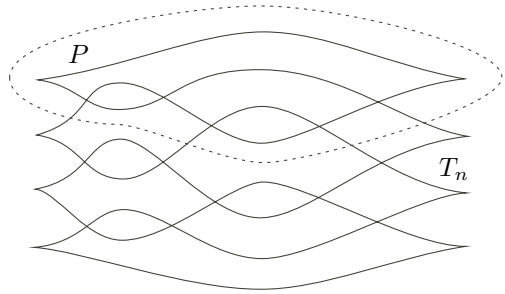}}
}
\caption{Decomposition of $K_n$}
\label{fig7}
\end{figure}

Now the decomposition of a knot as a sum of two tangles can be specified
by drawing a simple closed curve in a projection plane that intersects
the knot in four points, while proper arcs can be used to specify the
decomposition of a tangle as a partial sum of tangles.  So consider
the decomposition $K_n = P\oplus T_n$ coming from a simple closed curve
enclosing the top five crossings in $K_n$, as shown in Figure~\ref{fig7}
for the case $n=3$.

%%% FIG 7 %%%
%%%%%%%%%

\begin{figure}[ht!]\small
\centerline{
%\ShowGrid
\SetLabels
\L(0.05*-0.15) (a) Legendrian form \\
\L(0.60*-0.15) (b) topological form \\
\endSetLabels
\AffixLabels{\includegraphics[scale=0.5]{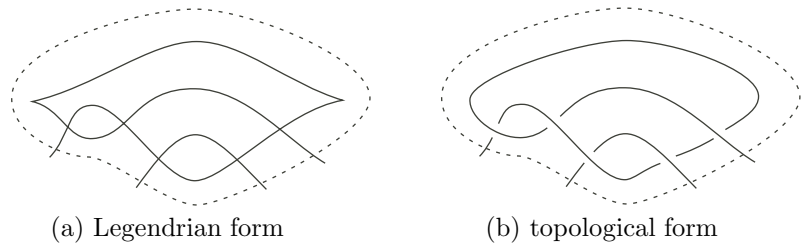}}
}
\vskip 10pt
\caption{The prime tangle $P$}
\label{fig8}
\end{figure}

It is easy to verify that the top tangle $P$, shown in its Legendrian
form in Figure~\ref{fig8}(a) or its topological equivalent in
Figure~\ref{fig8}(b), is prime
(see \cite[Figure~2a]{l}).

%%% FIG 8 %%%
%%%%%%%%%

\begin{figure}\small
\centerline{
%\ShowGrid
\SetLabels
(0.15*0.30) $T$ \\
(0.85*0.30) $T$ \\
\L(-0.10*-0.15) (a) adding a clasp \\
\L(0.60*-0.15) (b) adding a twist \\
\endSetLabels
\AffixLabels{\includegraphics[scale=0.5]{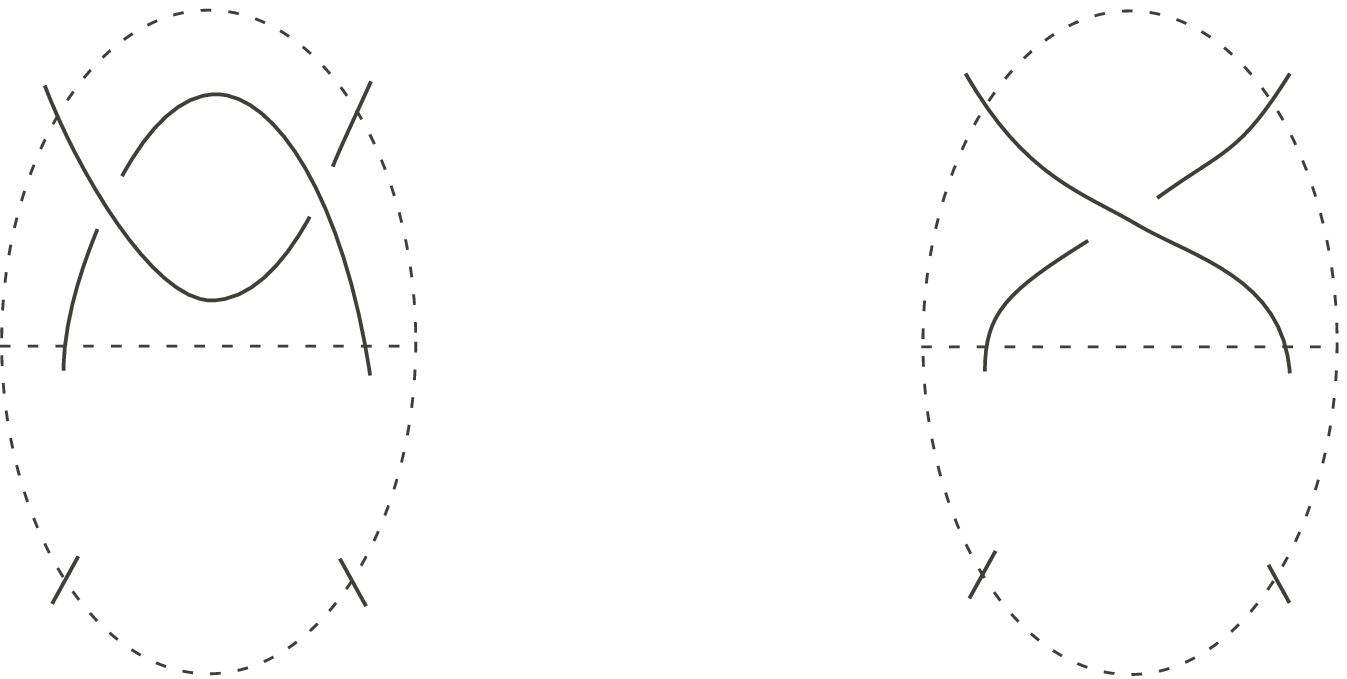}}
}
\vskip 10pt
\caption{Adding clasps and twists}
\label{fig9}
\end{figure}

The bottom tangle $T_n$ is also prime.  Indeed it can be written as
a partial sum $S_0+S_1+S_2+\cdots$, where $S_0$ is the prime tangle
(equivalent to $P$) encompassing the bottom five crossings in $T_n$,
and the $S_i$ are trivial (aka rational) tangles added so as to
appear as ``clasps" or ``twists", as shown in Figure~\ref{fig9}.

%%% FIG 9 %%%
%%%%%%%%%

\begin{figure}\small
\centerline{
%\ShowGrid
\SetLabels
(0.14*0.28) $S_0$ \\
(0.12*0.55) $S_1$ \\
(0.81*0.45) $S_2$ \\
(0.53*0.63) $S_3$ \\
(0.12*0.82) $S_4$ \\
(0.82*0.71) $S_5$ \\
(0.97*0.70) $T_4$ \\
\endSetLabels
\AffixLabels{\includegraphics[scale=0.5]{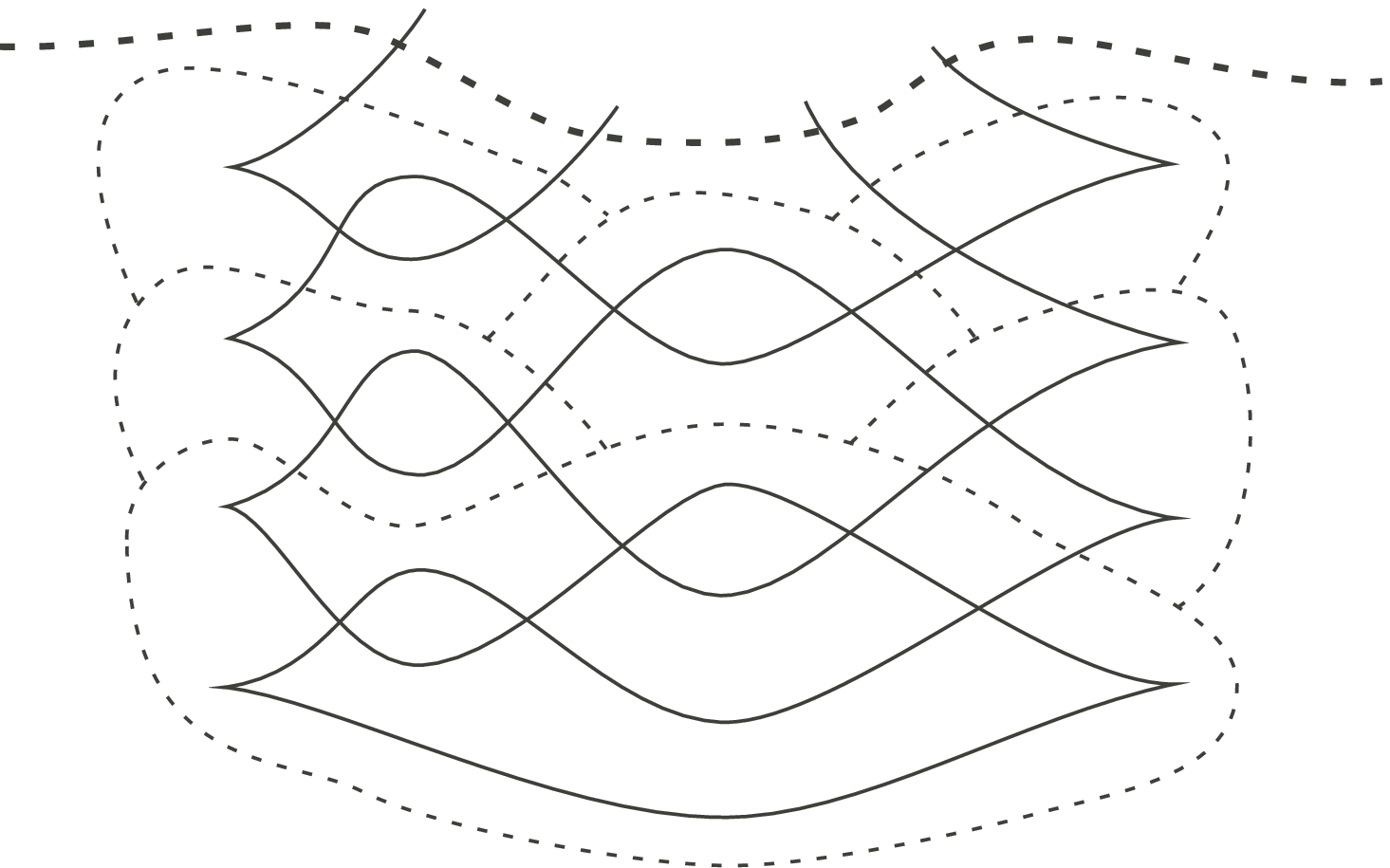}}
}
\caption{The decomposition of $T_4$}
\label{fig10}
\end{figure}

The case $n=4$ is illustrated in Figure~\ref{fig10}.

%%% FIG 10 %%%
%%%%%%%%%

\noindent Since adding a twist does not change the equivalence class of
the tangle, and straightforward (innermost disk) arguments show that
adding a clasp preserves primality, it follows that $T_n$ is prime.
Thus for $n>2$, the knot $K_n$ is a sum of two prime tangles, and is
therefore a prime knot.  This completes the proof of Theorem~\ref{thm1.1}.

%%% SECTION 4 : LEGENDRIAN TANGLE REPLACEMENT %%%

\section{Legendrian tangle replacement}
\label{sec4}

\begin{figure}\small
\centerline{
%\ShowGrid
\SetLabels
(0.25*0.22) $e$ \\
(0.36*0.17) $t$ \\
(0.36*0.41) $s$ \\
(0.60*0.48) $a$ \\
(0.68*0.47) $b$ \\
(0.86*0.48) $c$ \\
(0.725*0.37) $p$ \\
(0.815*0.37) $q$ \\
(0.96*0.66) $r$ \\
(0.98*0.36) $s$ \\
(0.98*0.12) $t$ \\
(0.88*0.18) $e$ \\
(0.13*0.14) $K$ \\
(0.73*0.12) $\tau(K)$ \\
(0.22*0.72) $T$ \\
(0.92*0.76) $P$ \\
\L(-0.10*-0.15) (a) Special Legendrian knot \\
\L(0.60*-0.15) (b) Tangle replacement \\
\endSetLabels
\AffixLabels{\includegraphics[scale=0.5]{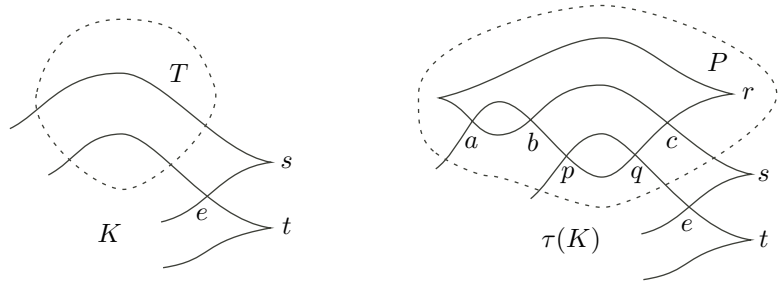}}
}
\vskip 10pt
\label{fig11}
\caption{Legendrian tangle replacement}
\end{figure}

The purpose of this section is to reinterpret the examples in \S\ref{sec3} as
instances of {\it iterated tangle replacement}.  In general one can
modify a Legendrian knot $K$ that is expressed as a tangle sum $S\oplus
T$ by replacing $T$ with another tangle $T'$, giving $S\oplus T'$.
For simplicity we study only the case in which $K$ has a simple front
with upper right corner as shown in Figure~\ref{fig11}(a).  Such a knot will be
called a {\it special Legendrian knot}.   Also assume that  $T$ is the
trivial tangle near the top right cusp $s$, and that $T'$ is the prime
tangle $P$ considered in the last section.  The Legendrian knot resulting
from this particular tangle replacement will be denoted $\tau(K)$.

%%% FIG 11 %%%
%%%%%%%%%

\begin{remark}
\label{rem4.1}
If $S$ is a prime tangle, then the knot $\tau(K)=S\oplus P$ is prime, by
Lickorish's theorem \cite{l}.  In fact {\sl all} the knots in the sequence
$$
 \tau(K),\ \tau^2(K),\ \tau^3(K),\ \ldots
$$
are then prime.  This can be seen by induction since $S\oplus P$ can be
written as $S'\oplus T$ where $S'$ is a prime tangle obtained by adding
two ``clasps" and one ``twist" to $S$ as in \S\ref{sec3}.
\end{remark}

\begin{theorem}
\label{thm4.2}
If $K$ is any special Legendrian knot with nonzero Chekanov number,
then the sequence of Chekanov numbers
$$
\ch(\tau(K)),\ \ch(\tau^2(K)),\ \ch(\tau^3(K)),\ \ldots
$$
grows without bound.
\end{theorem}

For example, applying the theorem to the mirror trefoil ($\bar 3_1$
in the table below) recovers the unboundedness of the Chekanov numbers
of the knots $\bar8_{21},\ldots$  discussed in \S\ref{sec3}, and the fact that
these knots are prime follows from Remark~\ref{rem4.1}.   The mirror
figure eight ($\bar 4_1$ in the table) produces another sequence
$\bar9_{45},\ldots$ of prime knots with increasing Chekanov numbers.
(It is suggestive that $\bar8_{21}$ and $\bar9_{45}$ are the only knots
in the table below with Chekanov number greater than one.)  In fact,
{\sl every} knot with nonzero Chekanov number is isotopic to a special
Legendrian knot, and thence yields such a sequence.

The theorem will be proved by analyzing the effect of a tangle replacement
on the the set of Chekanov polynomials of $K$.  This analysis depends
on a number of factors, the most basic of which is the difference
between the Maslov potentials of the upper and lower strands of $T$.
(See \S\ref{sec2} for the definition of the Maslov potential.)  This potential
difference will be denoted by $d$ below, and will be referred to as the
{\it Maslov number} of $K$.

Note that the Maslov number of $\tau^n(K)$ is equal to 1 for each $n>0$.
Thus for the purposes of proving the theorem, it suffices to consider
the case $d=1$ once it is shown in general that $\ch(\tau(K))\ne0$
(which follows from Lemma~\ref{lem4.4} when $d\ne0$ and
Remark~\ref{rem4.5} when $d=0$).
But for other purposes it may be of interest to study the general case.
Hence we do not at present put any restrictions on $d$.

\begin{figure}\small
\centerline{
%\ShowGrid
\SetLabels
(0.39*0.52) $s$ \\
(0.27*0.23) $e$ \\
(1.02*0.52) $s$ \\
\L(0.00*-0.20) (a) disks for $s_1$ \\
\L(0.65*-0.20) (b) disks for $s_2$ \\
\endSetLabels
\AffixLabels{\includegraphics[scale=0.4]{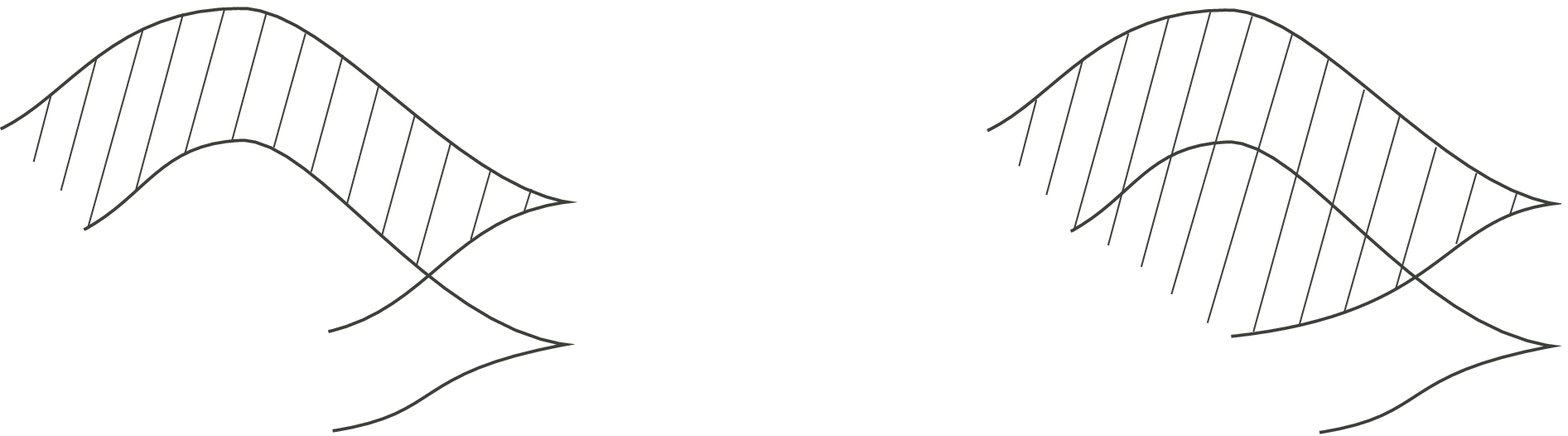}}
}
\vskip 10pt
\caption{Decomposition of $\df s$}
\label{fig12}
\end{figure}

The basic strategy of the proof is to identify a class of augmentations of
$K$ that are ``fertile" enough to cause the Chekanov number to grow under
iterated tangle replacements.  Writing $(\ca,\df)$ for the differential
graded algebra of $K$, consider the elements $s_1$ and $s_2$ in $\ca$
defined by
$$
\df s = 1 + s_1e + s_2.
$$
Thus $s_1e$ is the sum of all the monomials $\df D$ associated with
admissible disks $D$ (for $s$) that end in the crossing $e$, while $s_2$
is the sum of the remaining admissible disk monomials (see
Figure~\ref{fig12}).
Note that if $d\ne1$, then the term $s_1e$ can be ignored for the purposes
of computing Chekanov polynomials, since $e$ and $s_1$ are then both of
nonzero degree (their degrees add to zero and $|e|=1-d\ne 0$).

%%% FIG 12 %%%
%%%%%%%%%

\begin{definition}
\label{def4.3}
An augmentation $\e$ of $K$, viewed as a subset of the degree zero
crossings, is {\it fertile} if either \ (a) \,$d>1$, or \ (b) \,$d=1$,
$s_1=0$ (in $\bz_2$, after substituting $1$ for the crossings in $\e$
and $0$ for the remaining generators of $\ca$), and $\df_0 e$ is a
sum of differentials of some of the other degree zero crossings of
$K$.\foot{Here $\df_0$ is the differential $C_0\to C_{-1}$ in the chain
complex defined by $\e$.}
\end{definition}

Now what happens after the tangle replacement?  There are five new
crossings $a,b,c,p,q$ and one new right cusp $r$ (see Figure~\ref{fig11}).
The remaining crossings and cusps in $\tau(K)$ retain their labels from
$K$.   At this point it is convenient to impose a restriction on $d$,
which will remain in effect for most of what follows, with the exception
of Remark~\ref{rem4.5}.

\begin{assumption}
The {\it Maslov number} $d$ of $K$ is nonzero.
\end{assumption}

\noindent
Thus the crossings $a,b,c$, which are all of degree zero, are the only
new crossings that can potentially be augmented.  The crossings $p$
and $q$ are of nonzero degrees $\pm d$, respectively.

Writing $(\ca^\tau,\df^\tau)$ for the differential graded algebra
of $\tau(K)$ it is clear that $\df^\tau x = \df x$ for {\sl all} the
generators (crossings or right cusps) $x$ of $\ca$ {\sl except} $x=e$,
$s$ and $t$.  For these exceptional cases, an inspection of
Figure~\ref{fig11}
shows that
$$
\begin{aligned}
\df^\tau e &= \df e + q[ba]s_2 + \cdots
\\
\df^\tau s &= 1 + a(\df s -1) + c [ba]s_2 + \cdots
\\
\df^\tau t &= \df t + \cdots
\end{aligned}
$$
where $[ba] = 1 + ba$, and monomials with more than one factor of nonzero
degree are omitted since they do not contribute to the differential for
any augmentation of $\tau(K)$.   Similarly the nonzero differentials of
the new generators of $\ca^\tau$ are
$$
\begin{aligned}
\df^\tau p &= [ba]s_1
\\
\df^\tau c &= as_1q
\\
\df^\tau r &= 1 +  [abc] + \cdots
\end{aligned}
$$
where $[abc]=a+abc+c$.

\begin{lemma}
\label{lem4.4}
If $d\ne0$, then the augmentations of $\tau(K)$ are exactly the sets
$
\e \cup \alpha
$
where $\e$ is an augmentation of $K$ and $\alpha$ is an admissible subset
of $abc$, with the following restriction when $d=1$:
\begin{equation}
\text{$s_1=1$ (in $\bz_2$, after the usual substitution) $\implies \alpha$
is ample}
\tag{$\star$}
\end{equation}
$($Recall from \S\ref{sec3} that the {\rm admissible} subsets of $abc$ are
$abc,ab,a,bc,c$, of which the first two are {\rm ample} and the last
three are {\rm sparse}.$)$
\end{lemma}

\begin{proof}
Start with an augmentation $\e^\tau$ of $\tau(K)$, written as the union
$\e\cup\alpha$ where $\alpha = \e^\tau \cap abc$.  If $\e$ were not
an augmentation of $K$, then there would exist a crossing $x$ in $K$
for which $\df x =1$ (after substitution).  Since $\df^\tau x = 0$
(as $\e^\tau$ {\sl is} an augmentation) it would follow that $x=s$, and so
$
\df^\tau s = 1 + c[ba]s_2 = 0.
$
This would imply $c = [ba] = s_2 = 1$, which in turn would force $s_1=1$,
since $\df s = 1$, and thence contradict $\df^\tau p = 0$.  Thus $\e$
is an augmentation of $K$.

Now the condition $\df^\tau r = 0$ (after substitution) is equivalent to
the admissibility of $\alpha$, and the condition $\df^\tau p = 0$ imposes
the desired restriction $(\star)$ on $\alpha$ since $[ba]$ vanishes if and
only if $\alpha$ is ample.  (Note that this restriction only arises when
$d=1$ since $s_1=0$ otherwise.)   The remaining conditions $\df^\tau s =0$
and $\df^\tau t = 0$ impose no further restrictions on $\alpha$.  This is
clear for the latter, and for the former it follows from the observation
that $\df s = 0$ (since $\e$ is an augmentation) which shows that
$
\df^\tau s = 1 + a + c[ba]s_2 = 0
$
is equivalent to the pair of implications $s_2 = 1 \implies 1+a+c[ba]
= 1 + [abc] = 0$ (which is the admissibility of $\alpha$) and $s_2 =
0 \implies  a = 1$ (which follows from $(\star)$ since $\df s=0$ shows
that $s_1$ and $s_2$ cannot simultaneously vanish).
\end{proof}

\begin{remark}
\label{rem4.5}
If $d=0$, then the same argument shows that the augmentations of $\tau(K)$
{\sl that do not contain $p$ or $q$} are exactly the sets $\e\cup\alpha$,
where $\e$ is an augmentation of $K$ and $\alpha$ is an admissible subset
of $abc$.   However there may also exist augmentations that include one
or both of $p$ and $q$.
\begin{comment}
The terms $qp\df e$, $cps_2$, $qpt_2+q[ba]w$ and $apq$ in $\df^\tau e$,
$\df^\tau s$, $\df^\tau t$ and $\df^\tau r$, respectively, then become
relevant.  Here $t_2$ is defined by $\df t = 1 +et_1+t_2$, and $w$ is
the sum of the $\e$--pure monomials coming from paths joining the top
and bottom points of the part of $K$ shown in Figure~\ref{fig11}(a).
\end{comment}
\end{remark}

\begin{lemma}
\label{lem4.6}
If $d\ne0$, then the reduced Chekanov polynomial of $\tau(K)$
corresponding to an augmentation $\e\cup\alpha$ as above is of the form
$p_{\e\cup\alpha}(t) = p_\e(t) + \Delta(t)$, where
$$
\Delta(t) =
\begin{cases}
1+t^d &\text{for $\alpha$ ample} \\
1-t^{d-1} &\text{for $\alpha$ sparse}
\end{cases}
$$
{\rm if} $\e$ is {\rm fertile}.  If $\e$ is {\rm infertile} then
$\Delta(t) = 1+t^{|d|}$ or $0$ according to whether $s_1=0$ or $1$.
$($Note that $d<0\implies s_1=0$.$)$
 \end{lemma}

\begin{proof}
We must analyze how the chain complexes $(C,\df)$ and $(C^\tau,\df^\tau)$
for $K$ and $\tau(K)$, associated with the augmentations $\e$ and
$\e\cup\alpha$, are related.   As noted above, there are six new
generators in $C^\tau$, namely $a,b,c\in C_0^\tau$, $r\in C_1^\tau$,
$p\in C_d^\tau$ and $q\in C_{-d}^\tau$.   Write $n_k$ and $n_k^\tau$ for
the dimensions of $C_k$ and $C_k^\tau$, and $r_k$ and $r_k^\tau$ for the
ranks of $\df_k$ and $\df_k^\tau$.  It follows from the last remark in
\S\ref{sec2} that the coefficient $\Delta_k$ of $t^k$ in the difference $\Delta(t)
= p_{\e\cup\alpha}(t)-p_\e(t)$ of the reduced Chekanov polynomials is
$$
\Delta_k =
\begin{cases}
\Delta n_{-k}- \Delta r_{-k+1} - \Delta r_{-k-1} &\text{ for $k>0$} \\
1+(-1)^d - \textstyle\sum_{i>0} (-1)^i\Delta_i &\text{ for $k=0$}
\end{cases}
$$
where $\Delta n_k = n_k^\tau - n_k$ and $\Delta r_k = r_k^\tau - r_k$.
Thus it suffices to compute $\Delta n_j$ and $\Delta r_j$ for all $j\le0$.
Evidently $\Delta n_0 = 3$, $\Delta n_{-|d|} = 1$ and $\Delta n_j =
0$ for all other $j<0$.  To compute $\Delta r_j$ it is convenient to
consider the cases $d\ne1$ and $d=1$ separately.

First assume $d\ne 1$.  Then working in $\bz_2$ after the usual
substitution, $\df_0^\tau c=0$, since the degrees of $s_1$ and $q$ are
both nonzero, and so $\Delta r_0 = 0$.  It is clear that $\Delta r_j
= 0$ for all $j<0$, except possibly $j=1-d$ with $d>1$.  (Note that
$\df_{-|d|}^\tau p=0$ if $d<0$ since $s_1=0$.)  But in this case
$\df_{1-d}^\tau e = \df_{1-d} e + q[ba]s_2$ and so
$$
\Delta r_{1-d}=[ba]s_2.
$$
Note that $[ba]=0 \iff \alpha$ is ample, and $s_2=0 \implies s_1=1
\implies \alpha$ is ample (by $(\star)$), whence $\Delta r_{1-d} =
0 \iff \alpha$ is ample.  The formula for $\Delta(t)$ follows readily.

Now assume $d=1$.   Then $\Delta r_j = 0$ for all $j<0$, and so
the formulas above show that $\Delta(t) = 1+t$ or $0$ according to
whether $\Delta r_0 = 0$ or $1$.  To compute $\Delta r_0$, note that
(the matrix for) $\df_0^\tau$ is obtained from  $\df_0$ by adding three
new zero columns, corresponding to the generators $a$, $b$ and $c$, and
then adding a new row, corresponding to $q$.  All the entries in this
last row are zero except possibly the ones in the $e$ and $c$ columns,
denoted $\braket qe$ and $\braket qc$.  In particular (working in $\bz_2$)
$$
\braket qe =[ba]s_2 \quad\text{and}\quad \braket qc =  as_1
$$
and so $\braket qe = 0 \iff \alpha$ is ample (as shown above) and $\braket
qc = 0 \iff s_1=0$ (since $a=0 \implies s_1=0$ by $(\star)$).
Therefore, if $\braket qc=0$ (ie $s_1=0$) then $\Delta r_0 = 0$ if and
only if either $\braket qe=0$ (ie $\alpha$ is ample), or $\braket qe=1$
and the $e$ column in $\df_0$ is independent of the other columns (ie
$\alpha$ is sparse and $\e$ is infertile).  If $\braket qc=1$ (ie
$s_1=1$) then clearly $\Delta r_0 =1$.  Unravelling the definitions now
gives the desired result.
\end{proof}

All the ingredients are now in place to prove the theorem.

\medskip
{\bf Proof of Theorem~\ref{thm4.2}}\qua
Fix an augmentation $\e$ of $K$, and consider all words $w$ in the
letters $a$ (for ample) and $s$ (for sparse) satisfying
\begin{enumerate}
\item $w$ begins with $a$, and
\item $s$ never appears twice in a row in $w$.
\end{enumerate}
Any such $w$ naturally specifies a family of augmentations of $\tau^n(K)$,
where $n$ is the length of $w$, as follows:

Associate the one letter word $a$ to the family $\{\e\cup abc,\ \e\cup
ab\}$ of augmentations of $\tau(K)$ obtained by extending $\e$ by an
ample subset of $abc$.  (Lemma~\ref{lem4.4} and Remark~\ref{rem4.5}
show that these are indeed augmentations.)   Now $\tau(K)$ has Maslov
degree 1, and it is easily verified that the {\sl new} $s_1$ (the one
for $\tau(K)$ with either augmentation in $a$) is zero.  Therefore the
augmentations in $a$ can be extended to augmentations of $\tau^2(K)$
either amply (specified by $aa$) or sparsely (specified by $as$).

Inductively, if $w$ specifies a family of augmentations of $\tau^n(K)$
for some $n>1$, then the word $wa$ specifies the ample extensions of $w$
to augmentations of $\tau^{n+1}(K)$, while $ws$ specifies the sparse
extensions.  To see that condition (2) above is necessary, note that
the $s_1$ for any augmentation $\e'$ in $w$ is zero if and only if $w$
ends in $a$, and apply condition $(\star)$ in Lemma~\ref{lem4.4}.

Also, noting that an $\e'$ in $w$ is fertile if and only if $w$ ends in
$aa$, it follows from Lemma~\ref{lem4.6} that the reduced Chekanov polynomials
for the augmentations associated with $w$ are all equal.  This polynomial
will be denoted by $p_w(t)$.

To complete the proof it suffices for each $m>1$ to produce a collection
of $m$ words of equal length whose polynomials are distinct.  For example,
the polynomials associated with the words $w_k=(a^2)^k(sa)^{m-k+1}$
(for $k=1,\dots,m$) satisfy
$$
p_{w_k}(t)-p_{a^2}(t) = (m+k-2)(1+t)
$$
by direct calculation using Lemma~\ref{lem4.6}.  Thus these $m$ words, all of
length $2m+2$, have distinct polynomials, and the theorem is proved.
\qed

\begin{remark}
\label{rem4.7}
The method used in this proof shows how to compute the full list of
Chekanov polynomials of $\tau^n(K)$ for each $n>0$,  if $K$ is a special
Legendrian knot with nonzero Maslov number $d$ whose Chekanov polynomials
are known.
\end{remark}

The discussion above suggests two directions for further investigation:

\begin{problem}
\label{prob4.8}
(a)\qua Analyze the effect of a special tangle replacement on the Chekanov
polynomials (as in Lemma~\ref{lem4.6}) for knots of Maslov number zero.

(b)\qua Investigate the effects of arbitrary tangle replacements, at
different locations and using prime tangles other than $P$.
\end{problem}

%%% SECTION 5 : TABULATION %%%

\section{Tabulation}
\label{sec5}

The following is a list of maximal\foot{This means `of maximal
Thurston--Bennequin number in their knot type'.} Legendrian fronts for
each knot of nine or fewer crossings.  The knot types are ordered as
in the original tabulation of Alexander and Briggs: $n_k$ is the $k$th
$n$--crossing knot.  Taking the ambient orientation into account, each
chiral knot in the Alexander--Briggs table actually corresponds to a
pair of knot types (up to topological isotopy) that are mirror images of
each other.  These are distinguished by writing $n_k$ for the knot that
appears in Rolfsen's table \cite{r}, and $\bar n_k$ for its mirror image.
The fronts for these ``mirror pairs" occur side by side in the list
below, with $n_k$ appearing first.  The other columns in the table give
the classical invariants $(\tb,r)$ and the reduced Chekanov polynomials
of the knots and their mirrors.

The achiral knots in the table are identified with an asterisk.
For these knots, the mirror pairs $n_k,\bar n_k$ are topologically
isotopic, although they need not be Legendrian isotopic.  For example
$4_1,\bar 4_1$ are Legendrian isotopic by \cite{eh1}, while the pairs
$8_{3},\bar 8_{3}$ and $8_{12},\bar 8_{12}$ are not, since they are
distinguished by their Chekanov polynomials.  The analogous question
for the remaining achiral pairs in the table is left to the reader.

A special feature of the particular fronts shown below is that the
mirror pairs are ``rotationally related".  More precisely, the front for
$\bar n_k$ is  obtained by rotating the front for $n_k$ a quarter turn
clockwise, and then swapping cusps for caps.  Another special feature is
that all the fronts are of minimal crossing number, with the exception
of $9_{42}$ and its mirror (for which we have not been able to find
minimal crossing, rotationally related fronts).

Our indebtedness to the tabulation efforts of Yufa \cite{y} and Branson
\cite{br} (who first investigated the possibility of finding rotationally
related fronts for chiral pairs) is evident.   Most of the fronts drawn
below appear in one or both of their tables.  To verify that all the knots
shown are maximal (with the possible exception of $\bar 9_{42}$) we appeal
to Ng's computation of the maximal Thurston--Bennequin numbers for
knots up to nine crossings \cite{n1}, which builds on Yufa's tabulation.
The classical invariants and Chekanov polynomials for the listed knots
were computed using a Mathematica program written by the first author.

%%%%%%%%%%%
%\begin{comment}
%%%%%%%%%%%

\begin{table*}[!ht]
\caption{\scshape Invariants of Legendrian Knots}
\vskip .2in
\tabletop
\up{0_1^*}	&\up{(-1,0)}	&\up{0}		&\fig{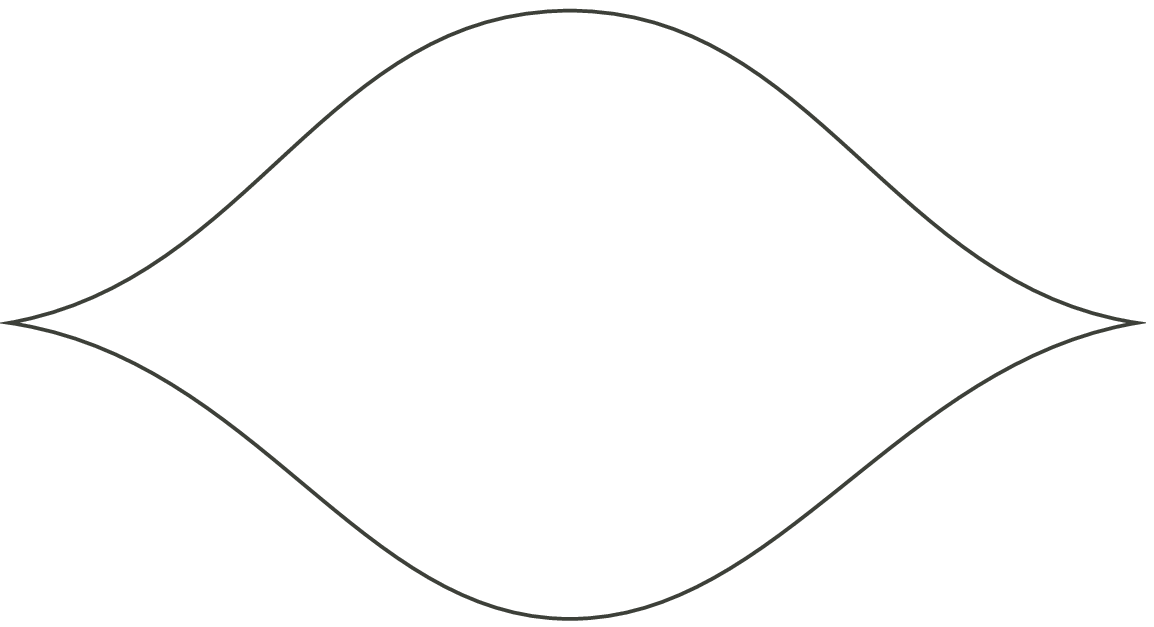}{.14}	&\up{0}
&\up{(-1,0)}	\\ [1ex]
\up{3_1}		&\up{(-6,1)}	&
&\fig{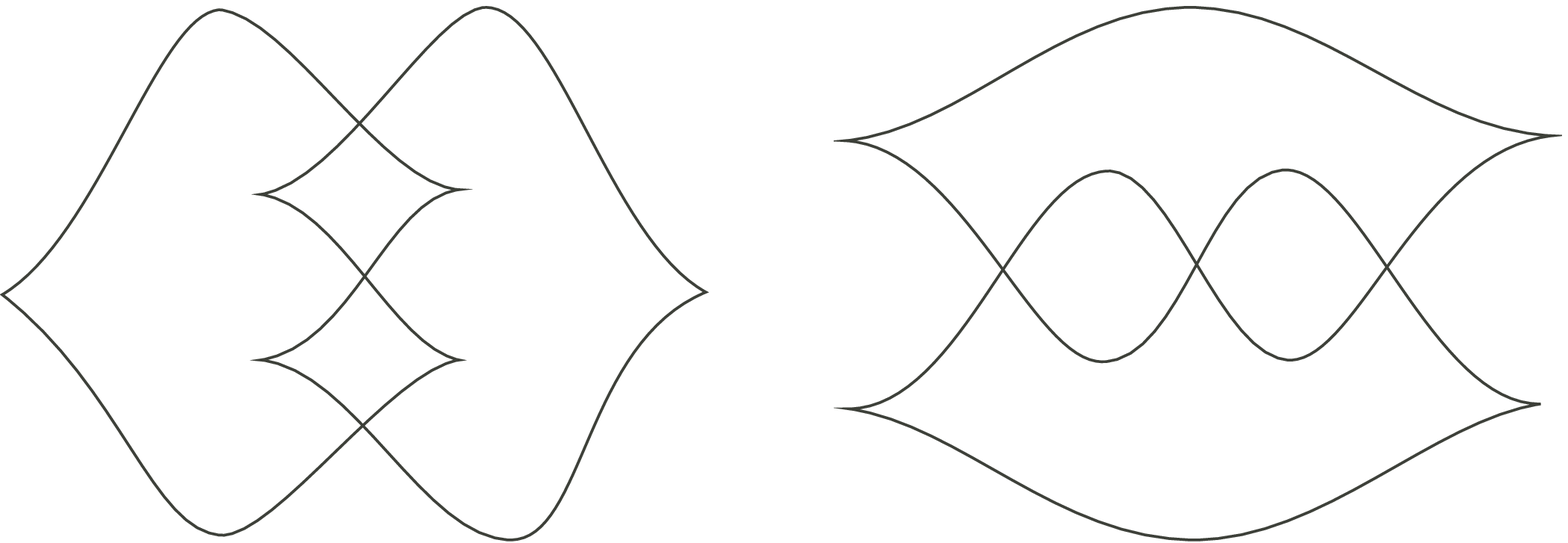}{.14}	&\up{1}		&\up{(1,0)}	\\ [1ex]
\up{4_1^*}	&\up{(-3,0)}	&\up{t}		&\fig{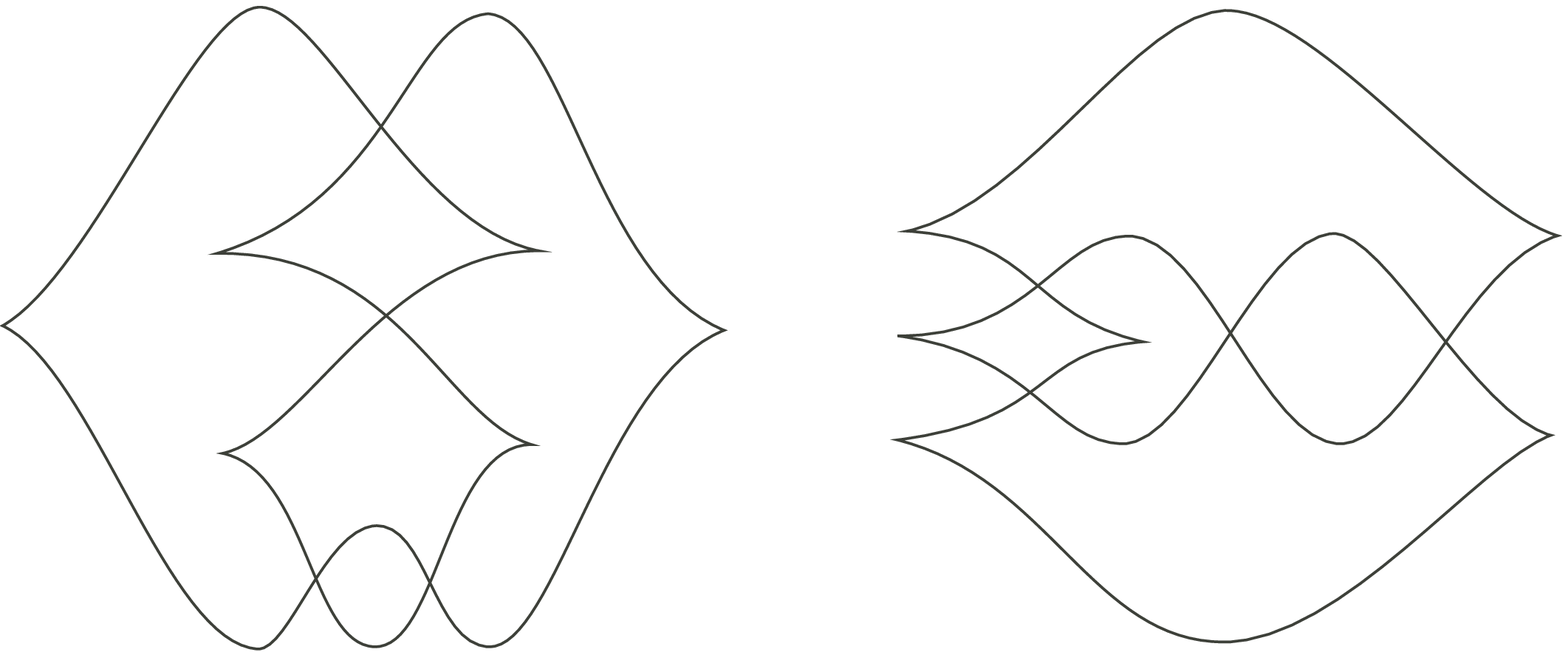}{.14}	&\up{t}
&\up{(-3,0)}	\\ [1ex]
\up{5_1}		&\up{(-10,3)}	&
&\fig{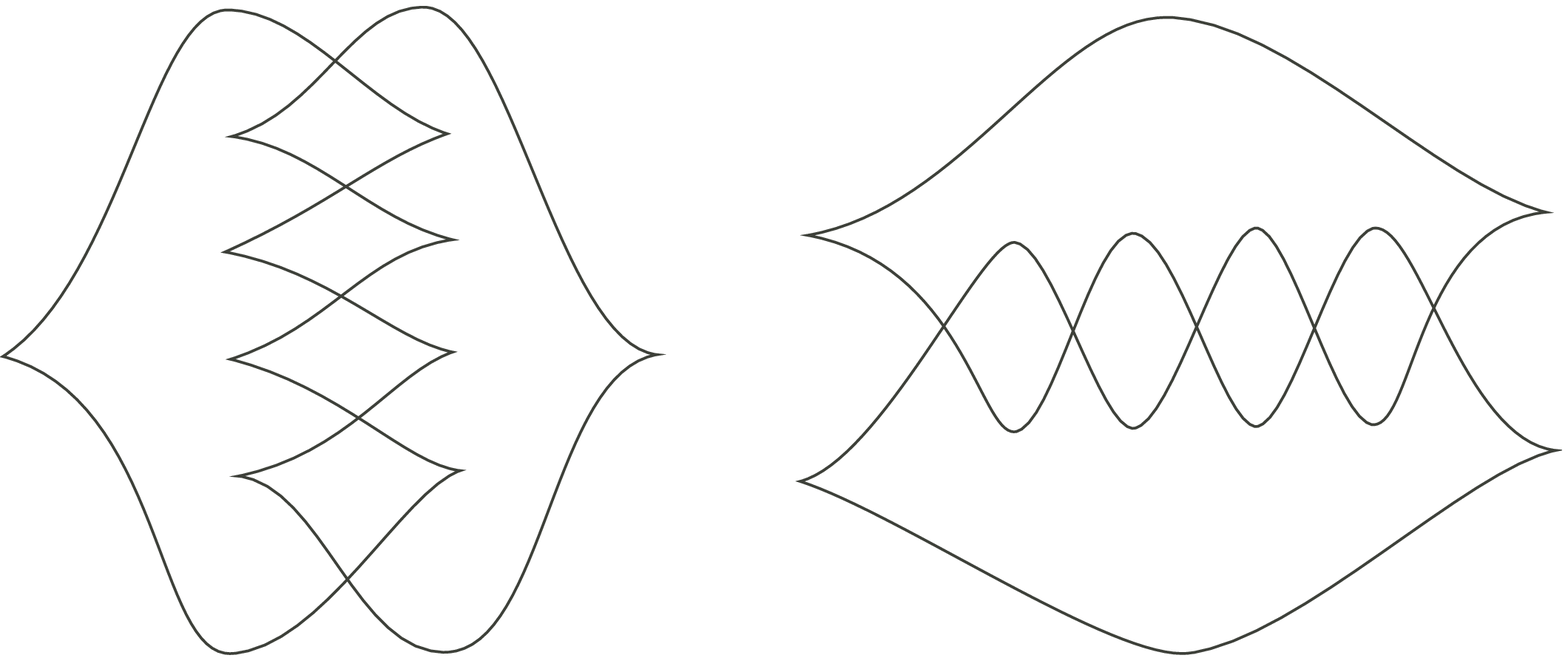}{.14}	&\up2		&\up{(3,0)}	\\ [1ex]
\up{5_2}		&\up{(-8,1)}	&
&\fig{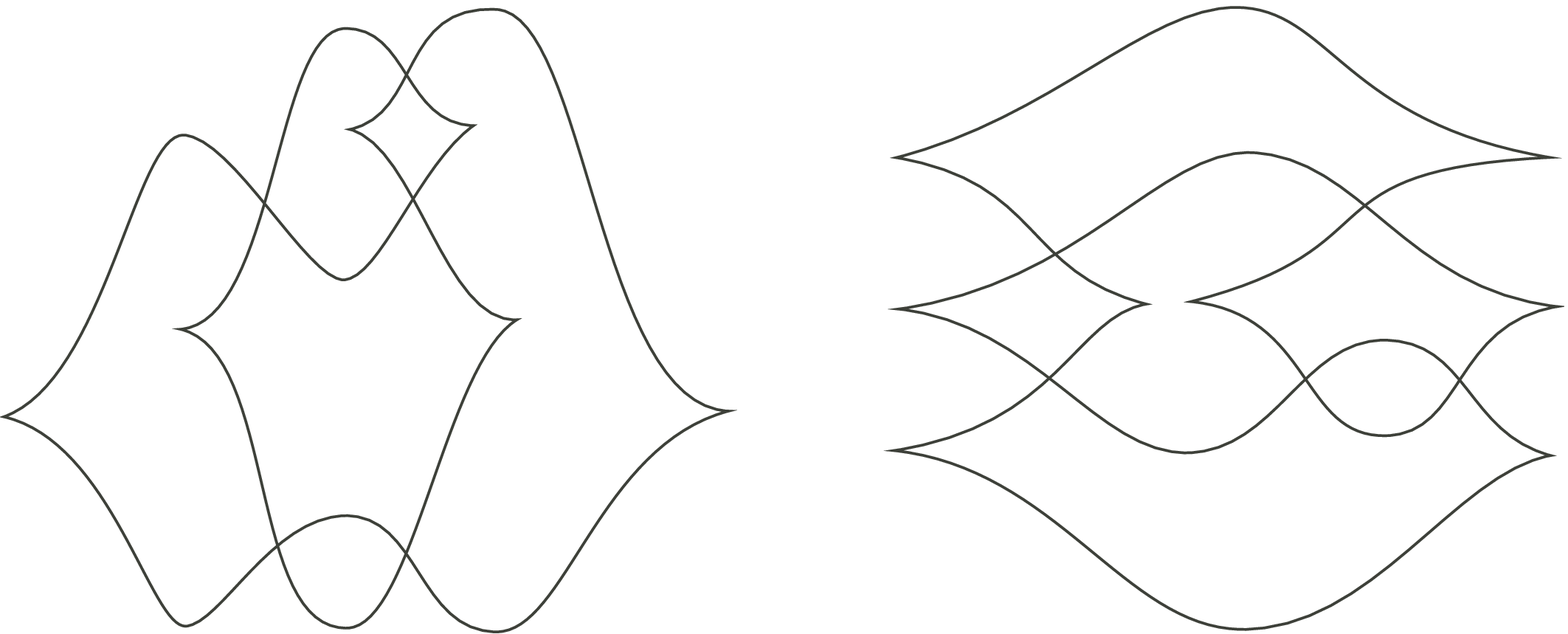}{.14}	&\up1		&\up{(1,0)}	\\ [1ex]
\up{6_1}		&\up{(-5,0)}	&\up{2t}
&\fig{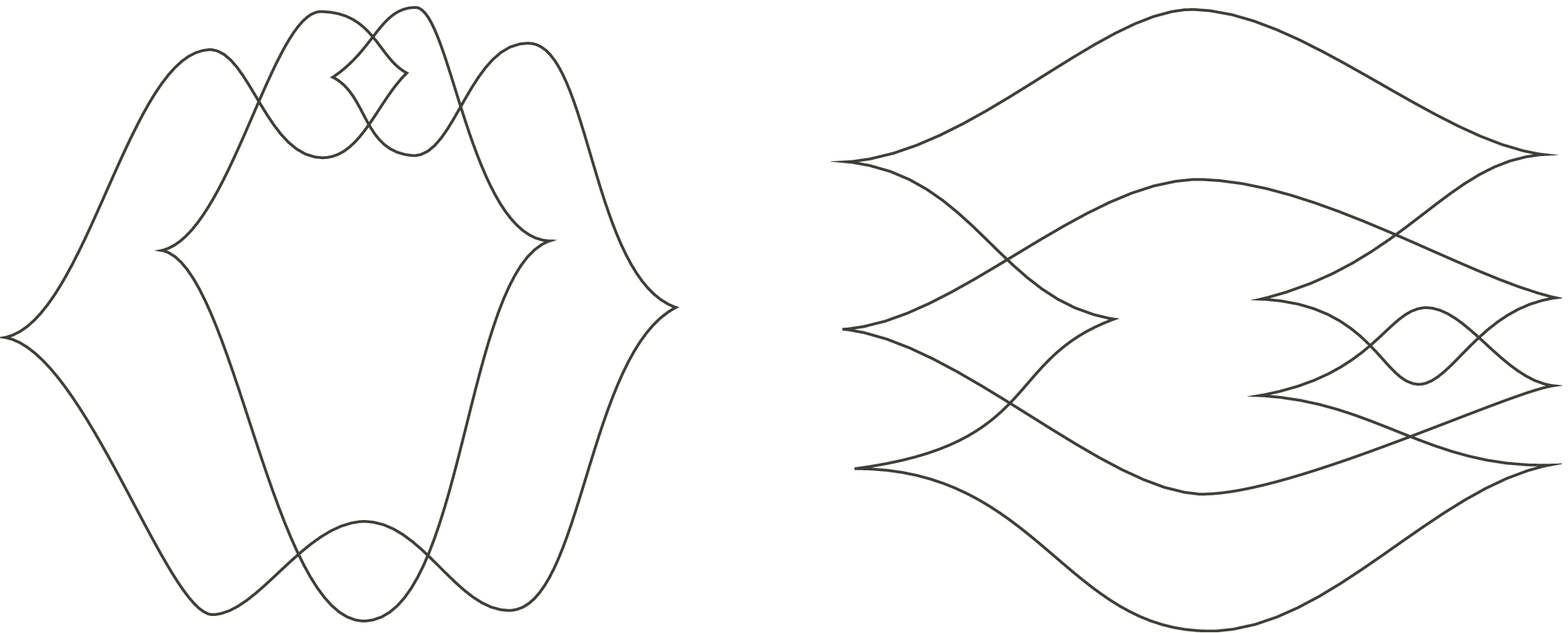}{.14}	&\up{t}		&\up{(-3,0)}	\\ [1ex]
\up{6_2}		&\up{(-7,0)}	&
&\fig{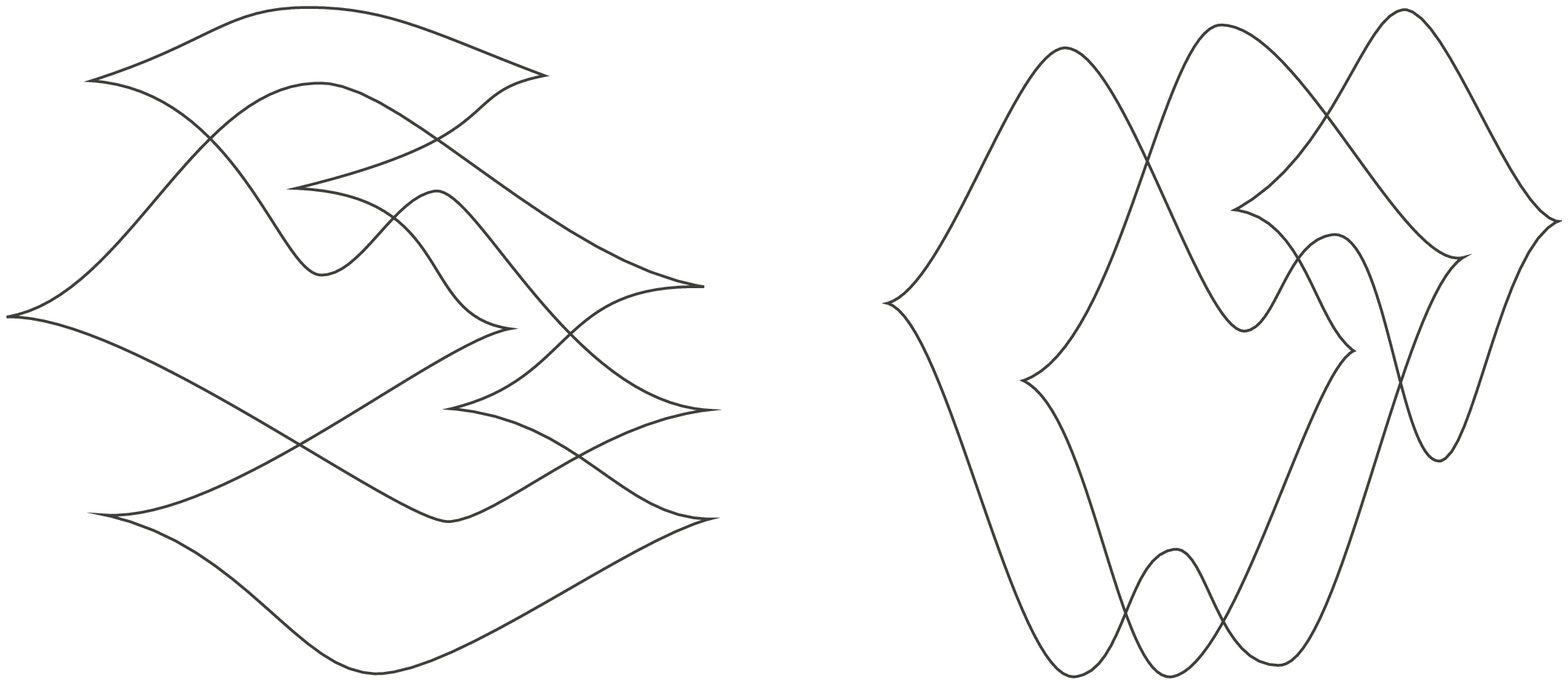}{.14}	&\up{1+t}		&\up{(-1,0)}	\\ [1ex]
\up{6_3^*}	&\up{(-4,1)}	&			&\fig{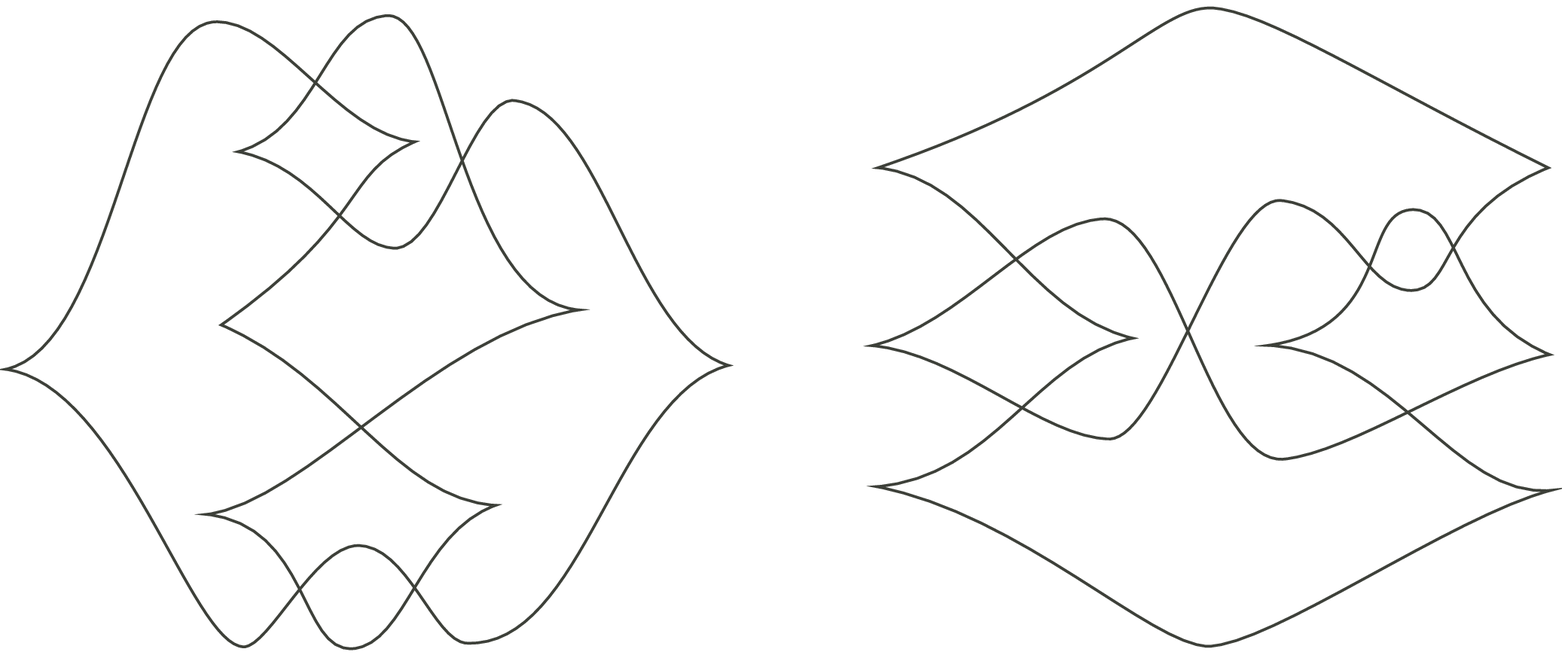}{.14}
&			&\up{(-4,1)}	\\ [1ex]
\up{7_1}		&\up{(-14,5)}	&
&\fig{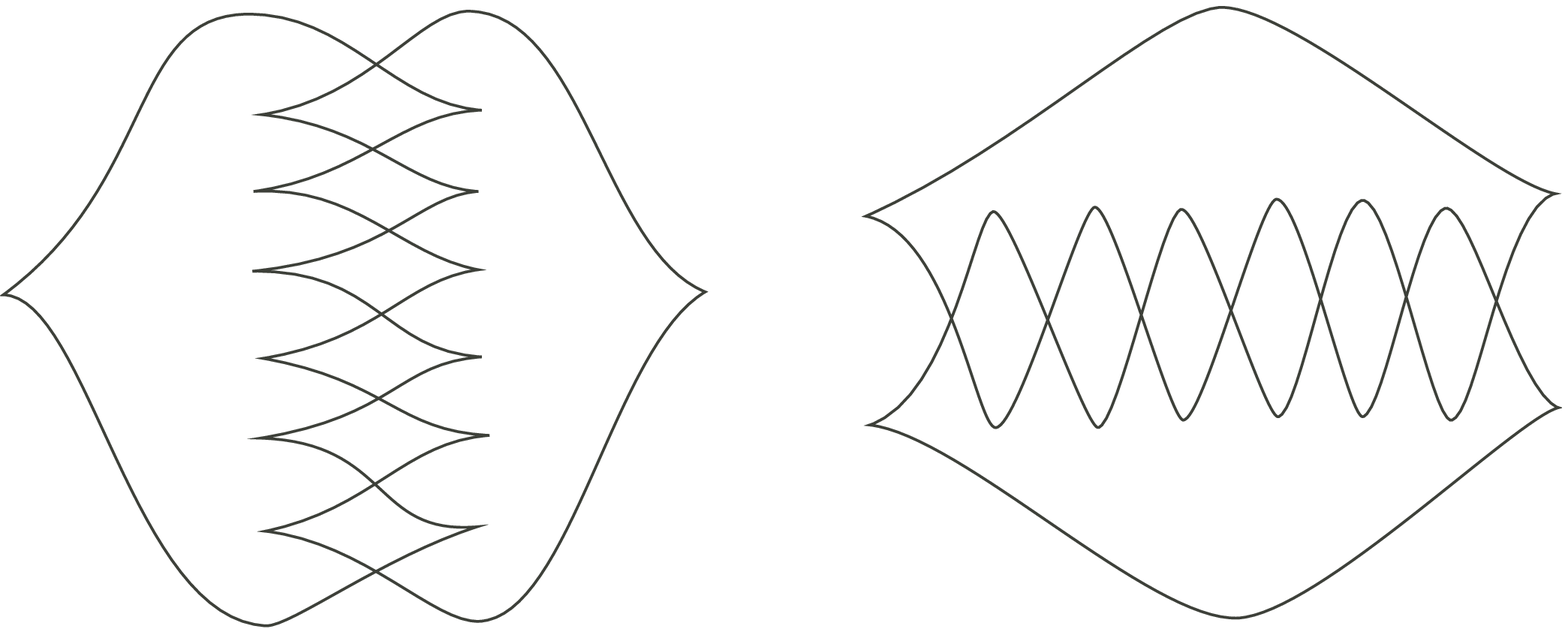}{.14}	&\up{3}		&\up{(5,0)}	\\ [1ex]
\hline
\end{tabular}
\end{table*}

\clearpage

\begin{table}[!ht]
\tabletop
\up{7_2}		&\up{(-10,1)}	&
&\fig{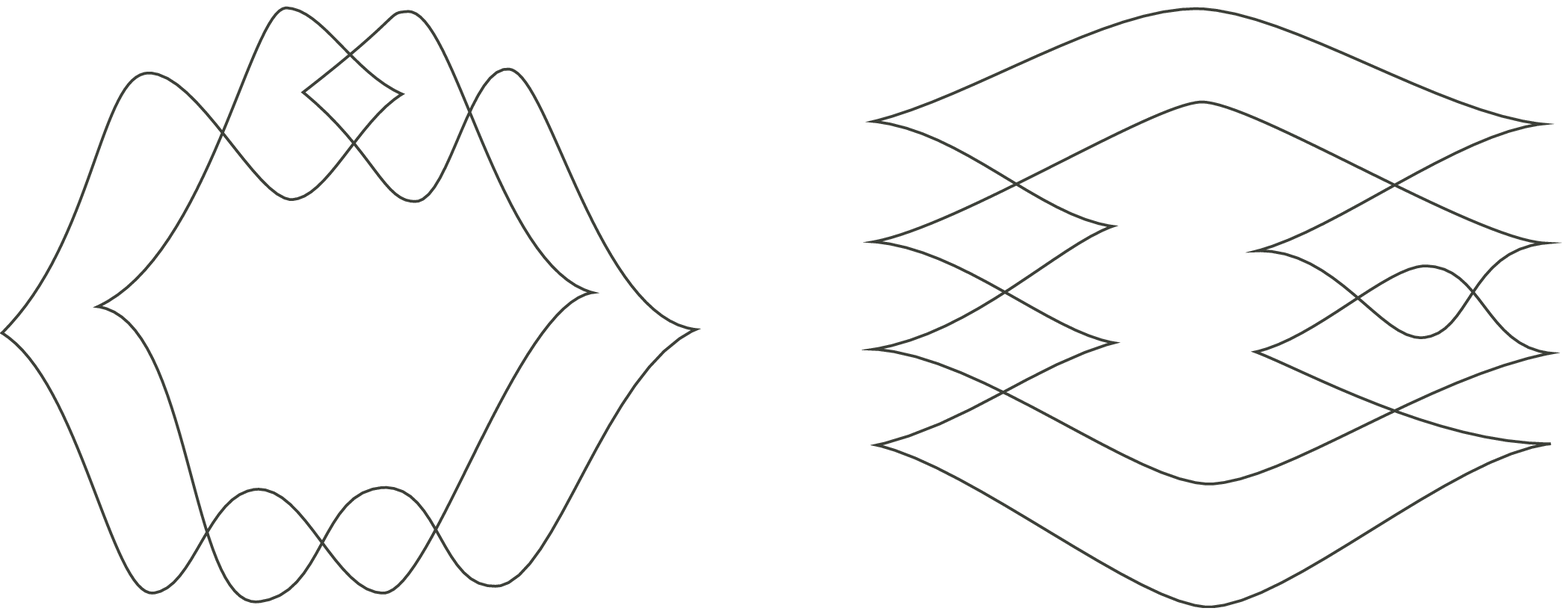}{.14}	&\up{1}		&\up{(1,0)}	\\ [1ex]
\up{7_3}		&\up{(3,0)}	&\up{2t^2}
&\fig{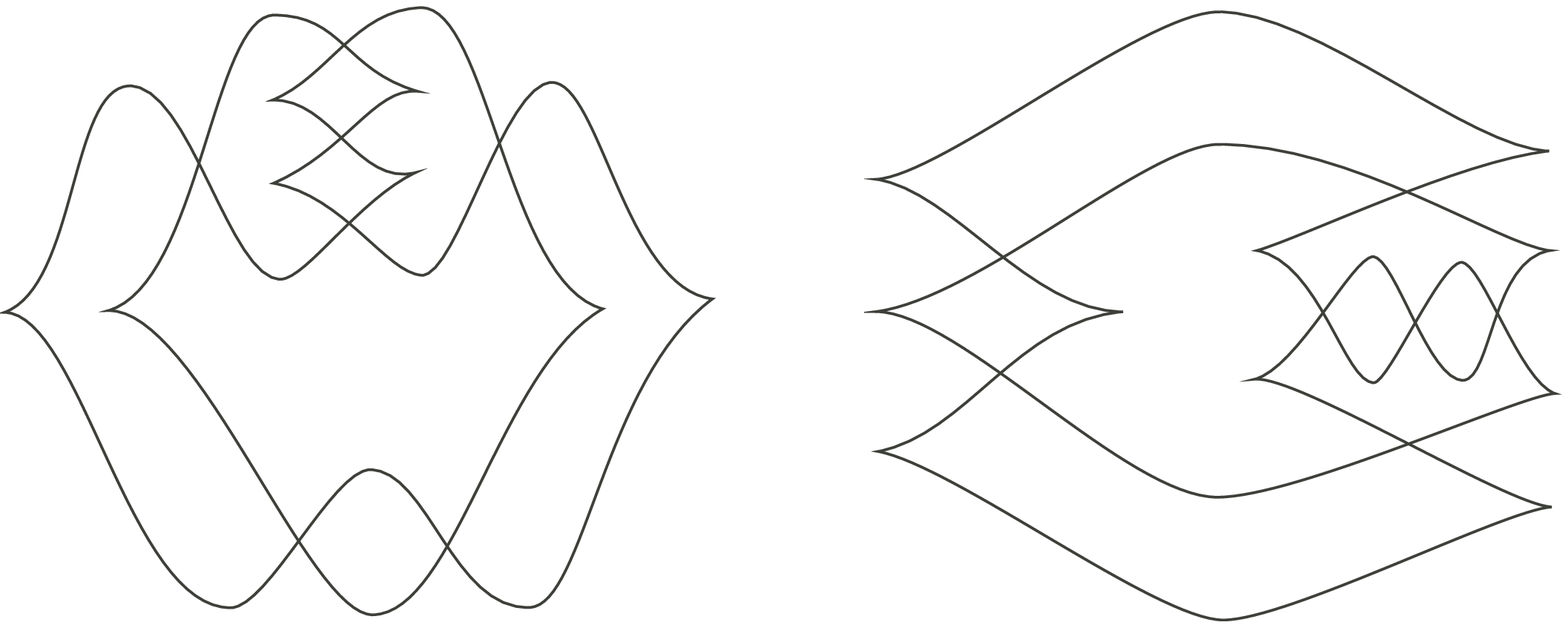}{.14}	&			&\up{(-12,1)}	\\ [1ex]
\up{7_4}		&\up{(1,0)}	&\up{1}		&\fig{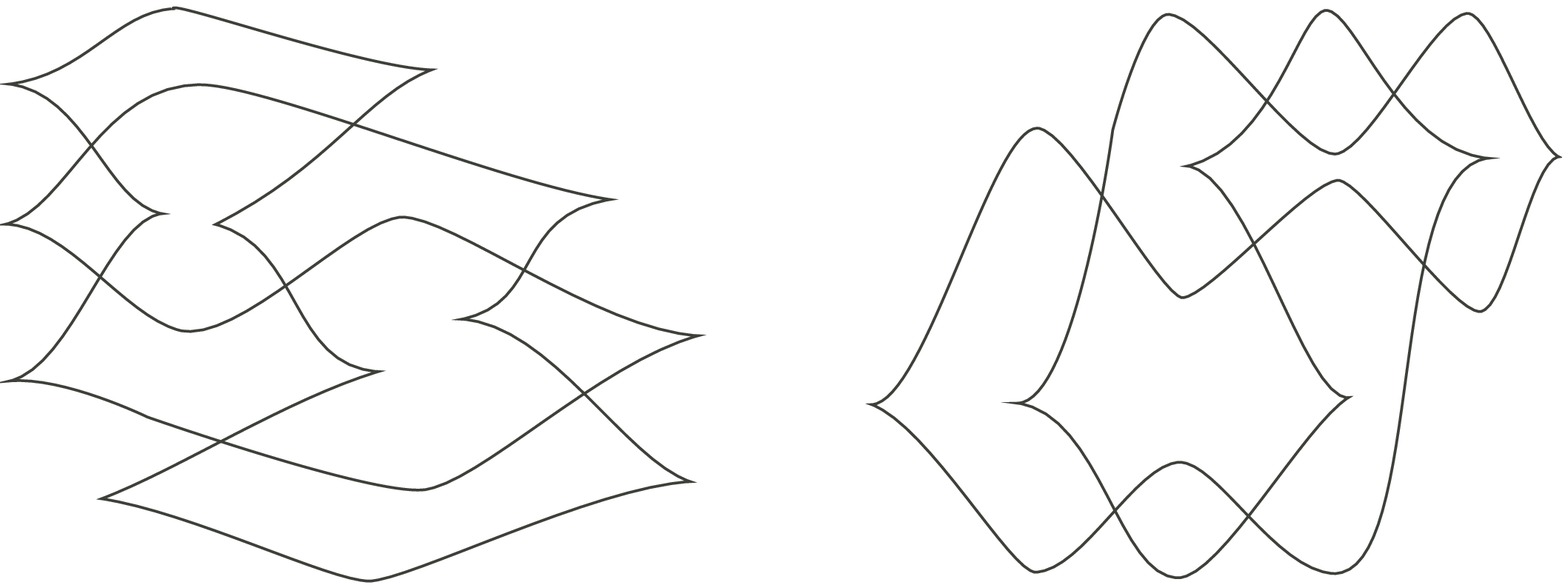}{.14}
&			&\up{(-10,1)}	\\ [1ex]
\up{7_5}		&\up{(-12,1)}	&
&\fig{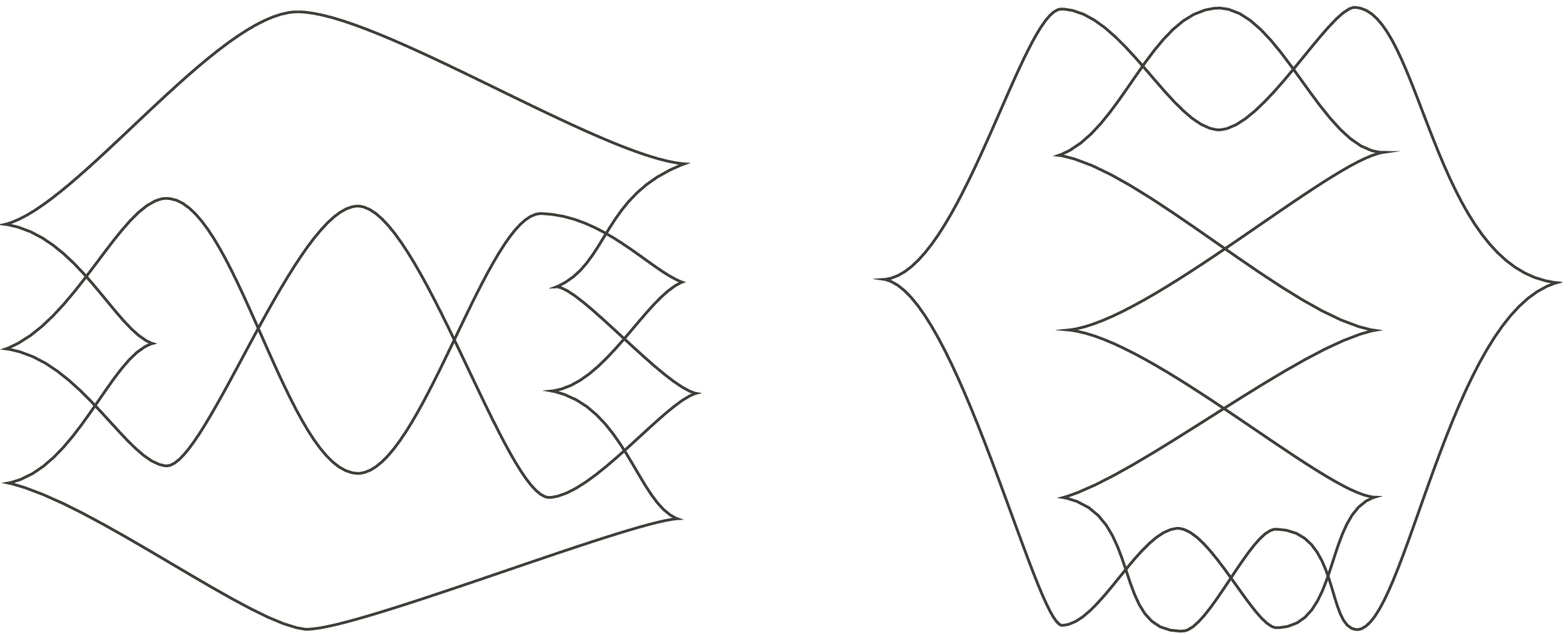}{.14}	&\up{1+t^2}	&\up{(3,0)}	\\ [1ex]
\up{7_6}		&\up{(-8,1)}	&
&\fig{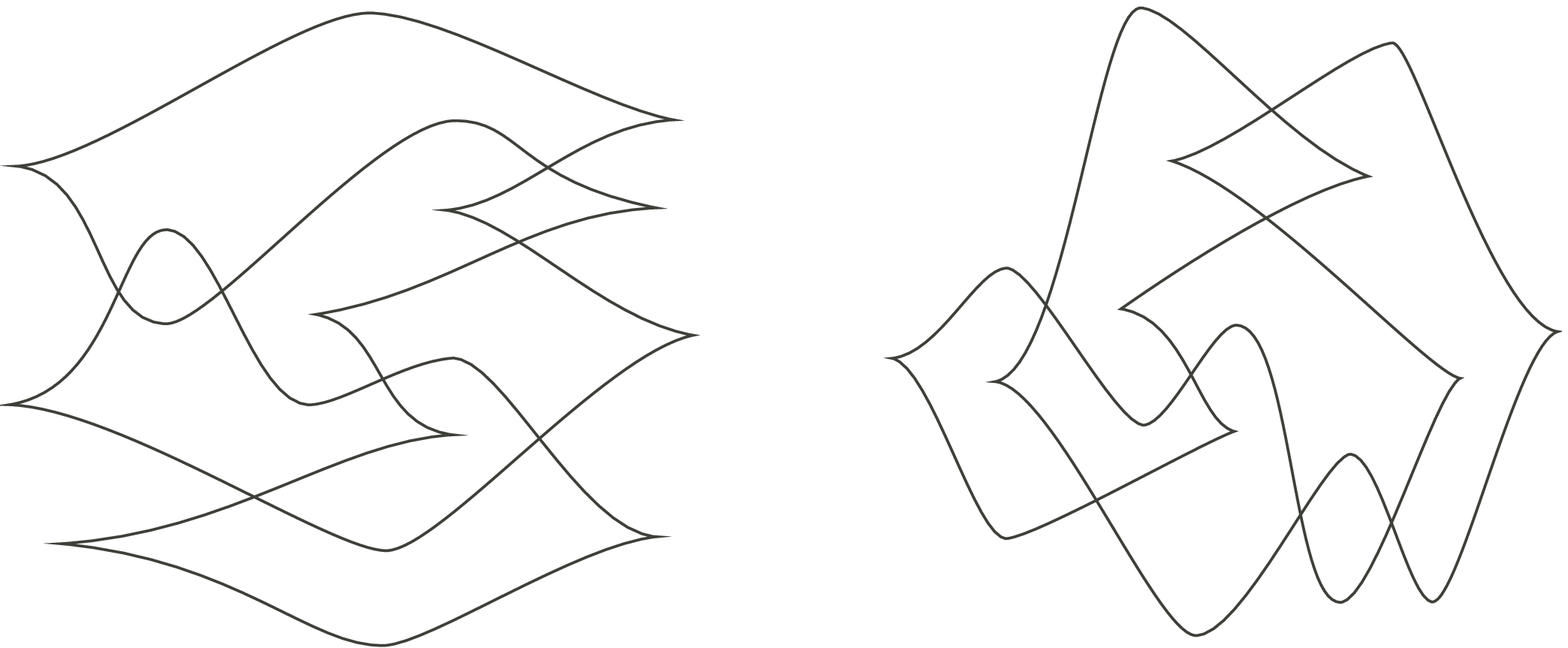}{.14}	&\up{t+t^2}	&\up{(-1,0)}	\\ [1ex]
\up{7_7}		&\up{(-4,1)}	&
&\fig{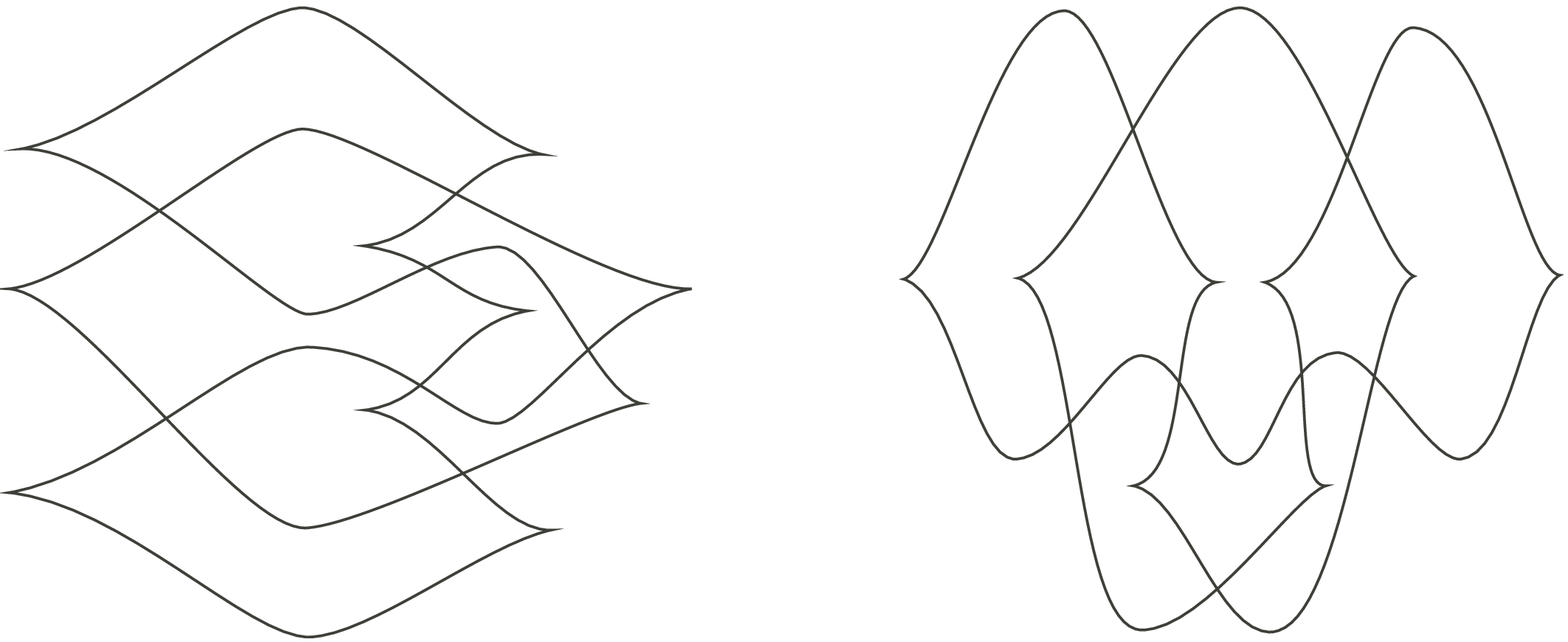}{.14}	&\up{2t}		&\up{(-5,0)}	\\ [1ex]
\up{8_1}		&\up{(-7,0)}	&\up{3t}
&\fig{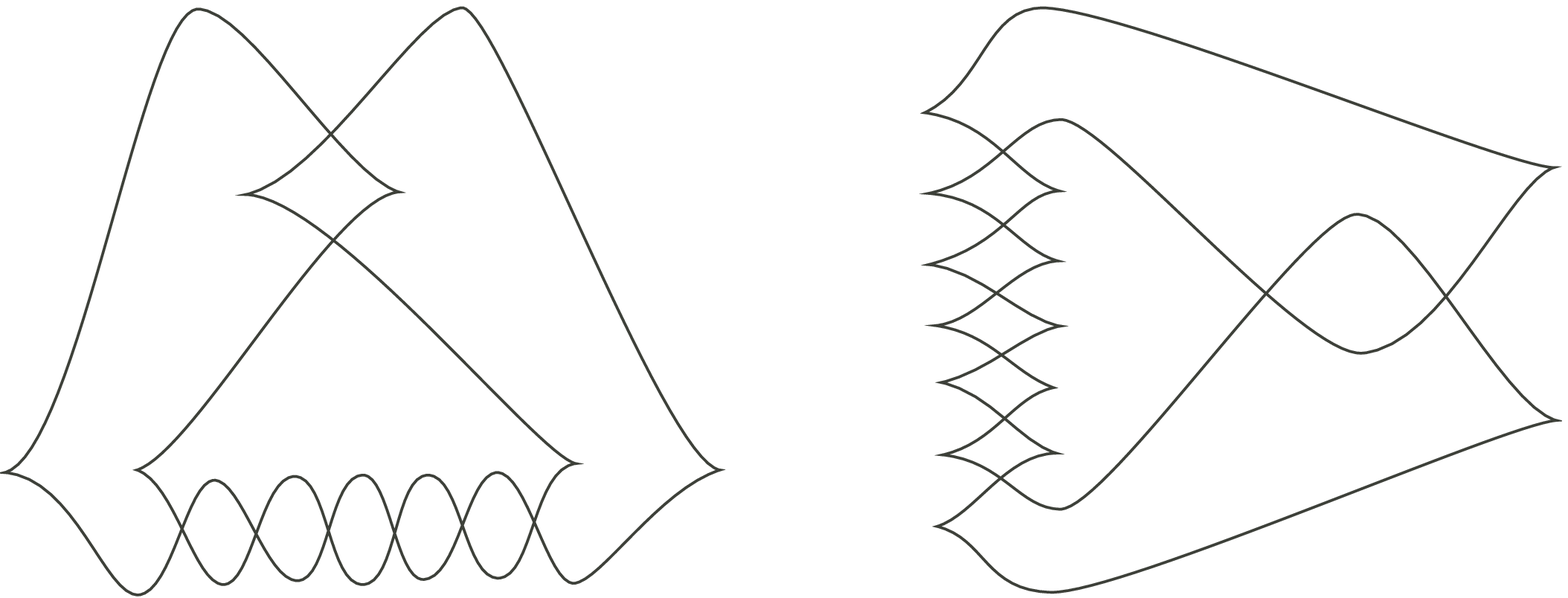}{.14}	&\up{t^5}		&\up{(-3,0)}	\\ [1ex]
\up{8_2}		&\up{(-11,4)}	&
&\fig{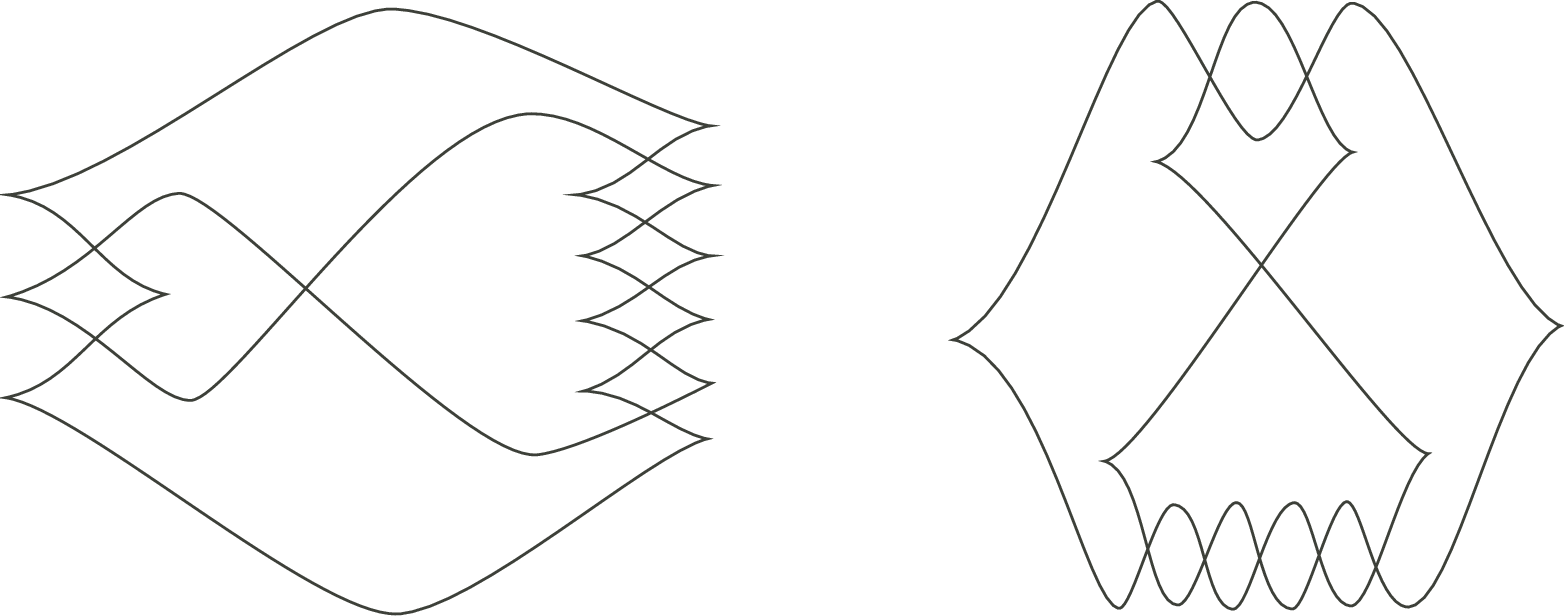}{.14}	&\up{2+t}		&\up{(1,0)}	\\ [1ex]
\up{8_3^*}	&\up{(-5,0)}	&\up{2t^3}		&\fig{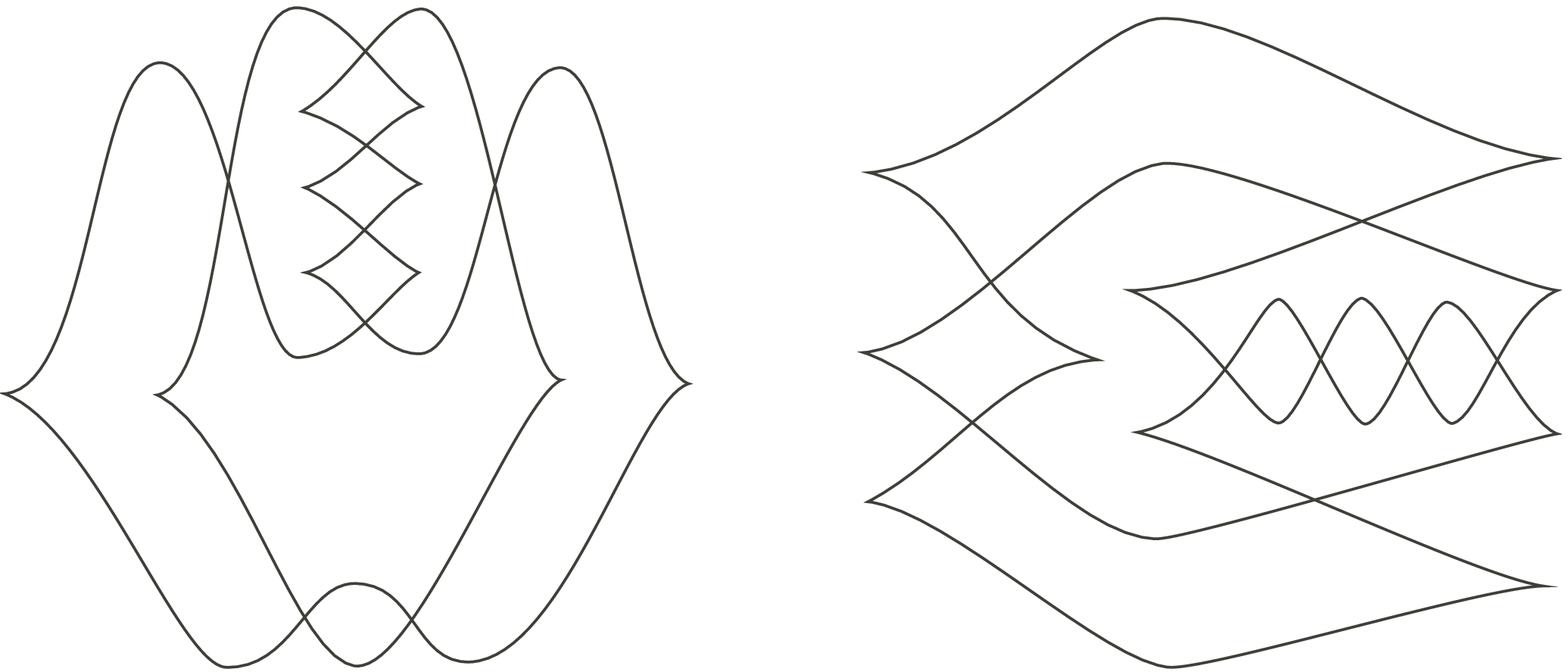}{.14}
&\up{2t}		&\up{(-5,0)}	\\ [1ex]
\up{8_4}		&\up{(-7,2)}	&
&\fig{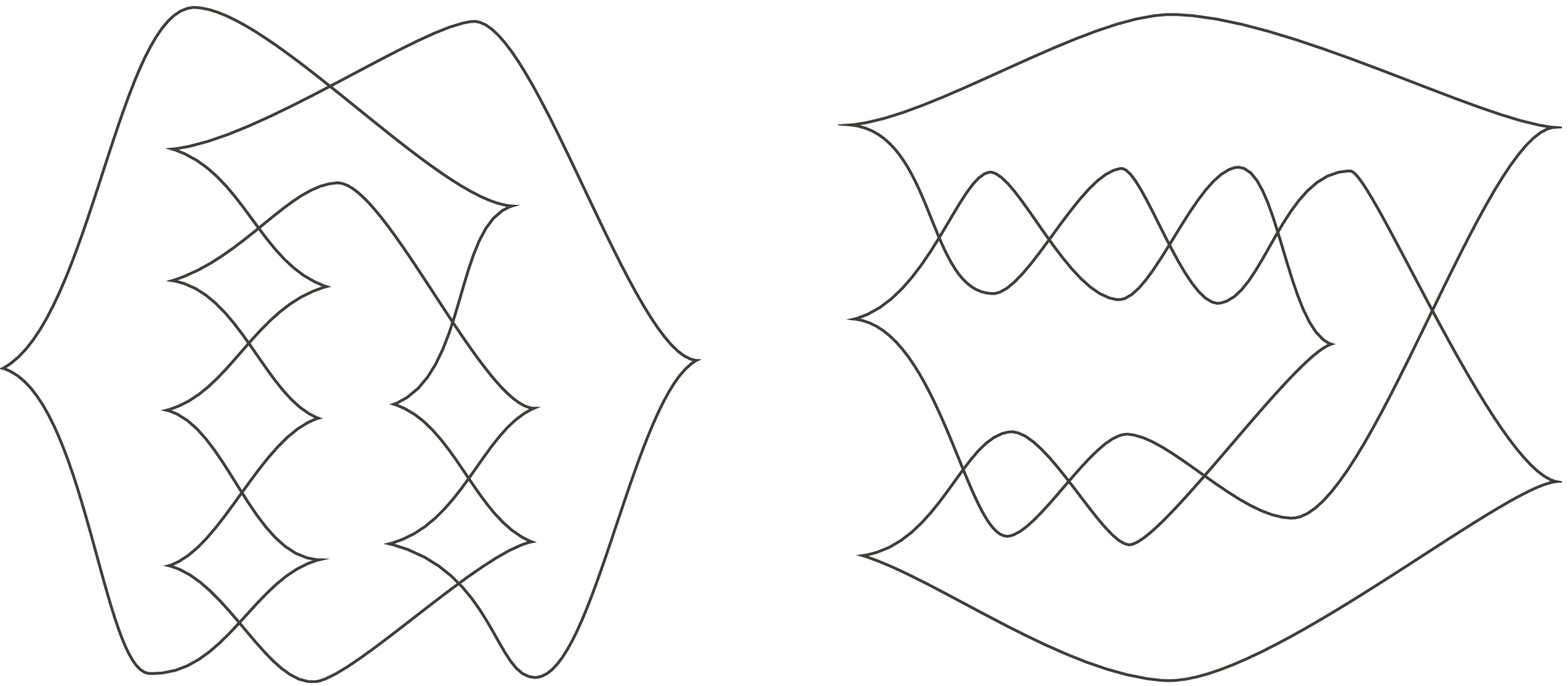}{.14}	&\up{1+2t}	&\up{(-3,0)}	\\ [1ex]
\up{8_5}		&\up{(1,0)}	&\up{2+t}
&\fig{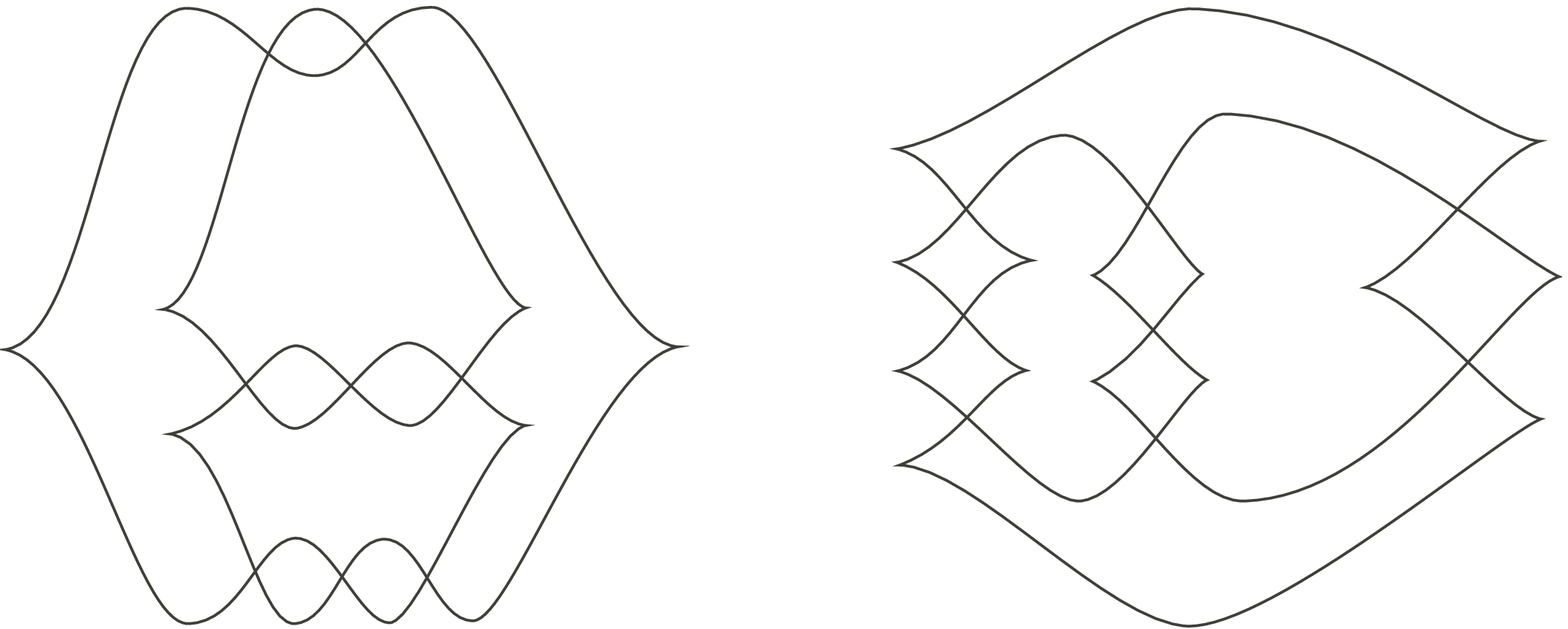}{.14}	&			&\up{(-11,0)}	\\ [1ex]
\hline
\end{tabular}
\end{table}

\clearpage

\begin{table}[!ht]
\tabletop
\up{8_6}		&\up{(-9,2)}	&
&\fig{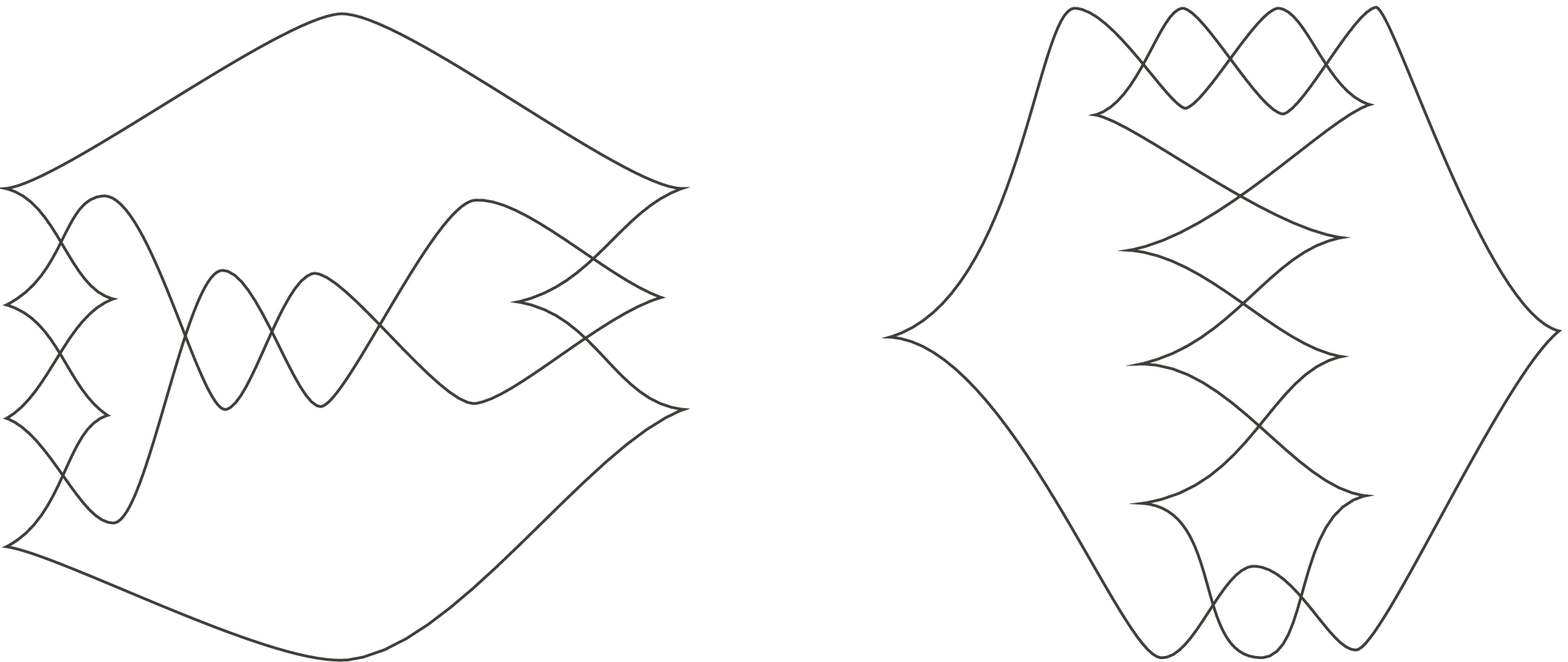}{.14}	&\up{1+t^3}	&\up{(-1,0)}	\\ [1ex]
\up{8_7}		&\up{(-2,1)}	&
&\fig{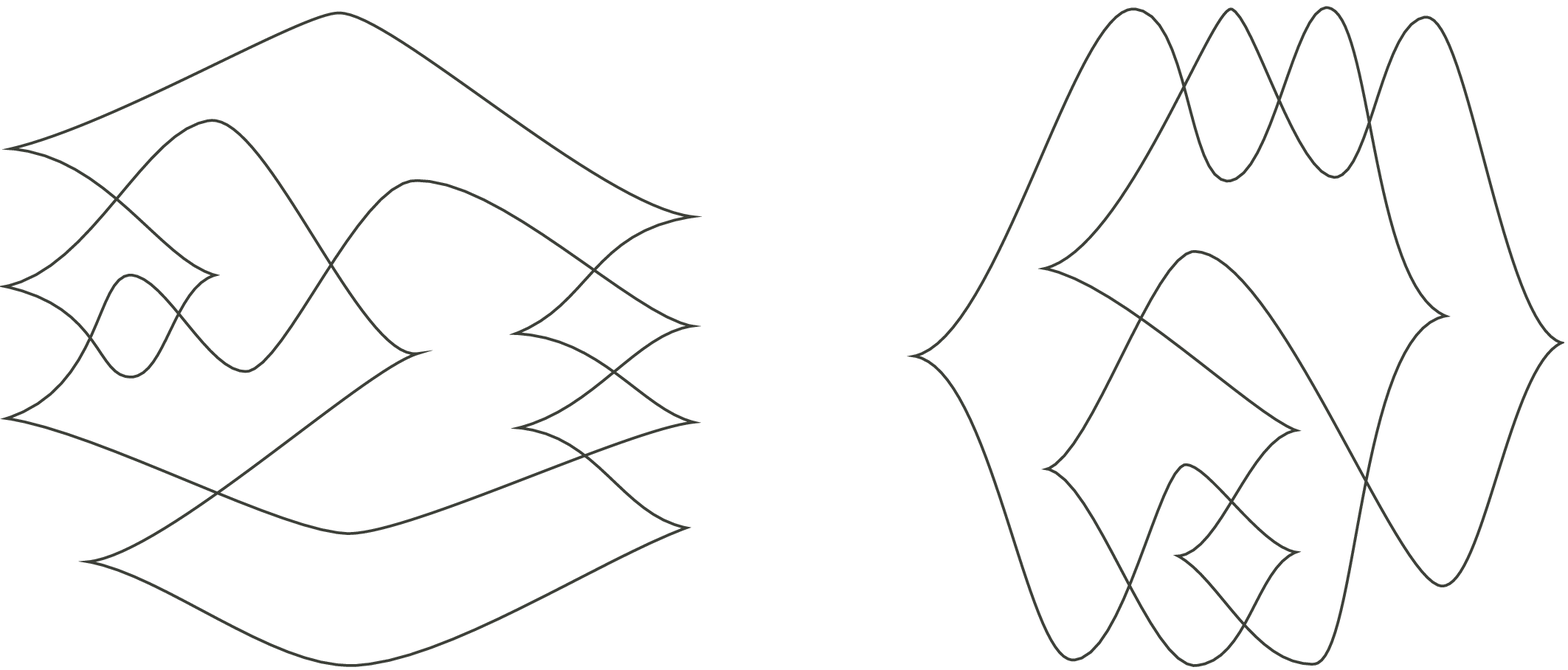}{.14}	&			&\up{(-8,1)}	\\ [1ex]
\up{8_8}		&\up{(-4,1)}	&
&\fig{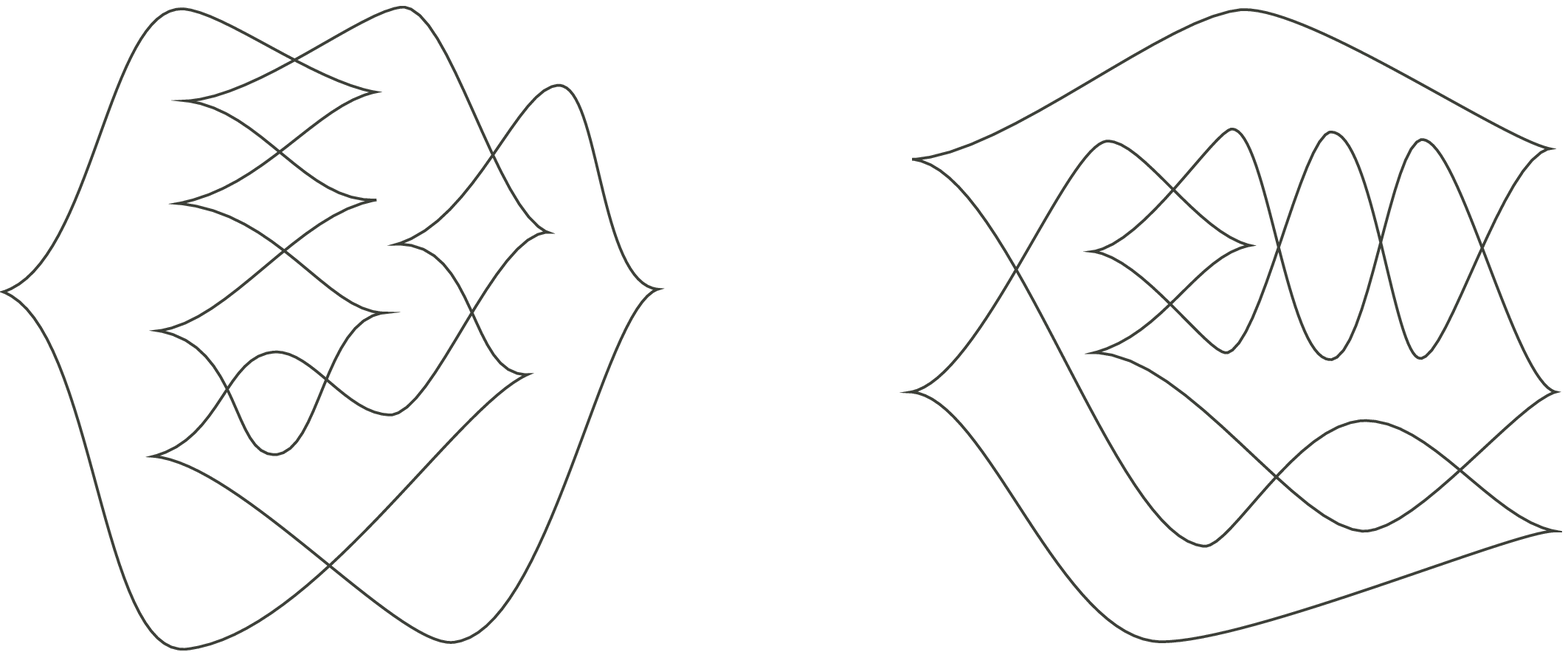}{.14}	&			&\up{(-6,1)}	\\ [1ex]
\up{8_9^*}	&\up{(-5,2)}	&			&\fig{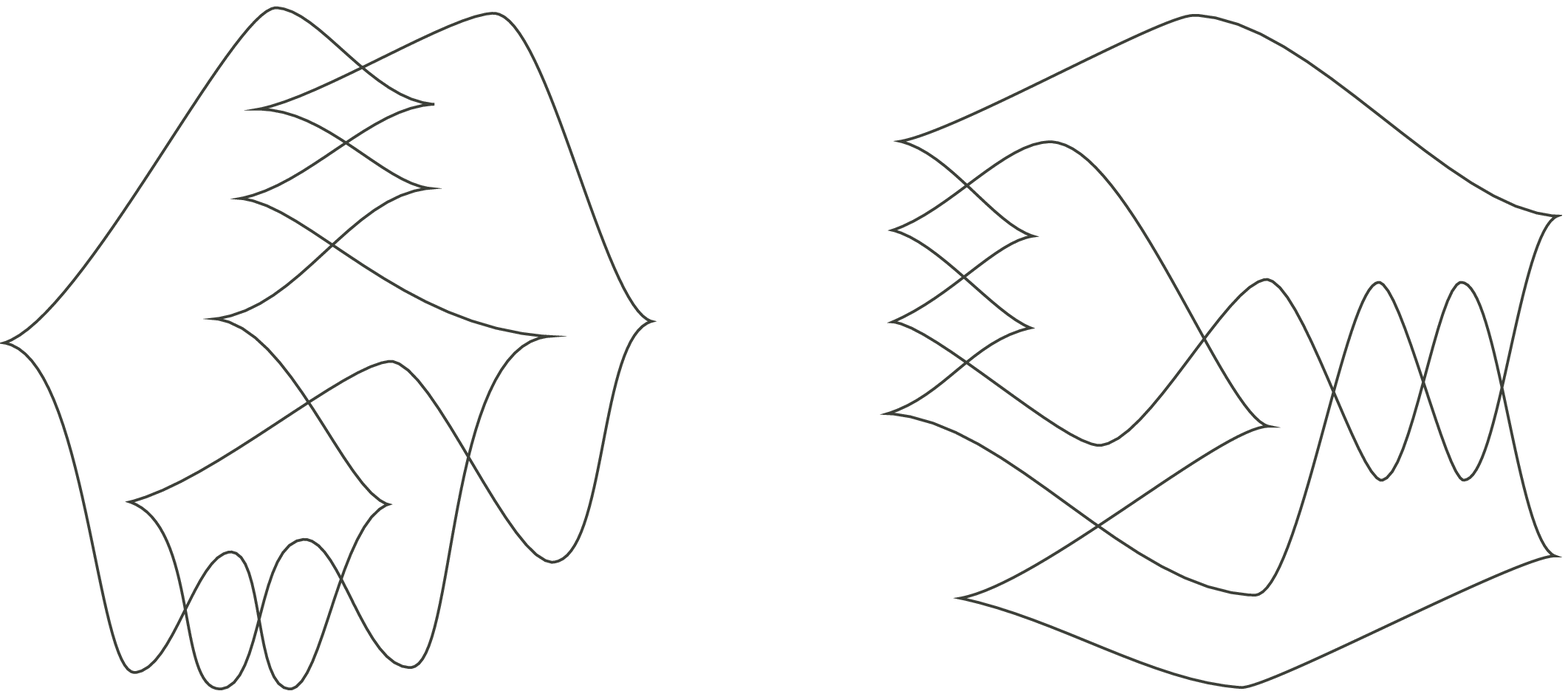}{.14}
&			&\up{(-5,2)}	\\ [1ex]
\up{8_{10}}	&\up{(-2,1)}	&			&\fig{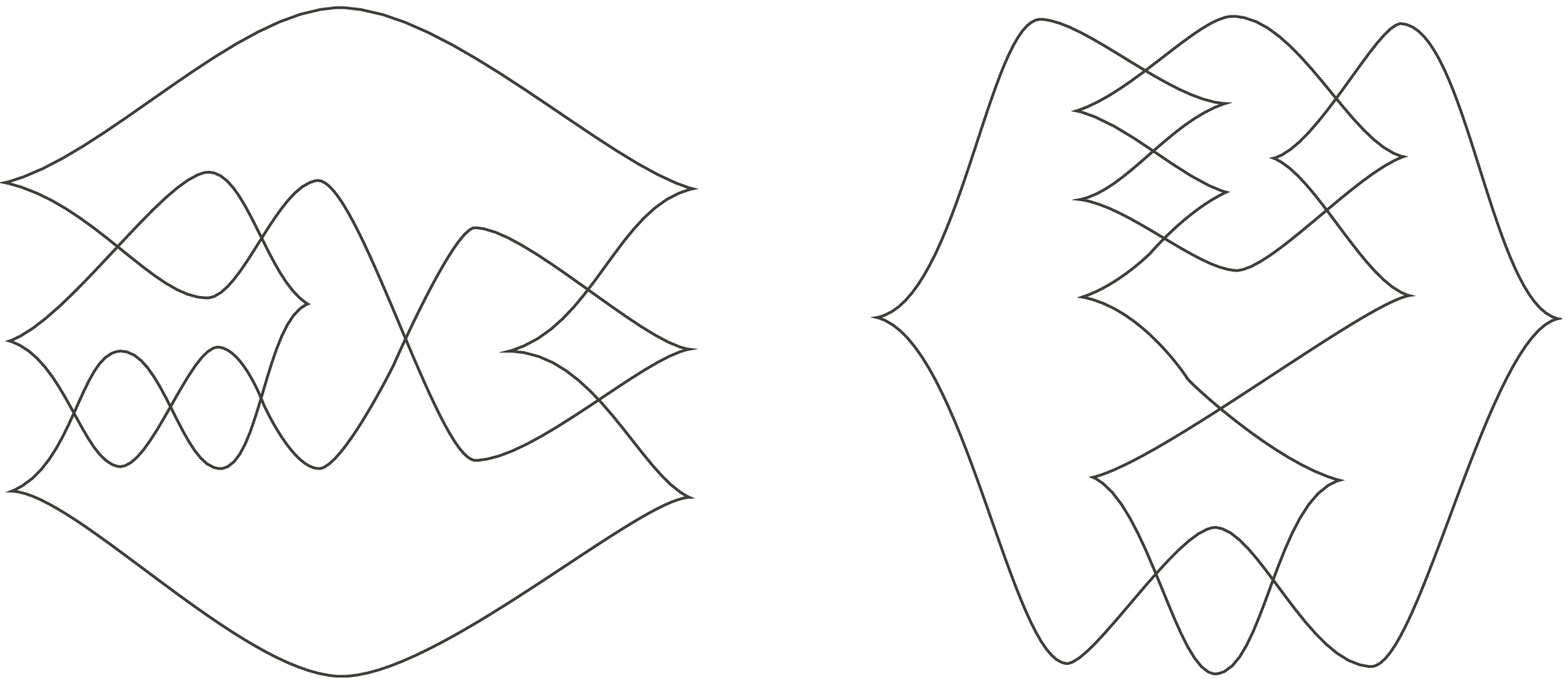}{.14}
&			&\up{(-8,1)}	\\ [1ex]
\up{8_{11}}	&\up{(-9,0)}	&			&\fig{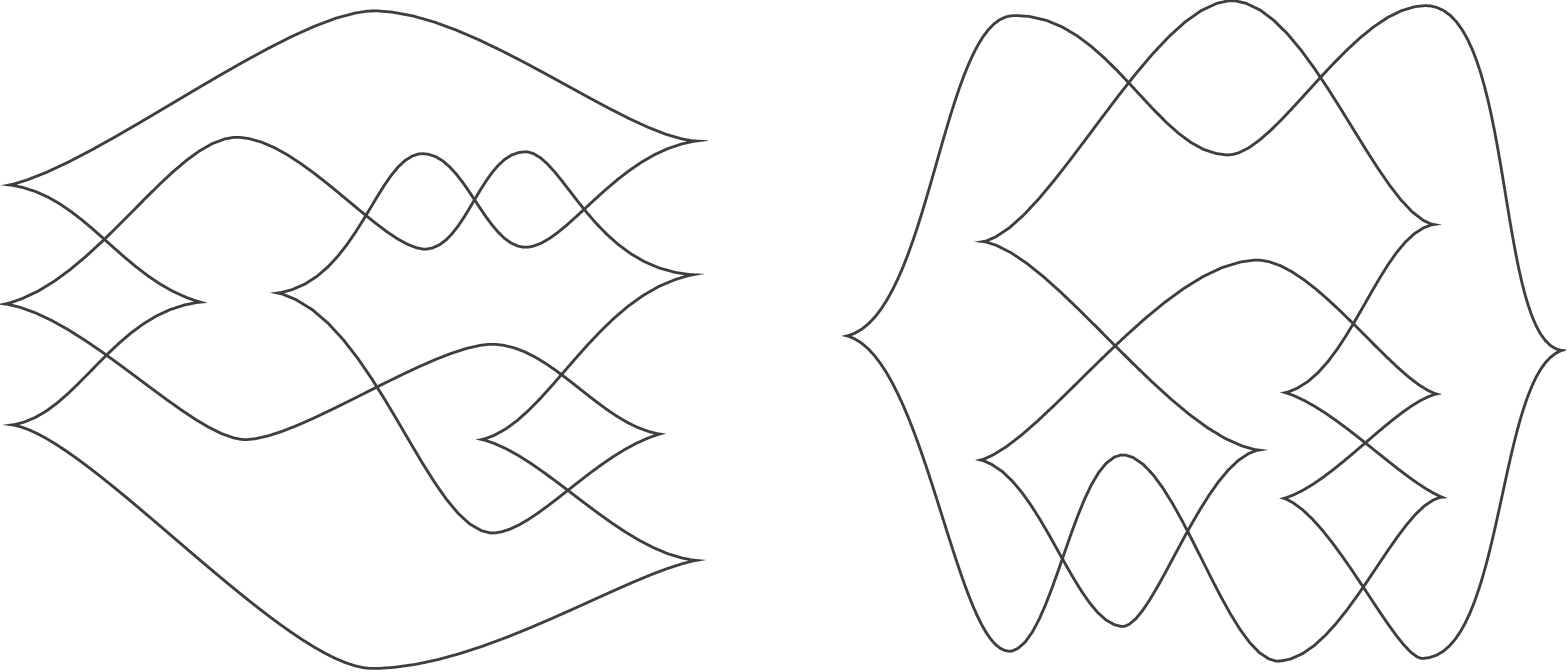}{.14}
&\up{t+t^2}		&\up{(-1,0)}	\\ [1ex]
\up{8_{12}^*}	&\up{(-5,0)}	&\up{t+t^3}	&\fig{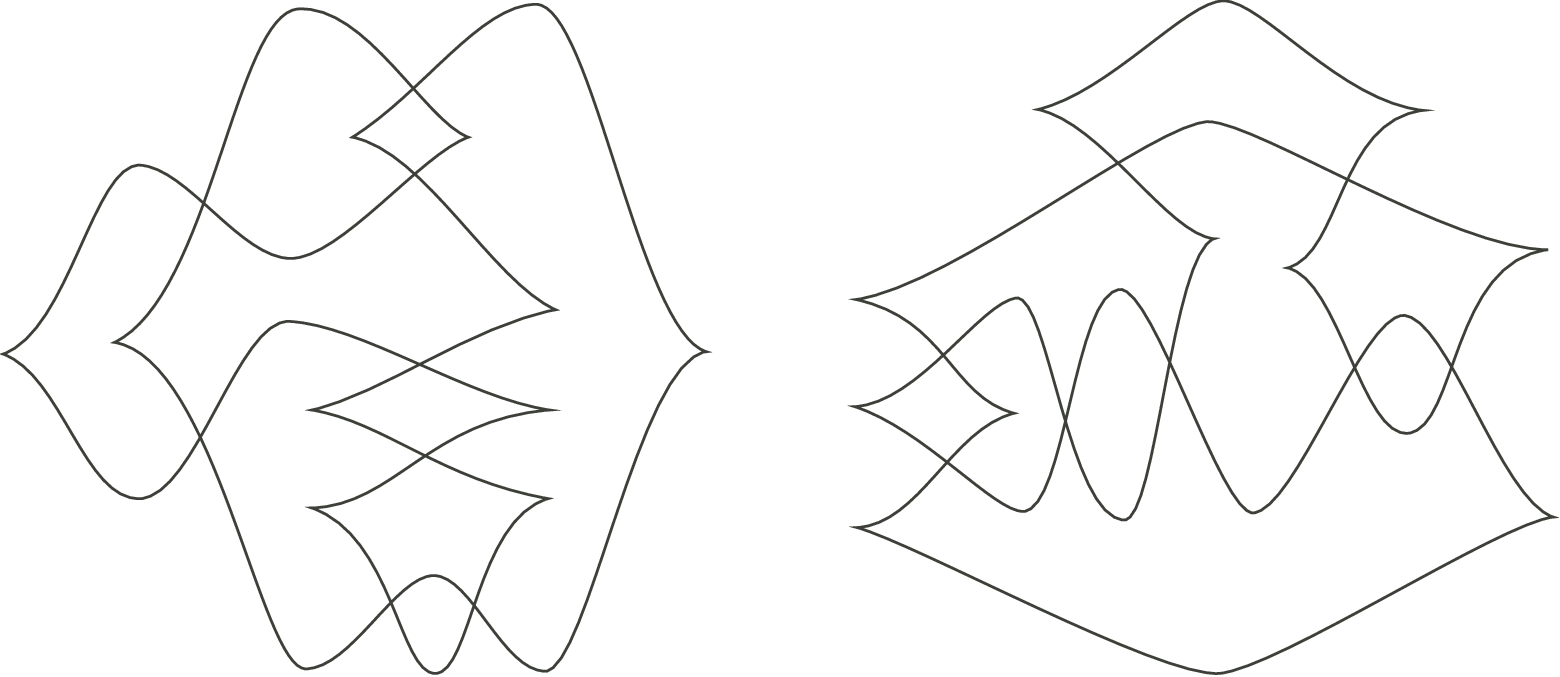}{.14}	&\up{2t}
&\up{(-5,0)}	\\ [1ex]
\up{8_{13}}	&\up{(-4,1)}	&			&\fig{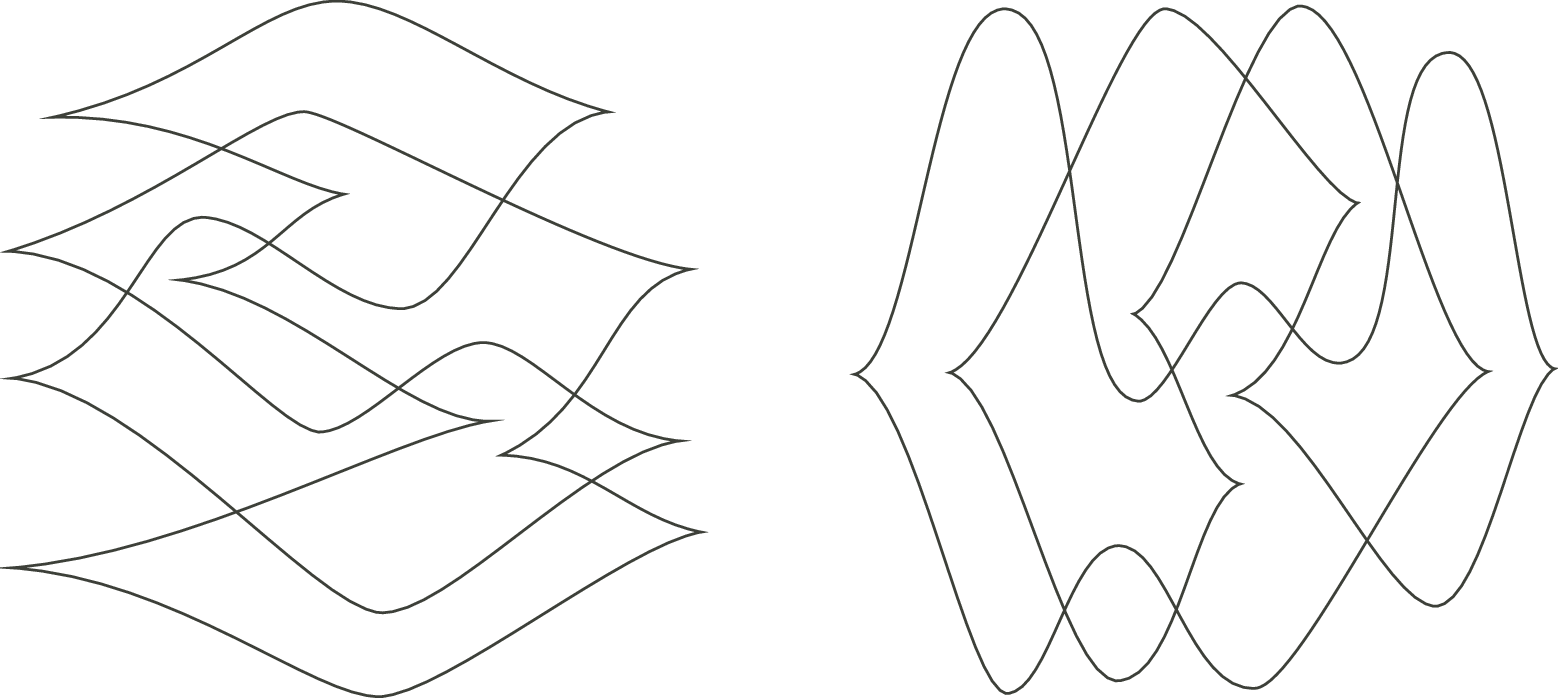}{.14}
&			&\up{(-6,1)}	\\ [1ex]
\up{8_{14}}	&\up{(-9,2)}	&			&\fig{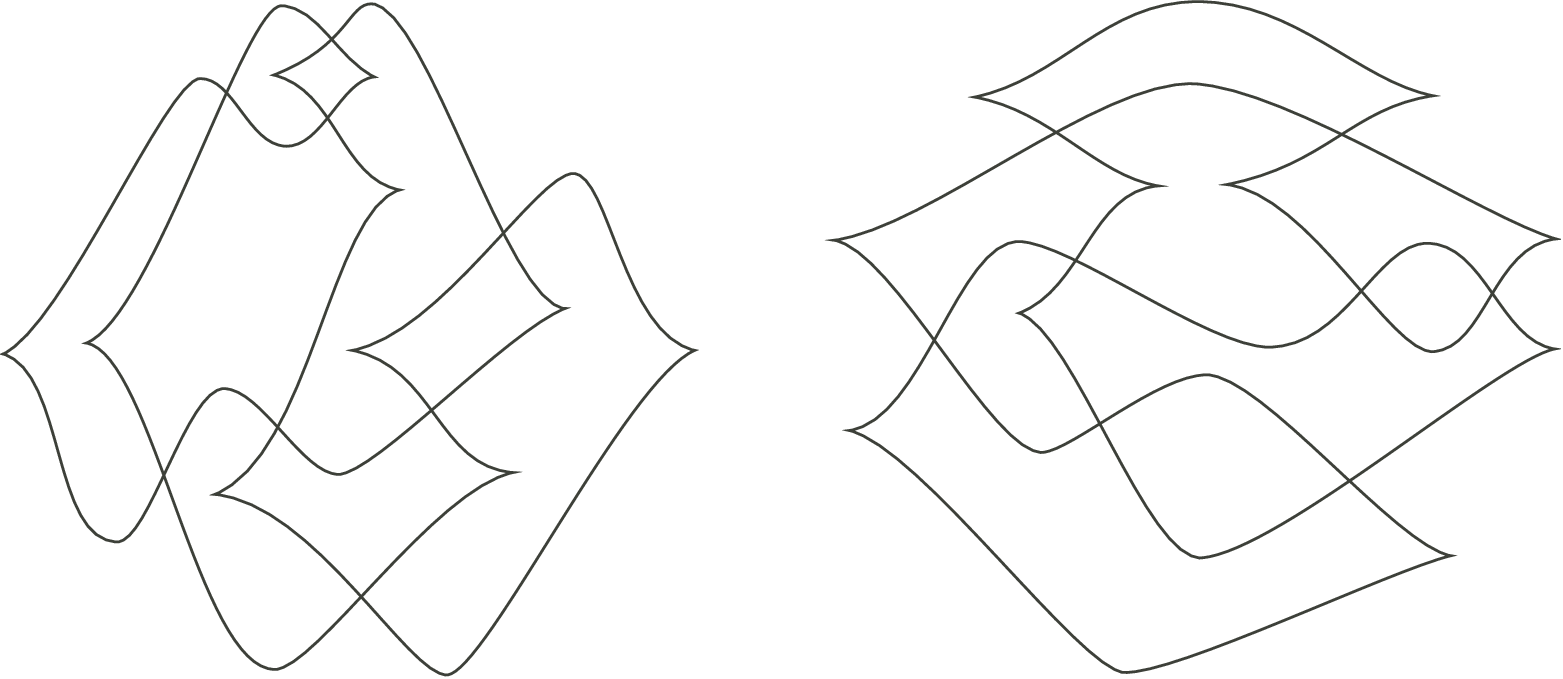}{.14}
&\up{1+t}		&\up{(-1,0)}	\\ [1ex]
\up{8_{15}}	&\up{(-13,0)}	&			&\fig{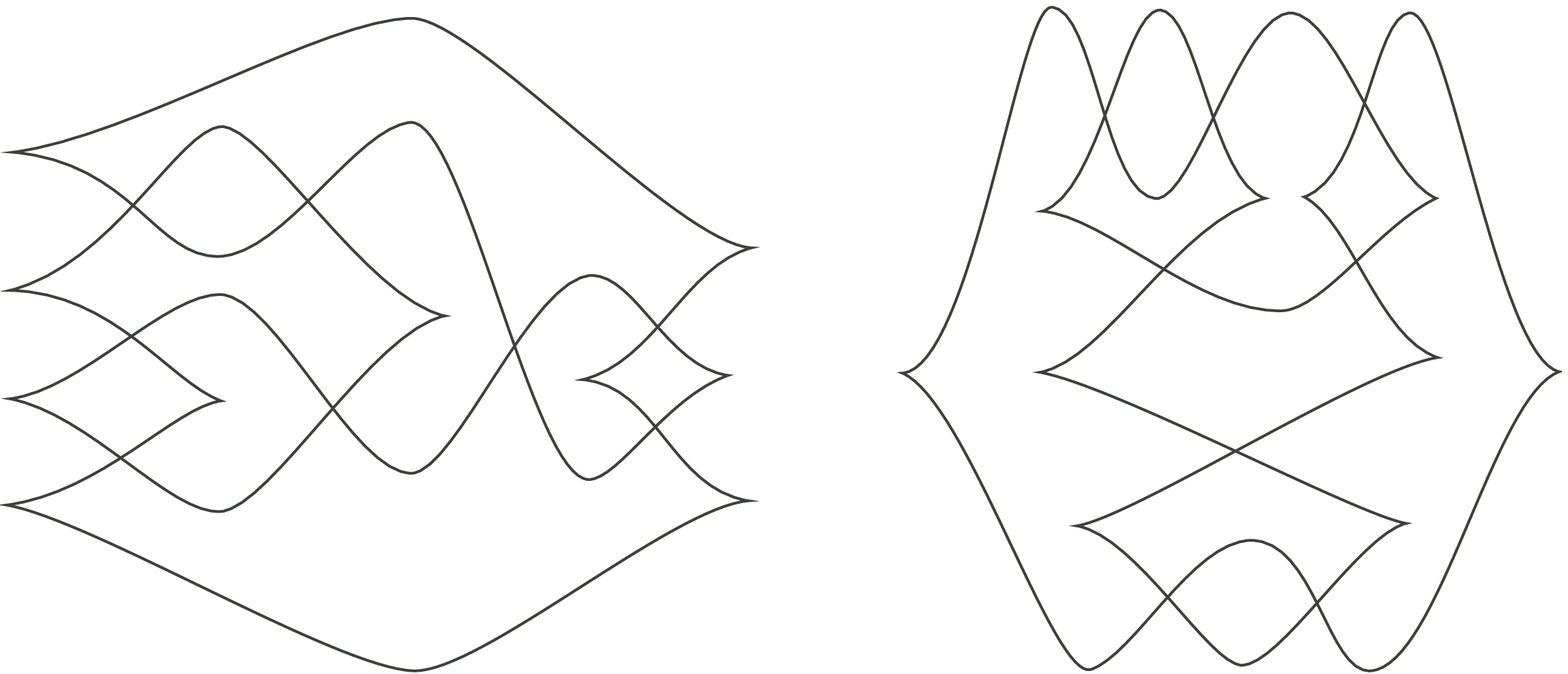}{.14}
&\up{1+t^2}	&\up{(3,0)}	\\ [1ex]
\hline
\end{tabular}
\end{table}

\clearpage

\begin{table}[!ht]
\tabletop
\up{8_{16}}	&\up{(-8,1)}	&			&\fig{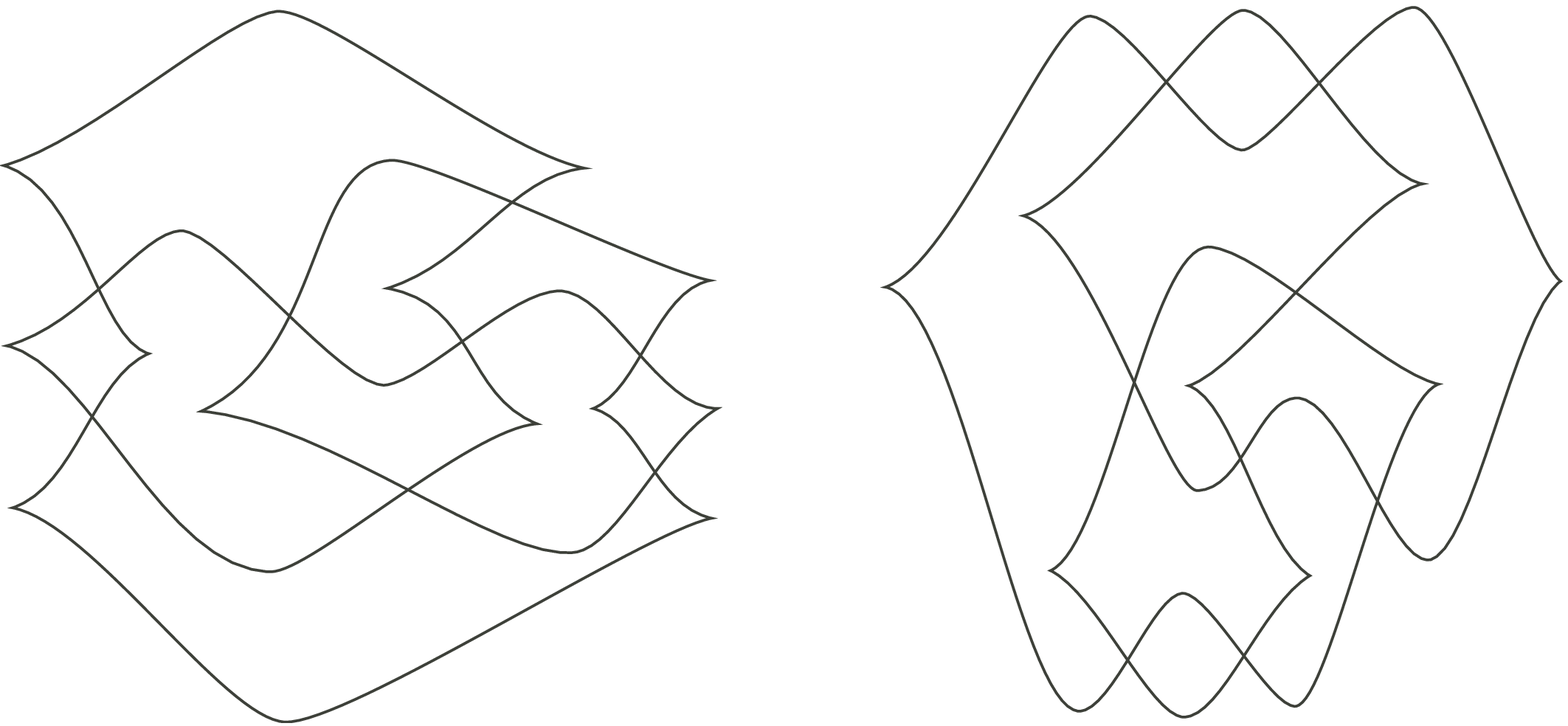}{.14}
&			&\up{(-2,1)}	\\ [1ex]
\up{8_{17}^*}	&\up{(-5,0)}	&			&\fig{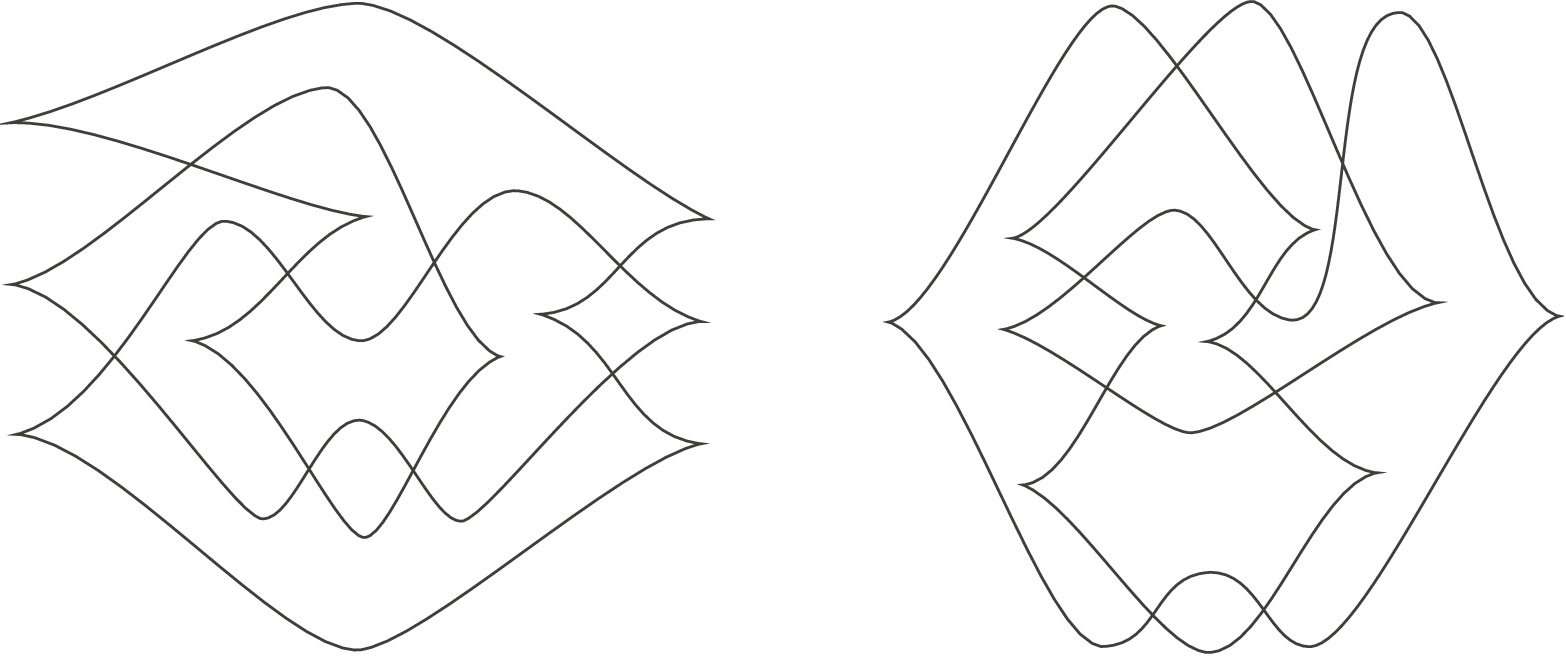}{.14}
&			&\up{(-5,0)}	\\ [1ex]
\up{8_{18}^*}	&\up{(-5,0)}	&			&\fig{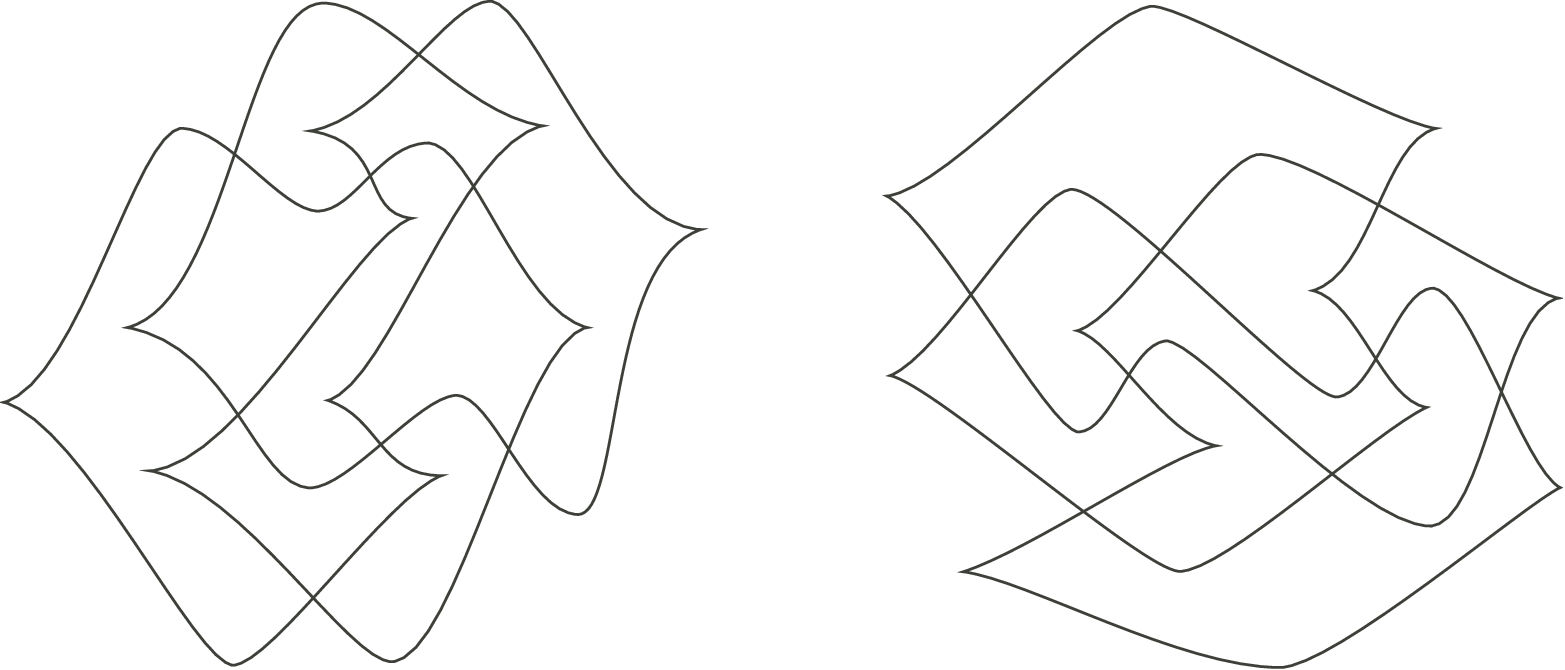}{.14}
&			&\up{(-5,0)}	\\ [1ex]
\up{8_{19}}	&\up{(5,0)}	&\up{3}		&\fig{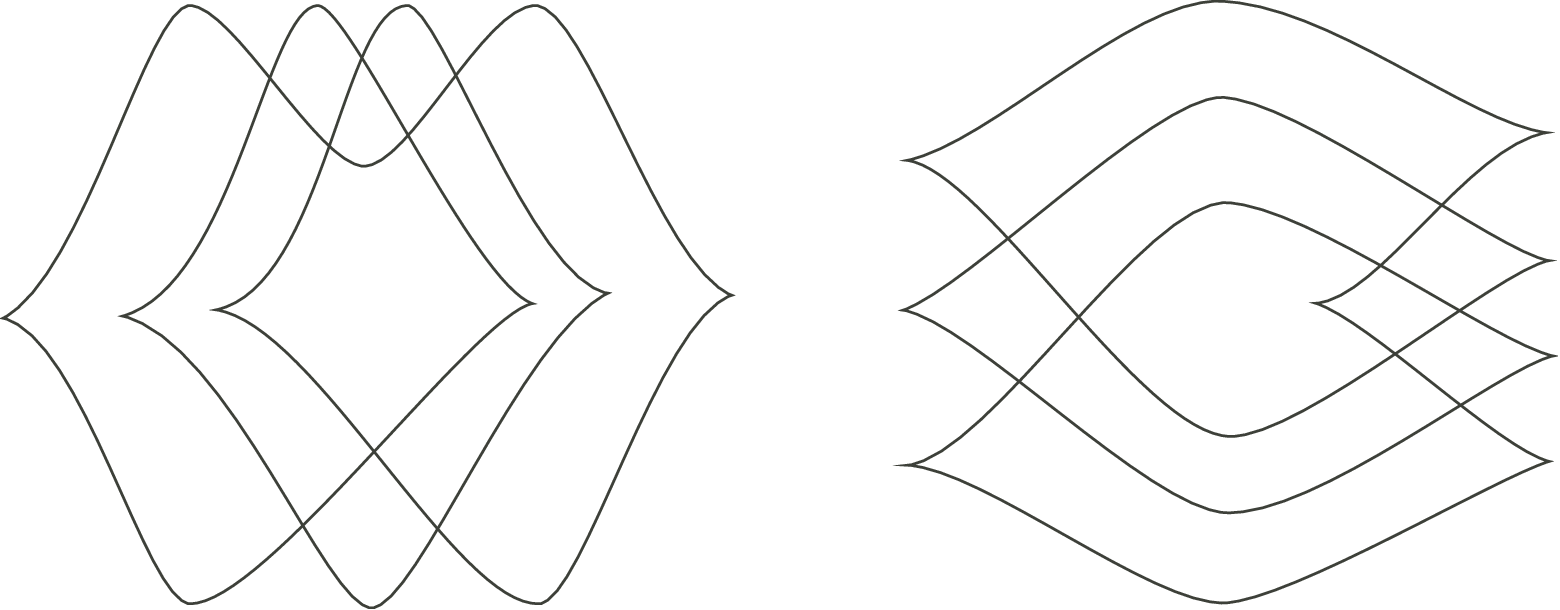}{.14}	&
&\up{(-12,1)}	\\ [1ex]
\up{8_{20}}	&\up{(-6,1)}	&			&\fig{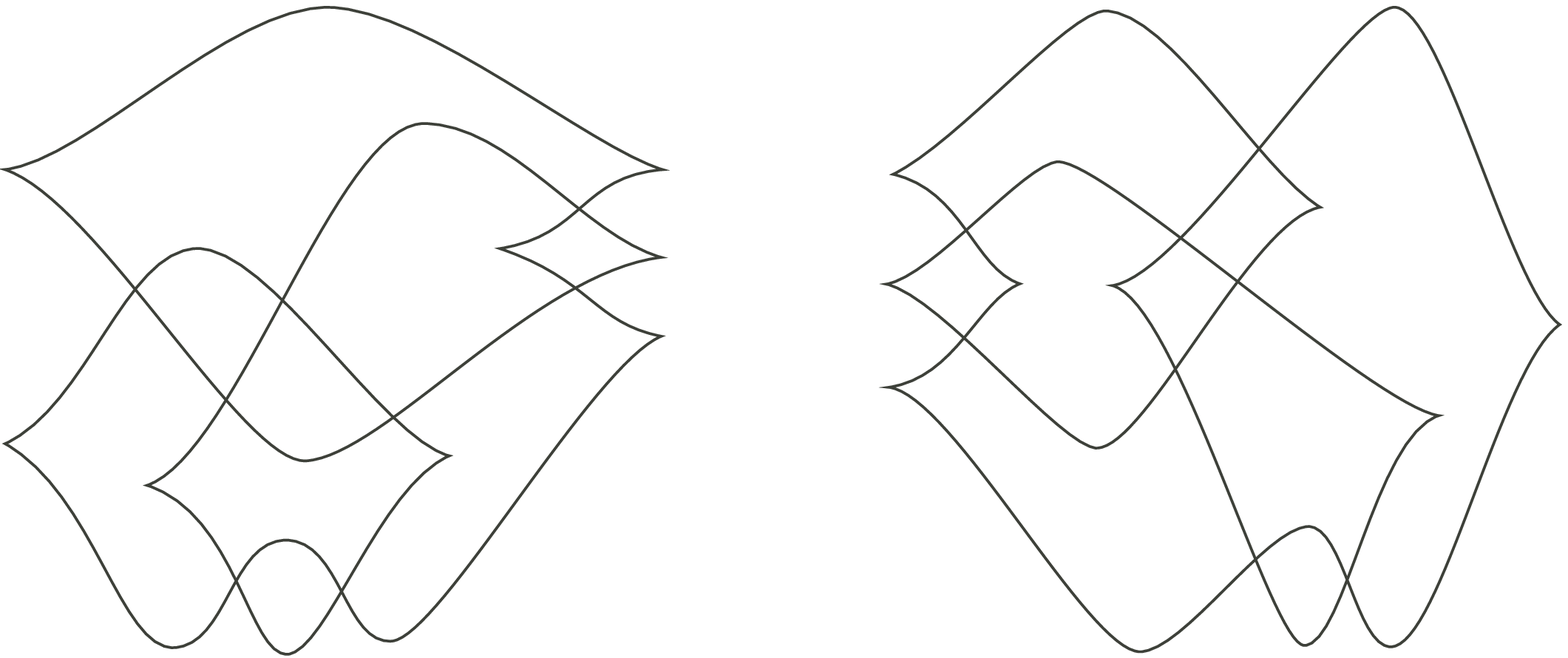}{.14}
&			&\up{(-2,1)}	\\ [1ex]
\up{8_{21}} &\up{(-9,0)}	&	&\fig{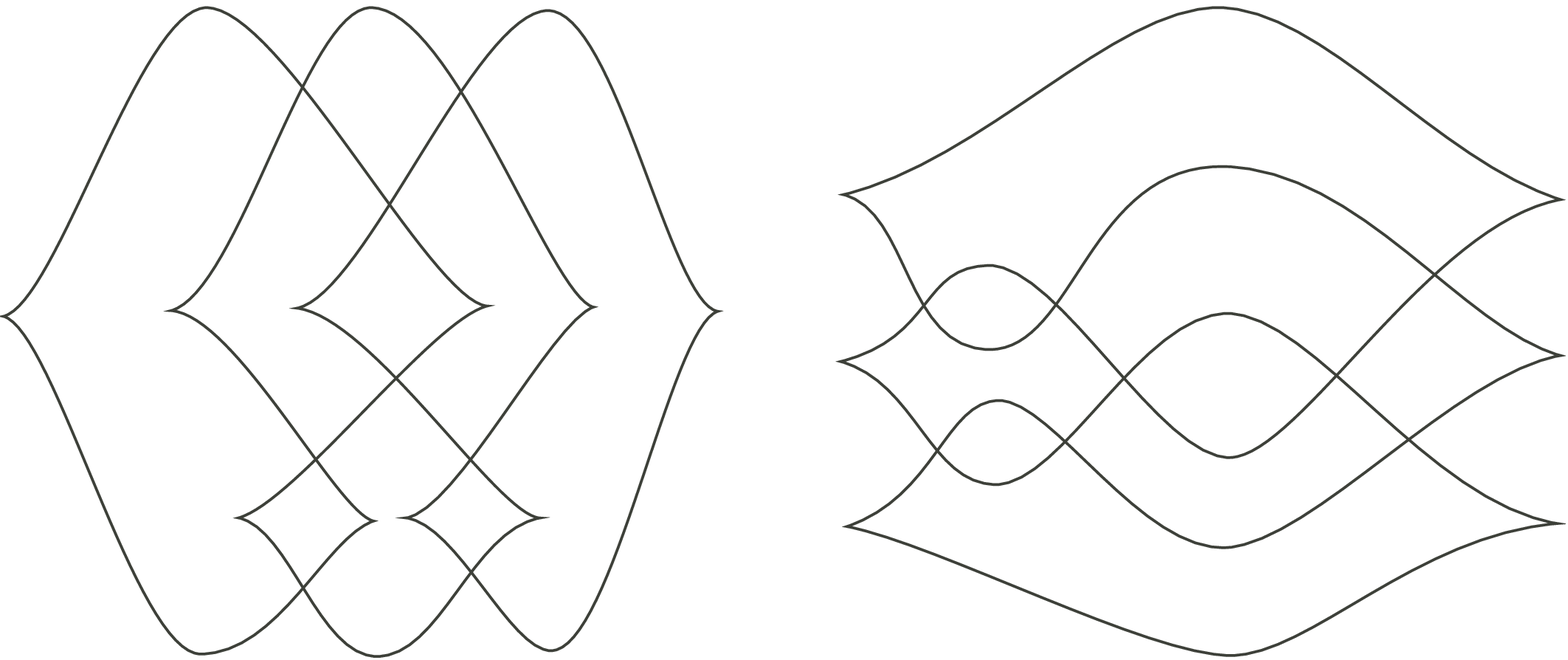}{.14}
&\up{\underset{\textstyle 2+t}{\overset{\textstyle 1}{}}}
															&\up{(1,0)}
															\\
															[1ex]
\up{9_1}		&\up{(-18,7)}	&
&\fig{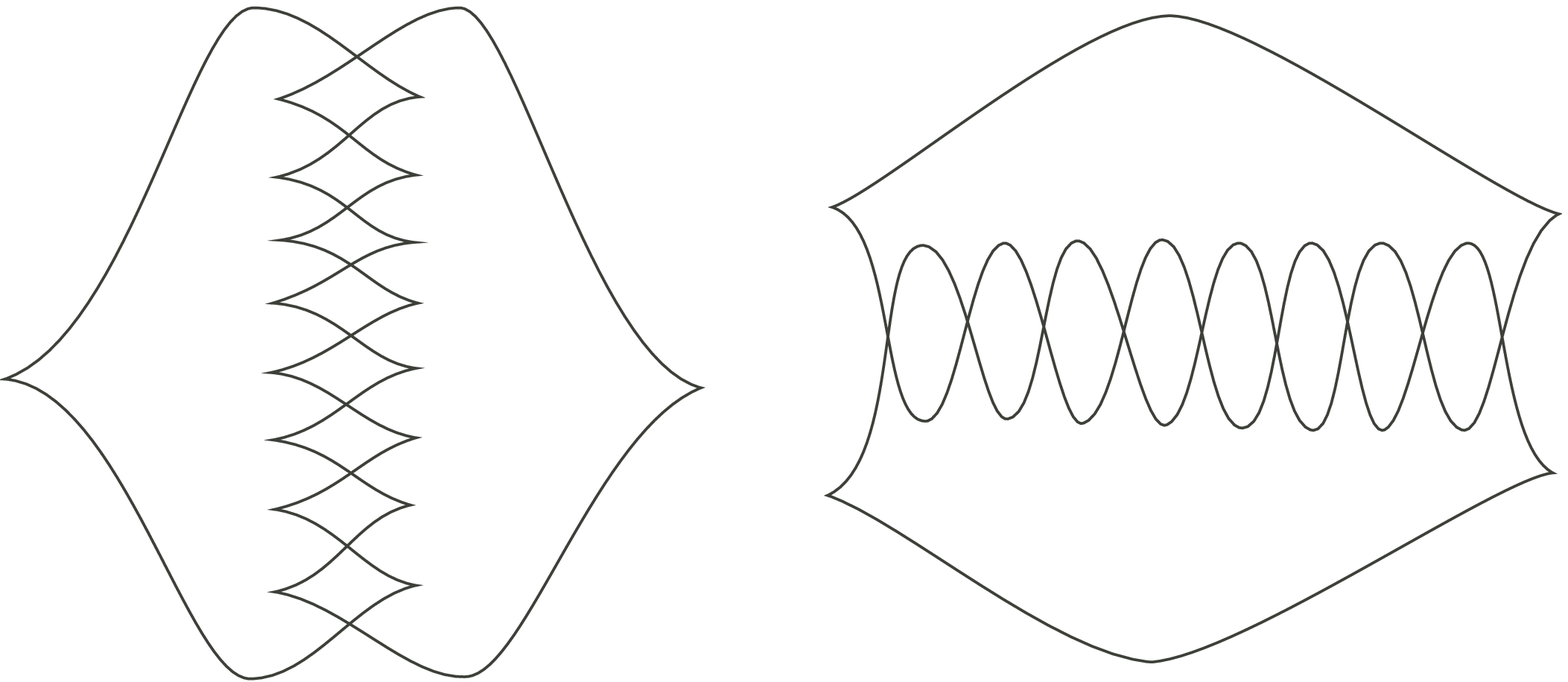}{.14}	&\up{4}		&\up{(7,0)}	\\ [1ex]
\up{9_2}		&\up{(-12,1)}	&
&\fig{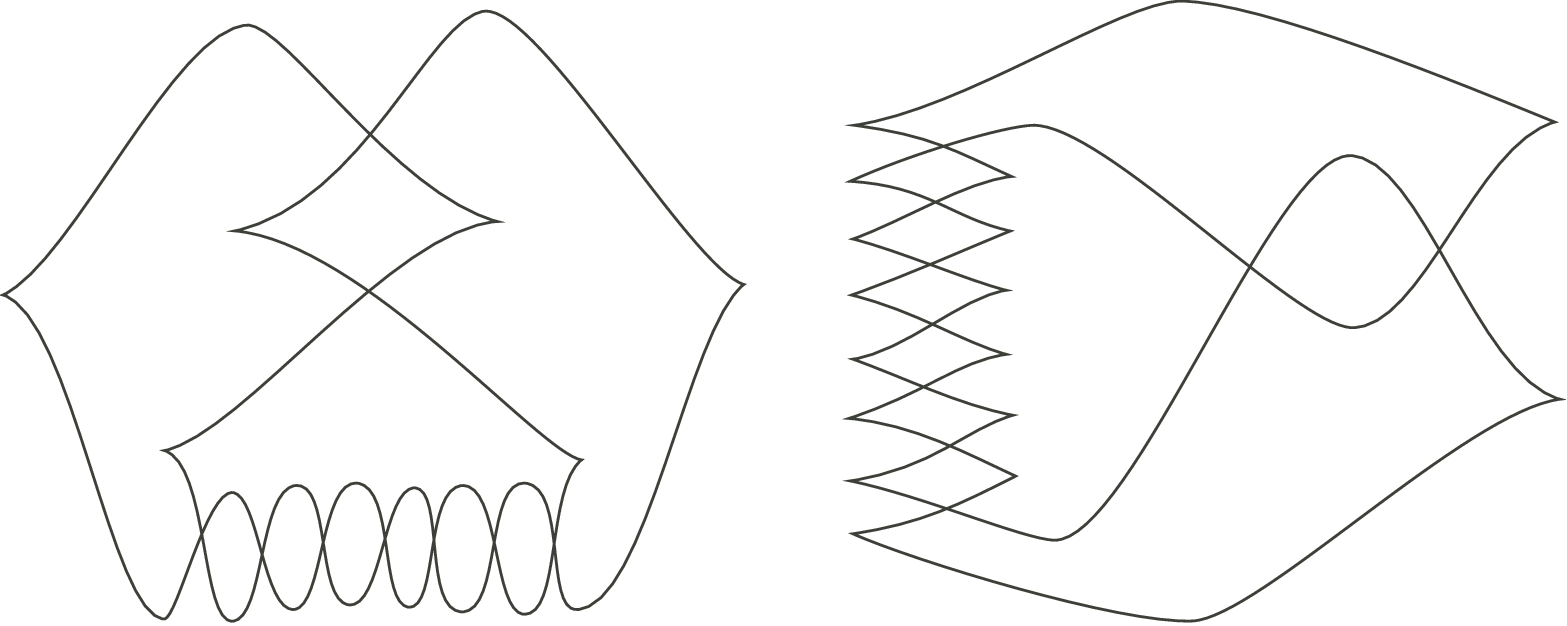}{.14}	&\up{t^6}		&\up{(1,0)}	\\ [1ex]
\up{9_3}		&\up{(5,0)}	&\up{3t^2}
&\fig{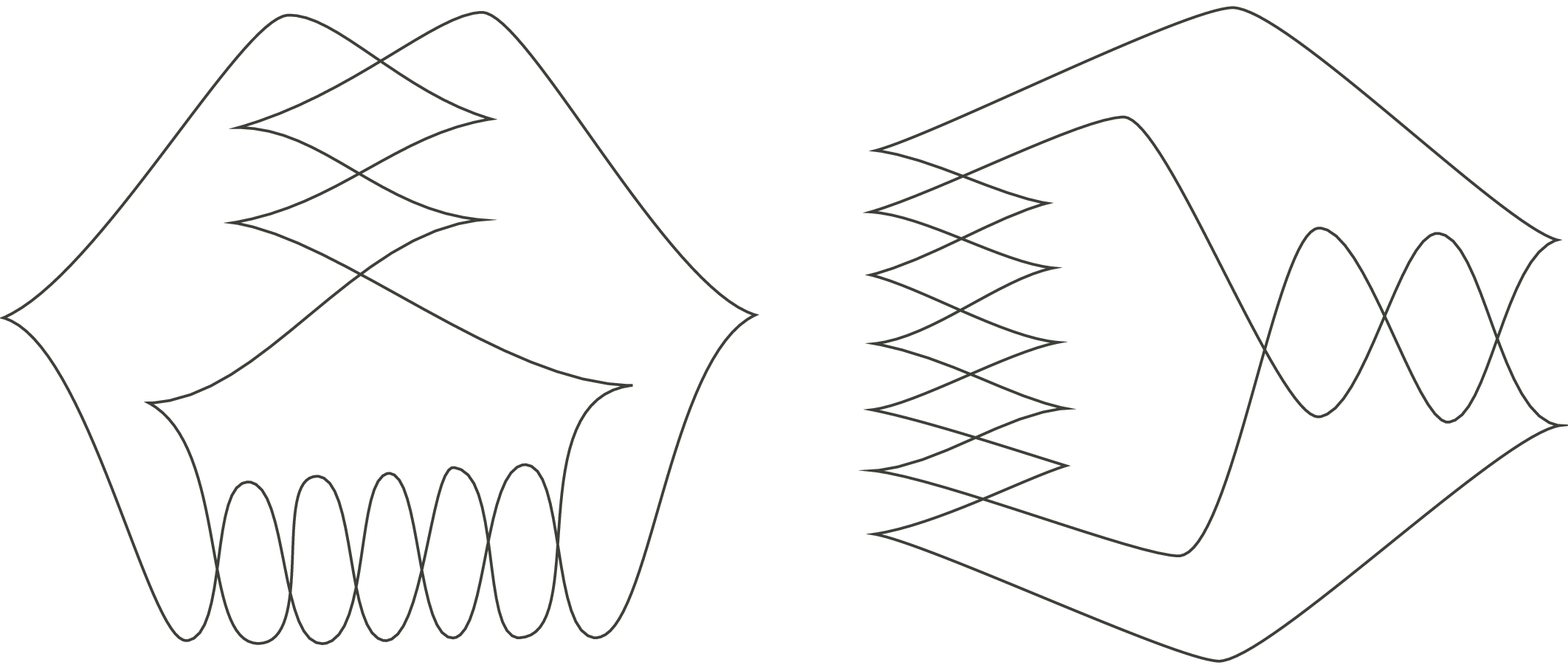}{.14}	&			&\up{(-16,5)}	\\ [1ex]
\up{9_4}		&\up{(-14,3)}	&
&\fig{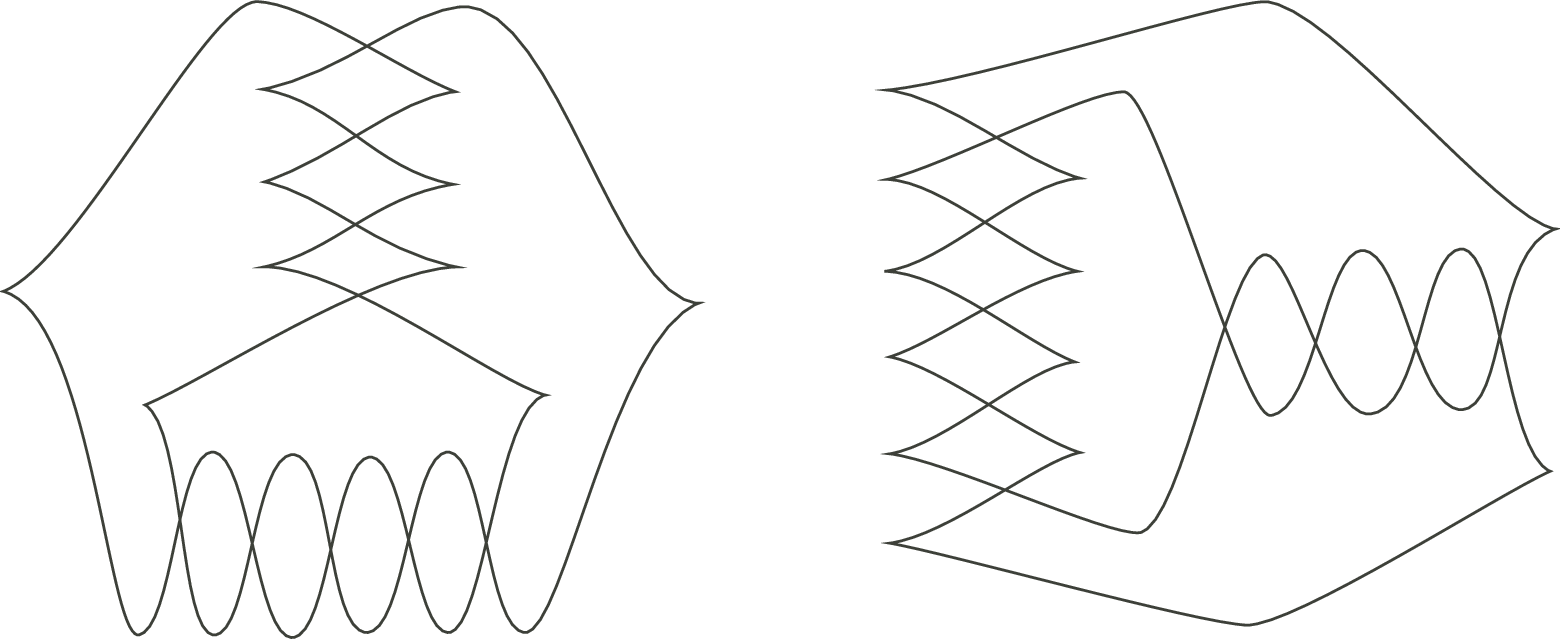}{.14}	&\up{2t^4}	&\up{(3,0)}	\\ [1ex]
\hline
\end{tabular}
\end{table}

\clearpage

\begin{table}[!ht]
\tabletop
\up{9_5}		&\up{(1,0)}	&
&\fig{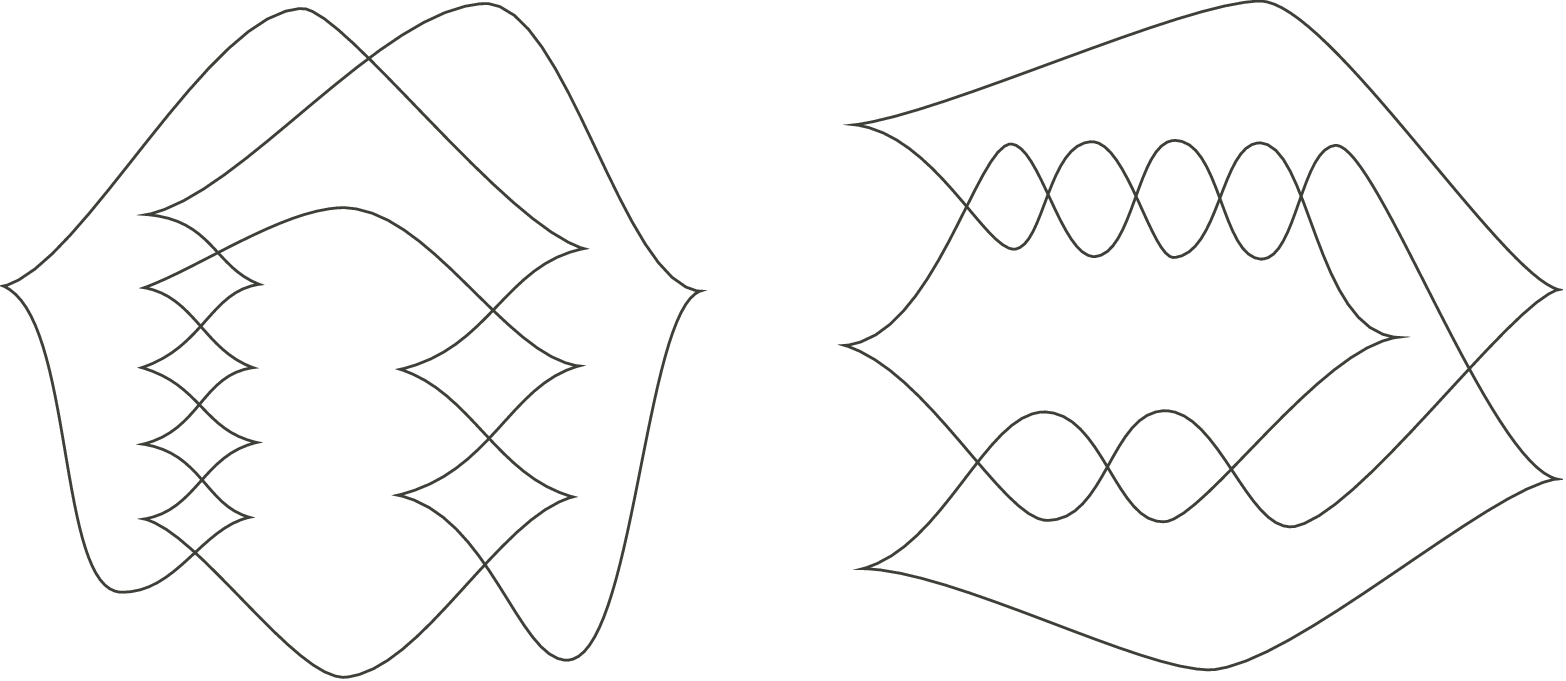}{.14}	&			&\up{(-12,1)}	\\ [1ex]
\up{9_6}		&\up{(-16,3)}	&
&\fig{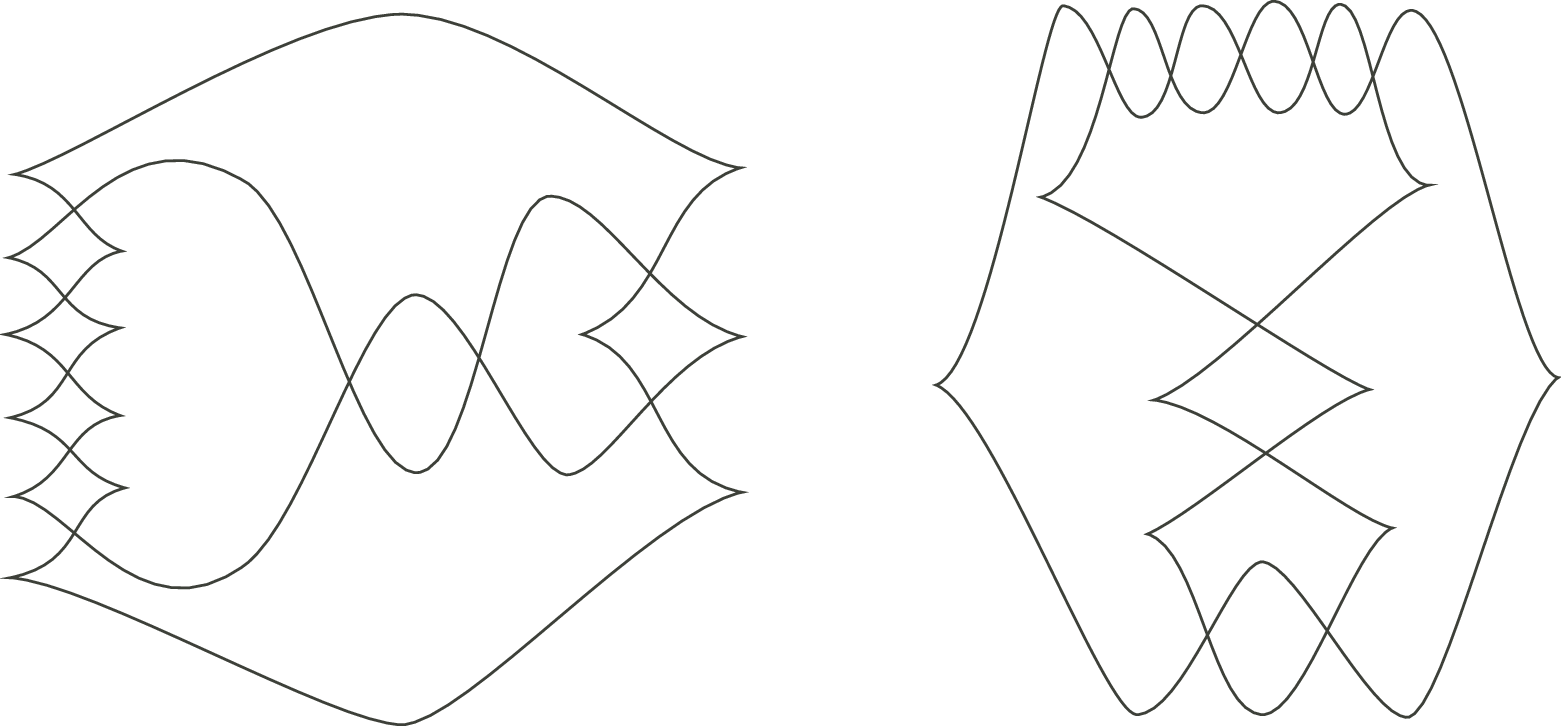}{.14}	&\up{2+t^2}	&\up{(5,0)}	\\ [1ex]
\up{9_7}		&\up{(-14,1)}	&
&\fig{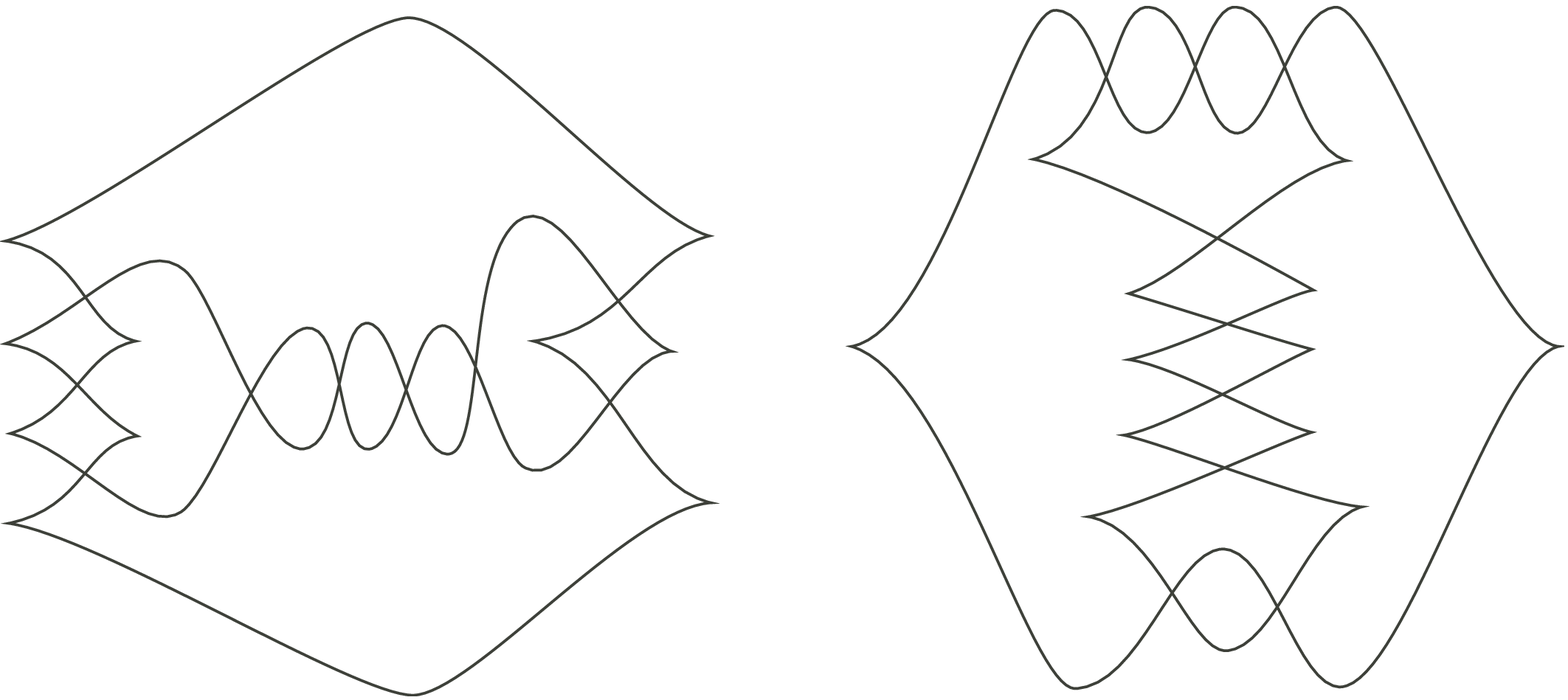}{.14}	&\up{1+t^4}	&\up{(3,0)}	\\ [1ex]
\up{9_8}		&\up{(-8,1)}	&
&\fig{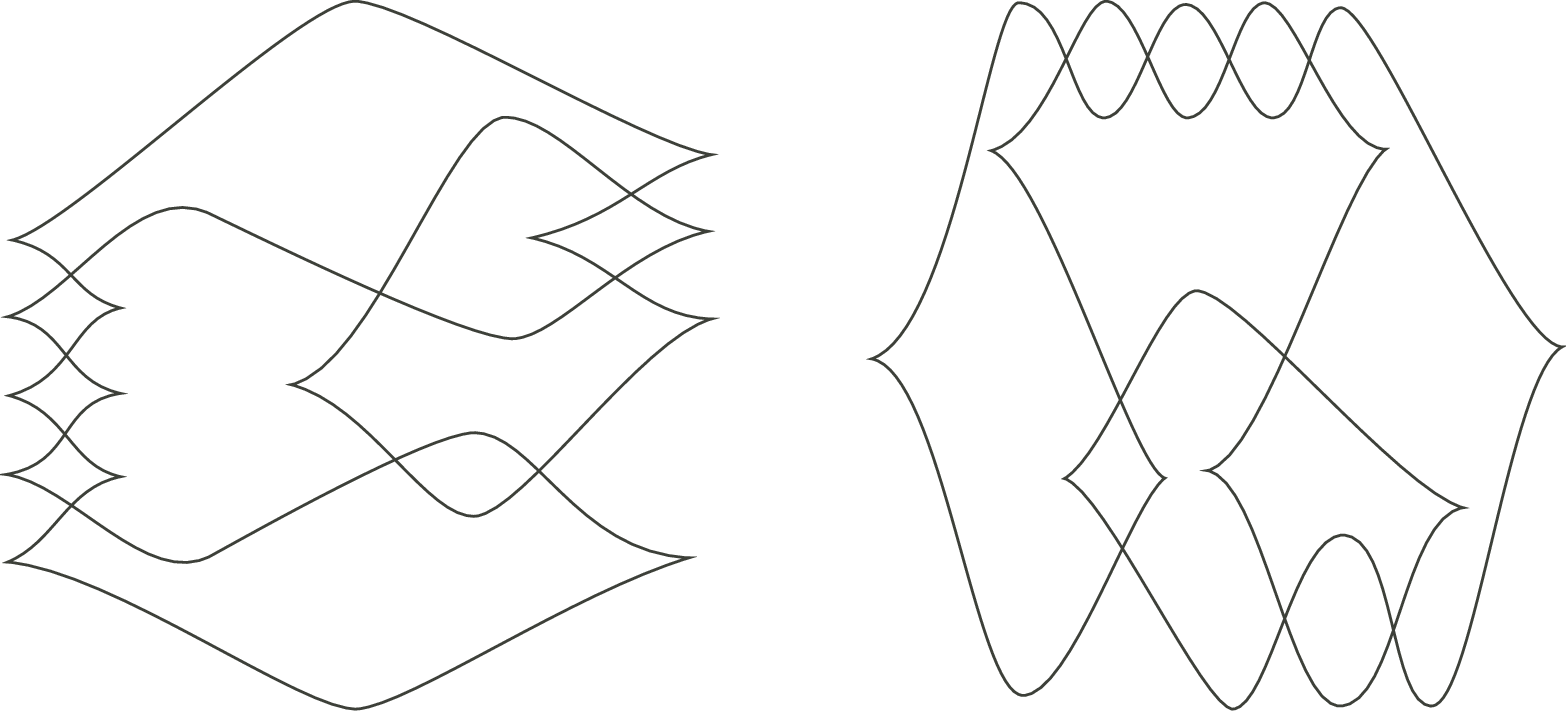}{.14}	&\up{1+2t}	&\up{(-3,0)}	\\ [1ex]
\up{9_9}		&\up{(-16,1)}	&
&\fig{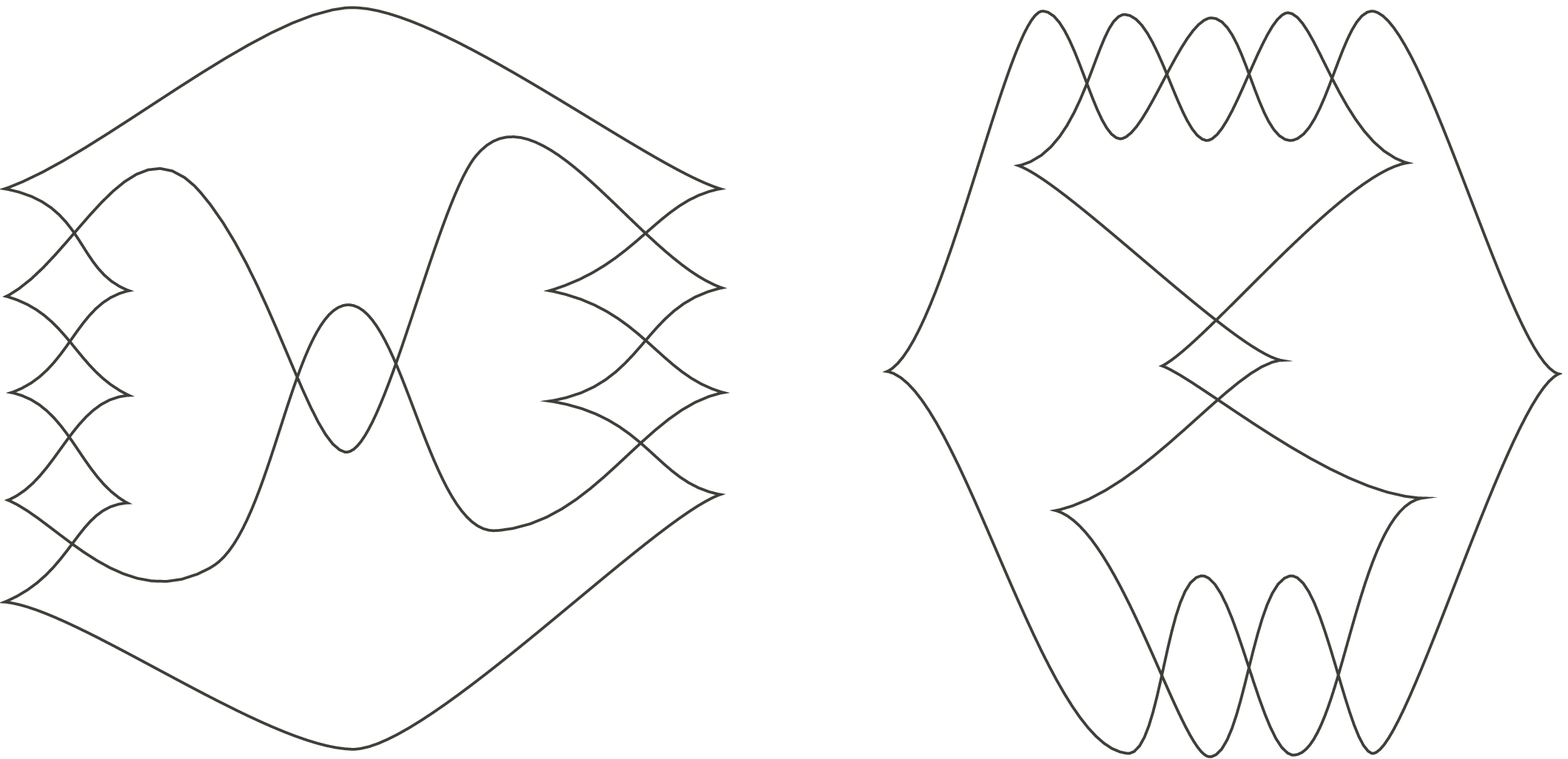}{.14}	&\up{1+2t^2}	&\up{(5,0)}	\\ [1ex]
\up{9_{10}}	&\up{(3,0)}	&			&\fig{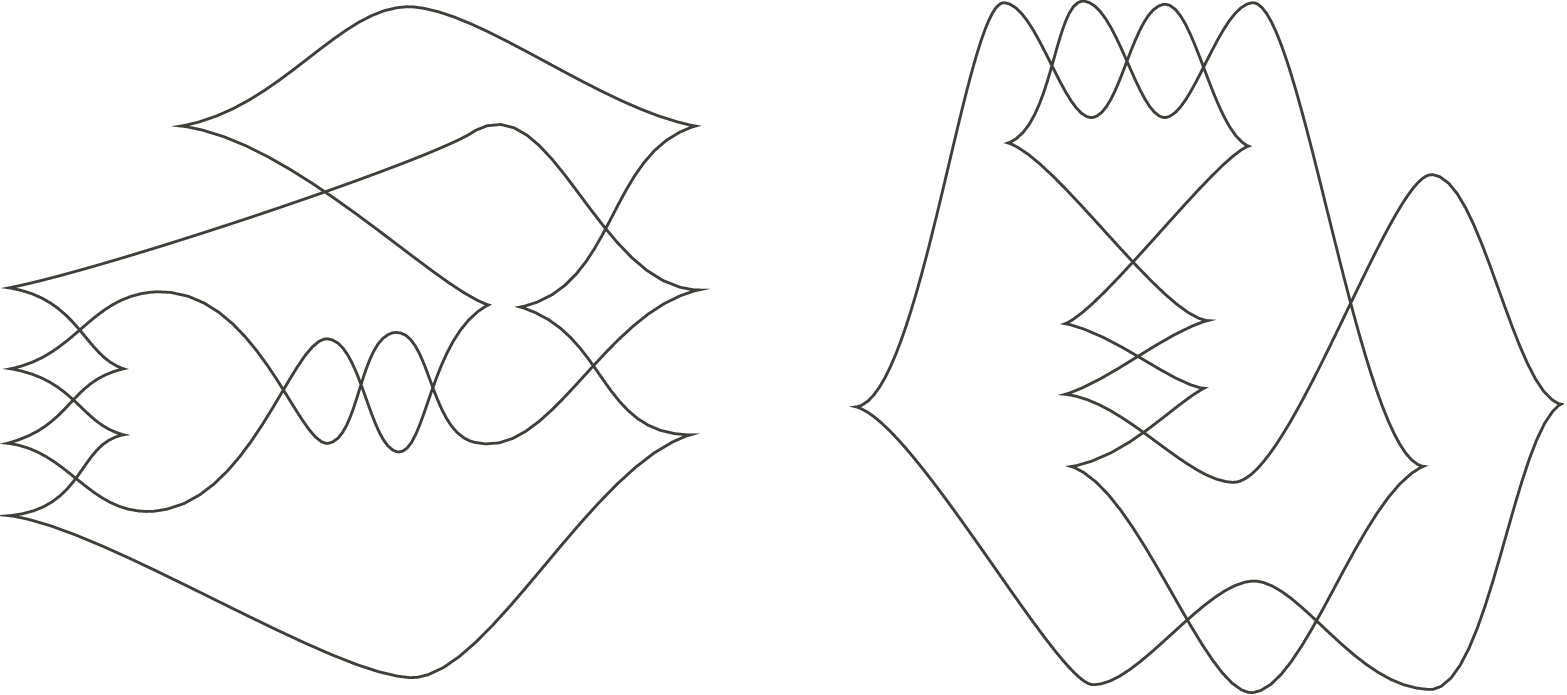}{.14}
&			&\up{(-14,3)}	\\ [1ex]
\up{9_{11}}	&\up{(1,0)}	&\up{t+2t^2}	&\fig{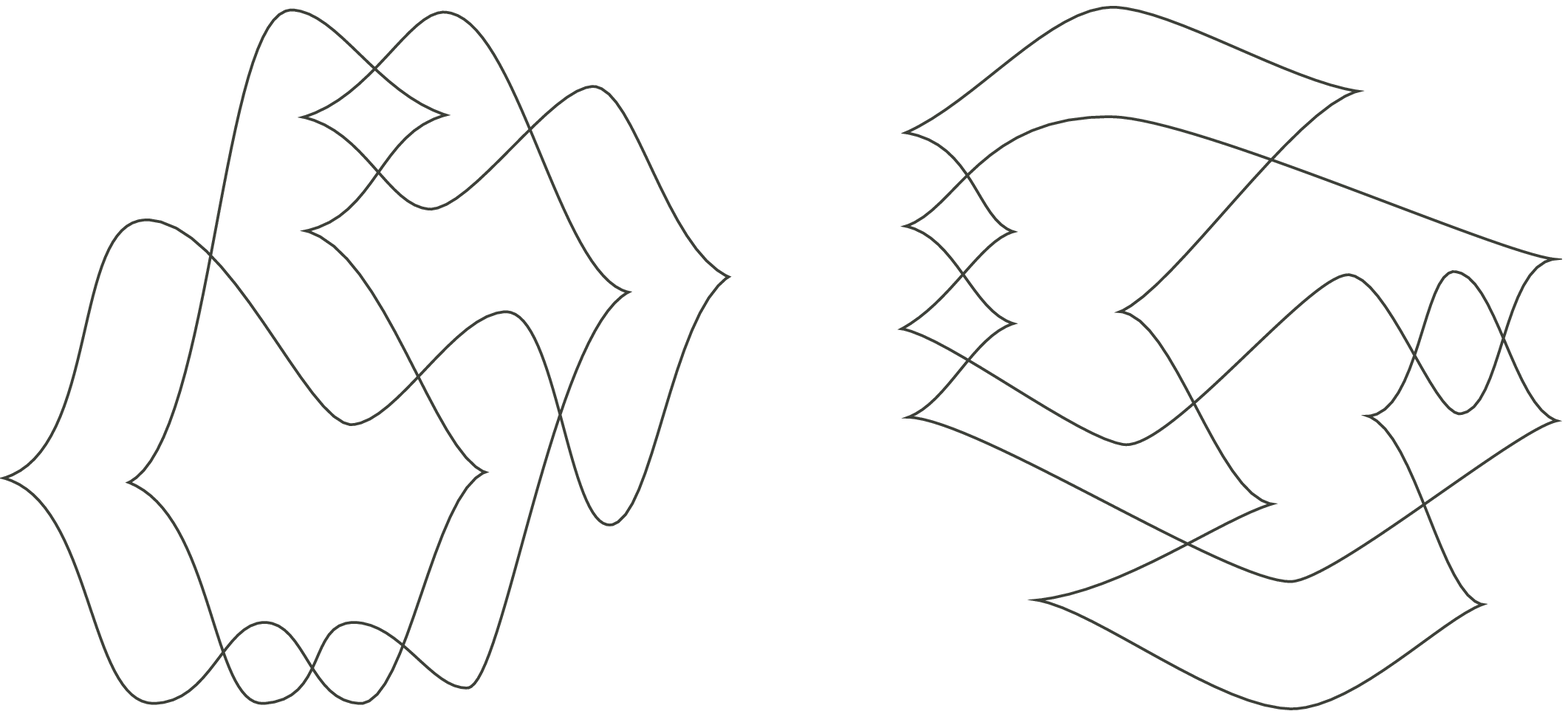}{.14}	&
&\up{(-12,1)}	\\ [1ex]
\up{9_{12}}	&\up{(-10,1)}	&			&\fig{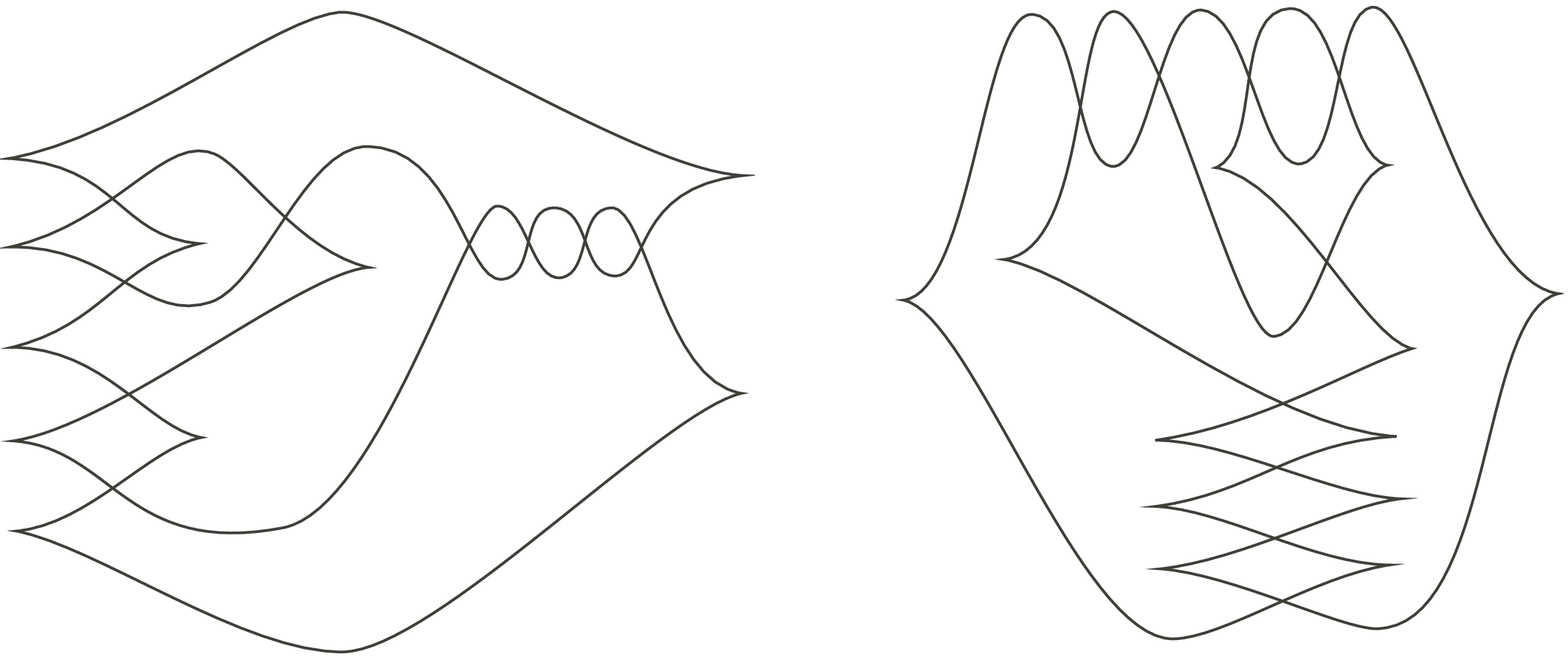}{.14}
&\up{t^3+t^4}	&\up{(-1,0)}	\\ [1ex]
\up{9_{13}}	&\up{(3,0)}	&			&\fig{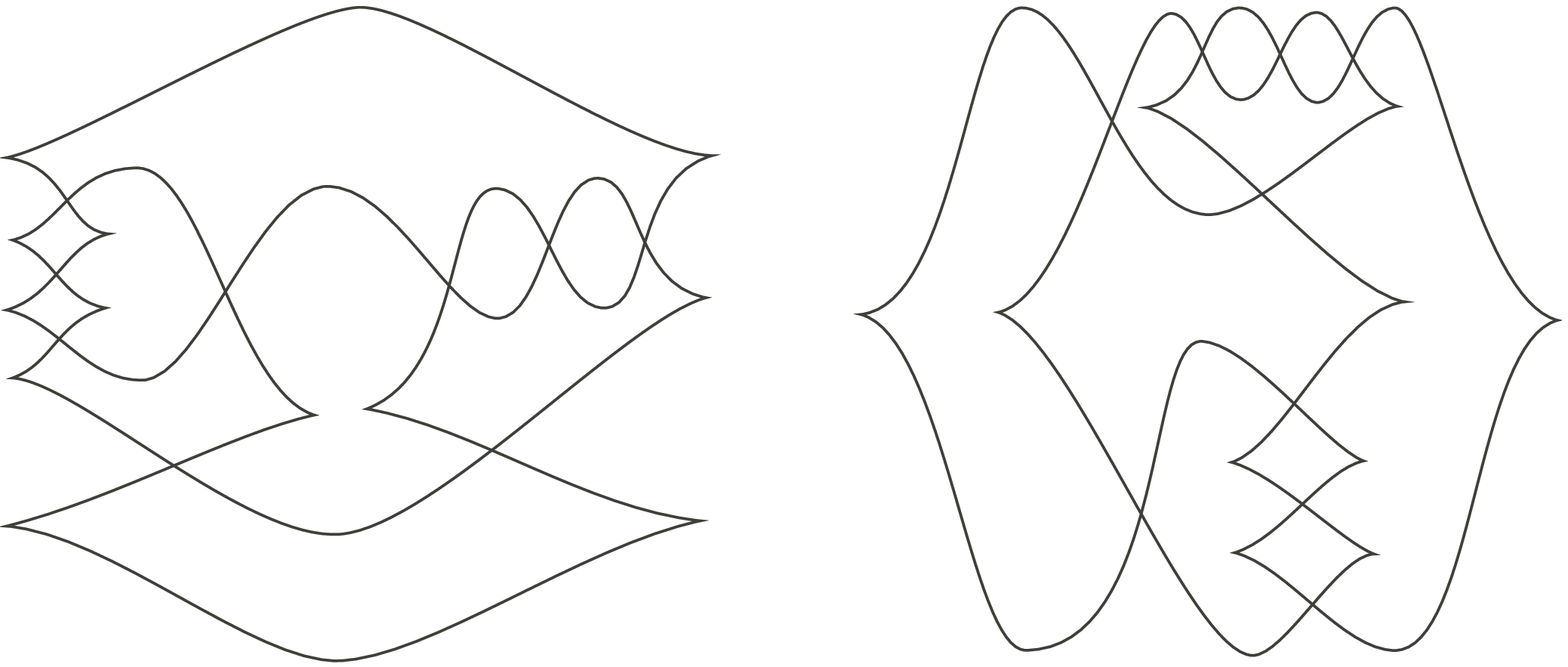}{.14}
&			&\up{(-14,3)}	\\ [1ex]
\up{9_{14}}	&\up{(-4,1)}	&			&\fig{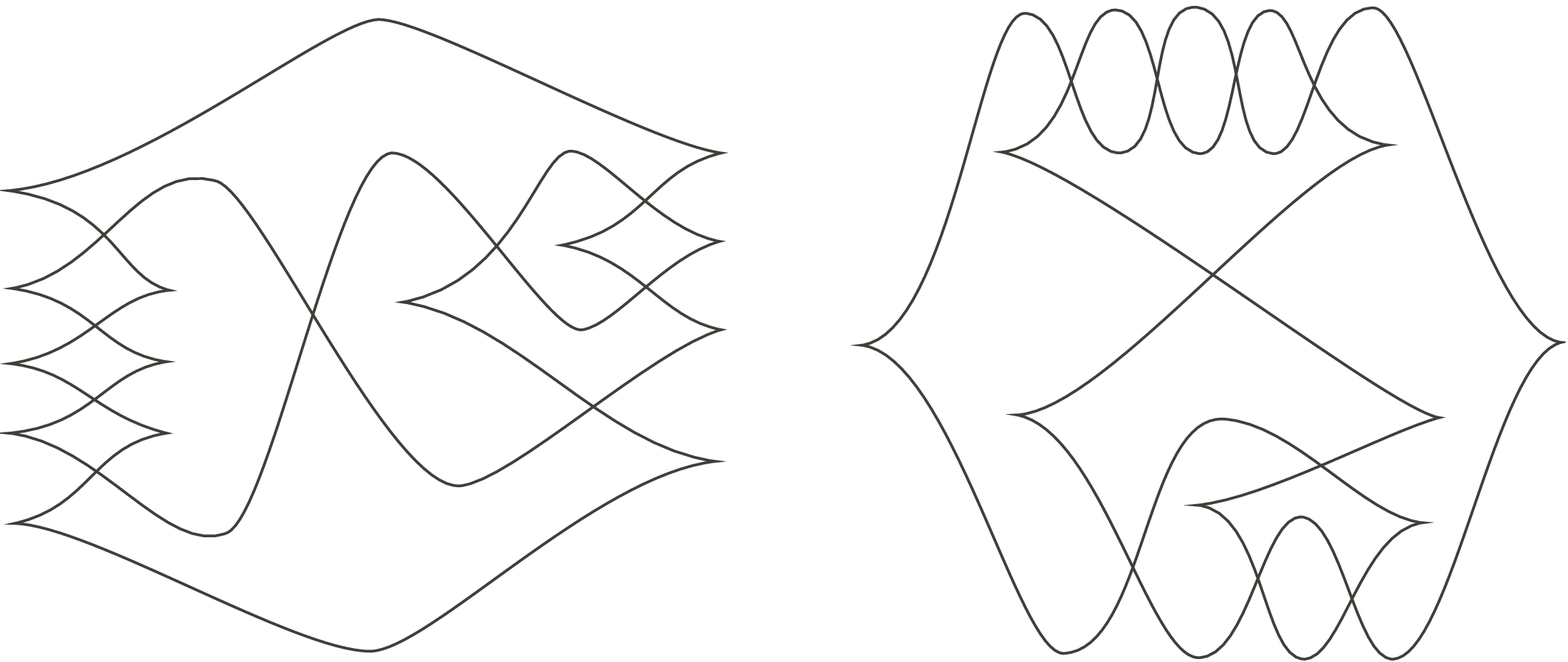}{.14}
&\up{3t}		&\up{(-7,0)}	\\ [1ex]
\hline
\end{tabular}
\end{table}

\clearpage

\begin{table}[!ht]
\tabletop
\up{9_{15}}	&\up{(-1,0)}	&\up{t+t^2}	&\fig{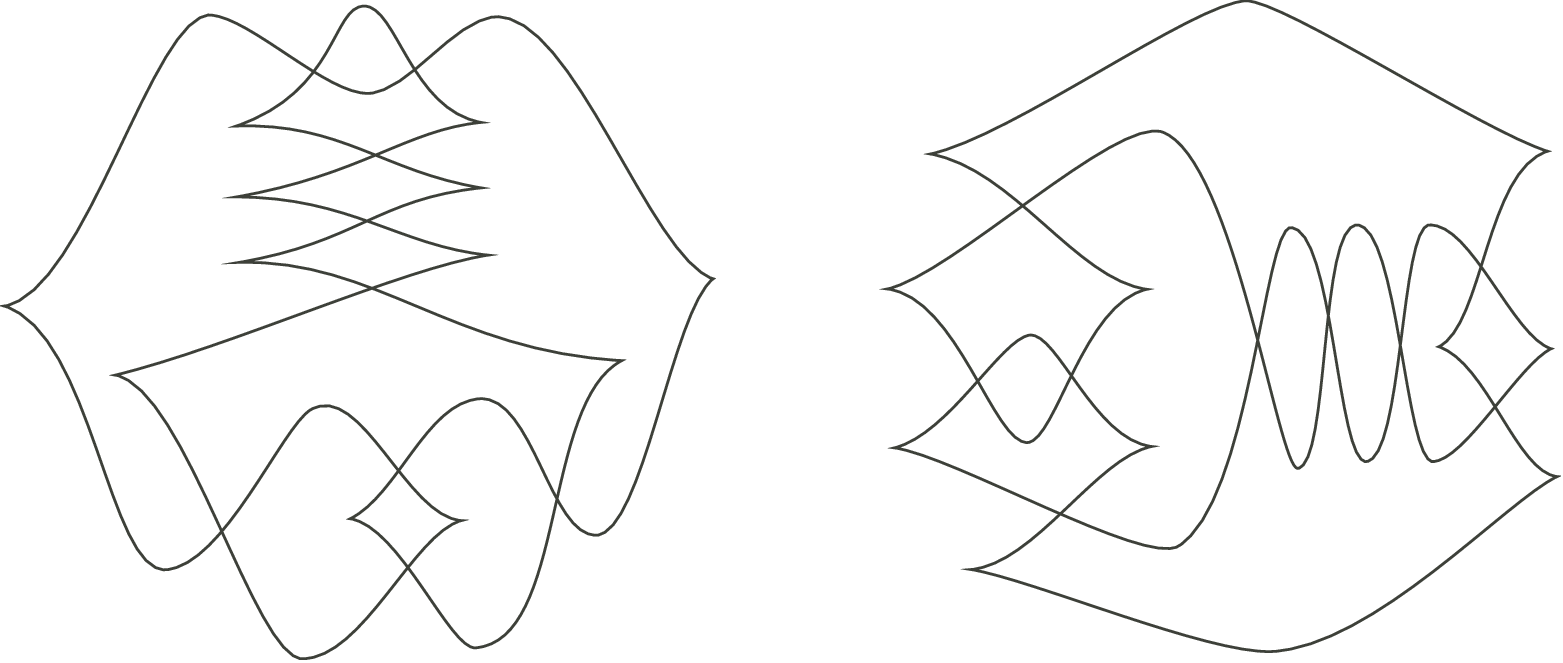}{.14}	&
&\up{(-10,1)}	\\ [1ex]
\up{9_{16}}	&\up{(5,0)}	&\up{3}		&\fig{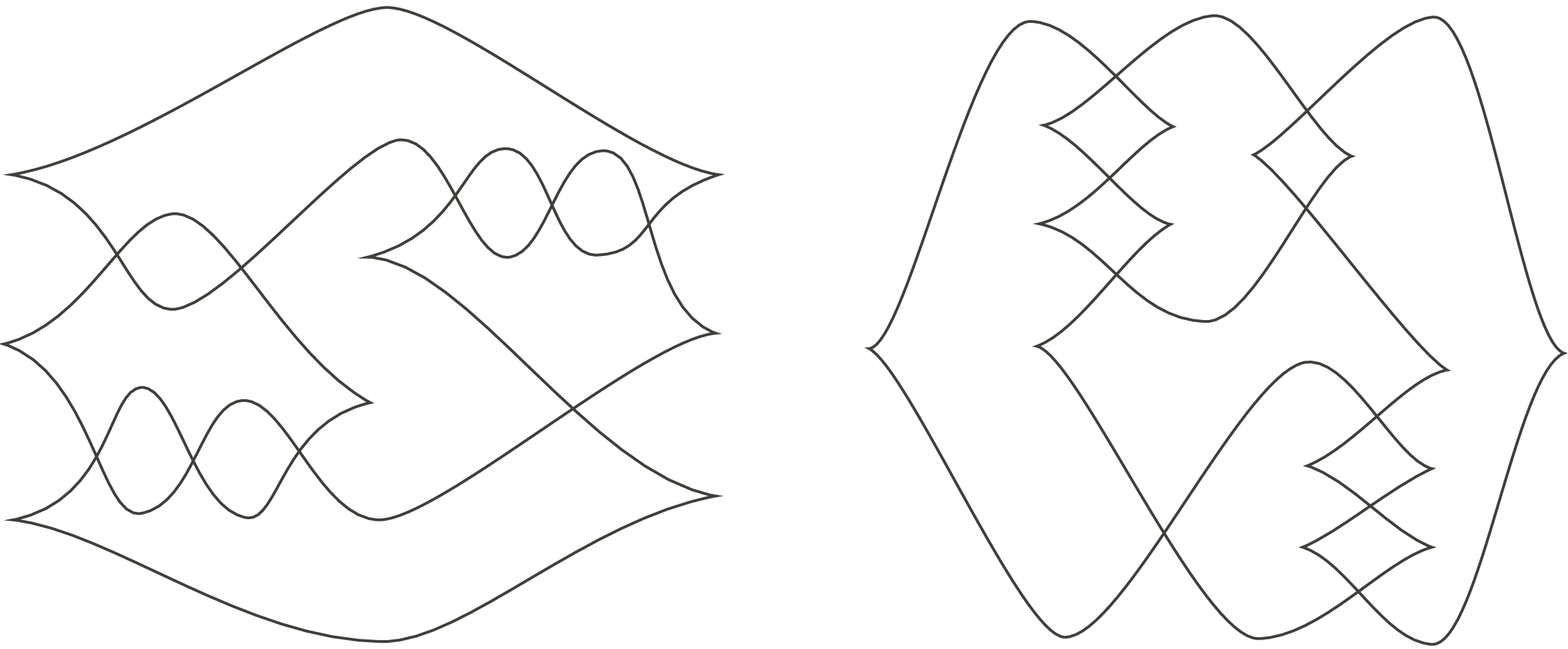}{.14}	&
&\up{(-16,1)}	\\ [1ex]
\up{9_{17}}	&\up{(-8,3)}	&			&\fig{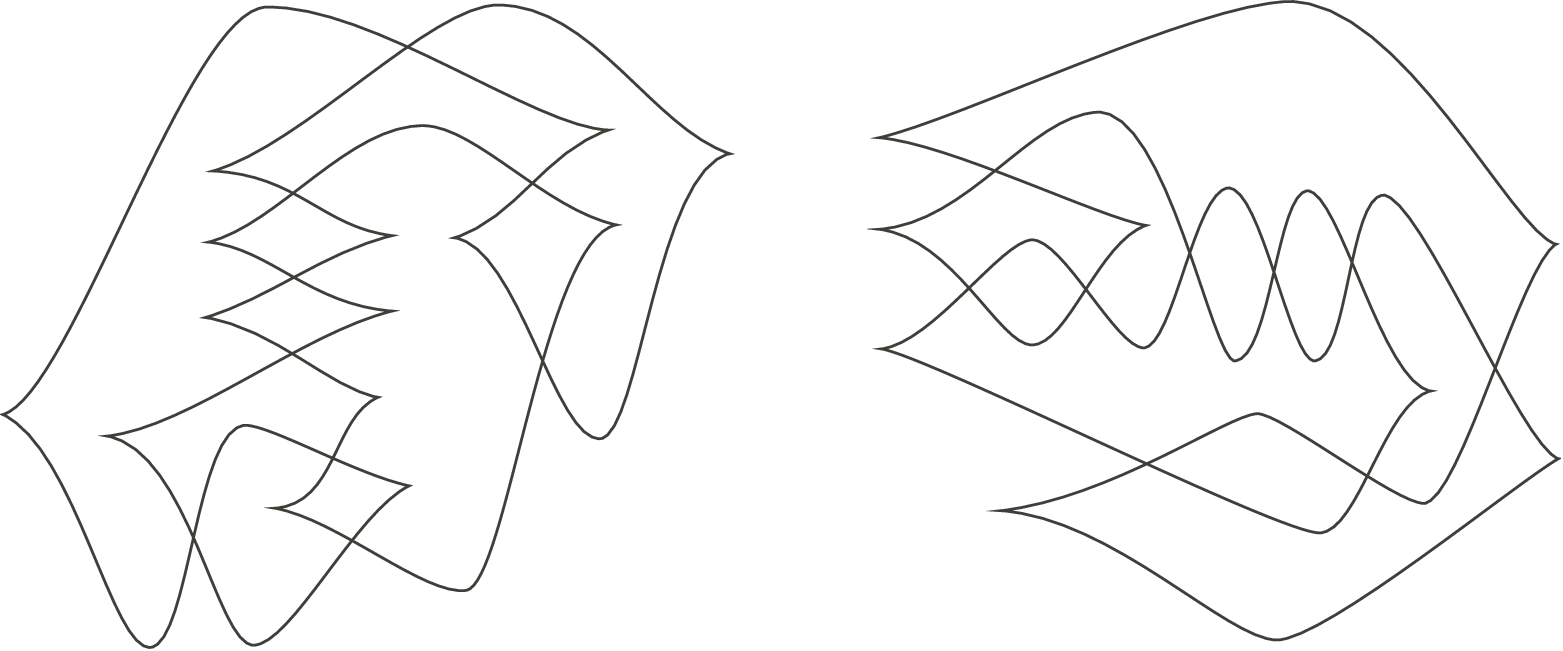}{.14}
&\up{1+2t}	&\up{(-3,0)}	\\ [1ex]
\up{9_{18}}	&\up{(-14,1)}	&			&\fig{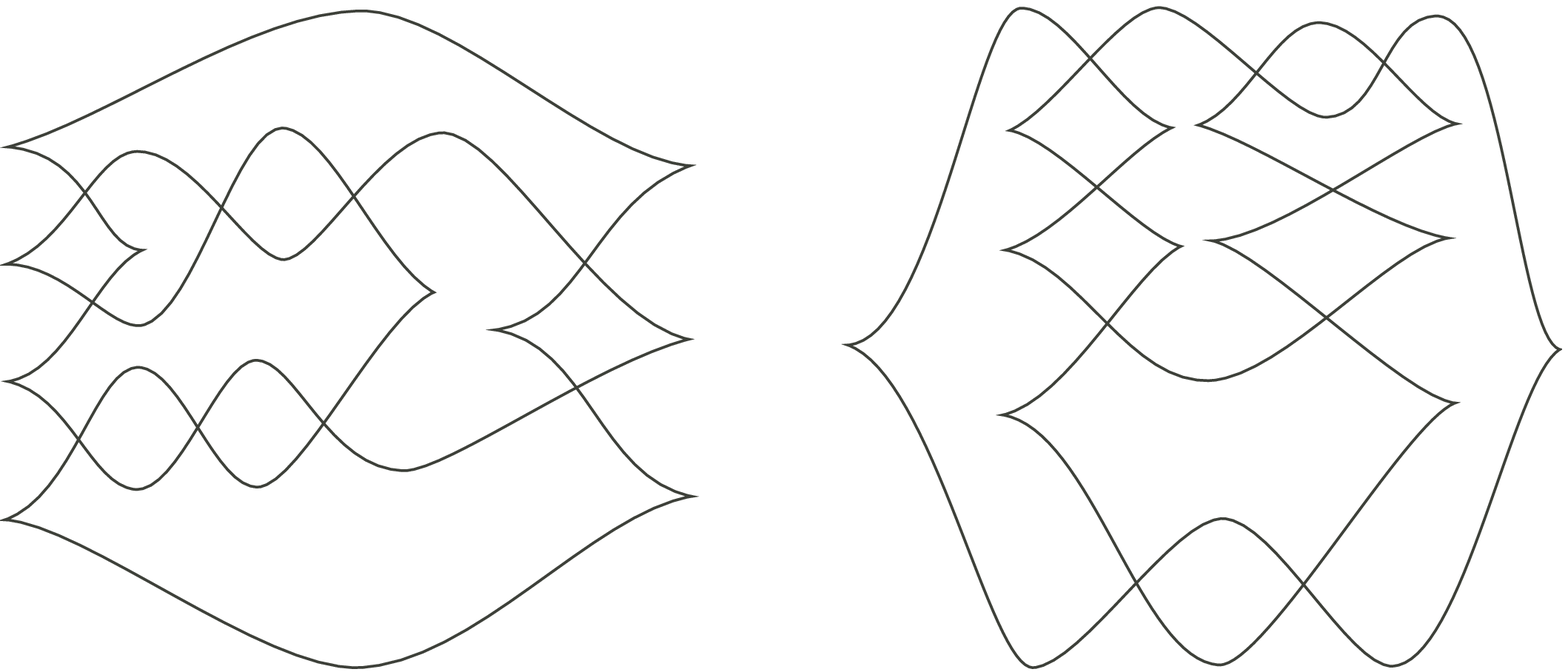}{.14}
&\up{1+t^2}	&\up{(3,0)}	\\ [1ex]
\up{9_{19}}	&\up{(-6,1)}	&			&\fig{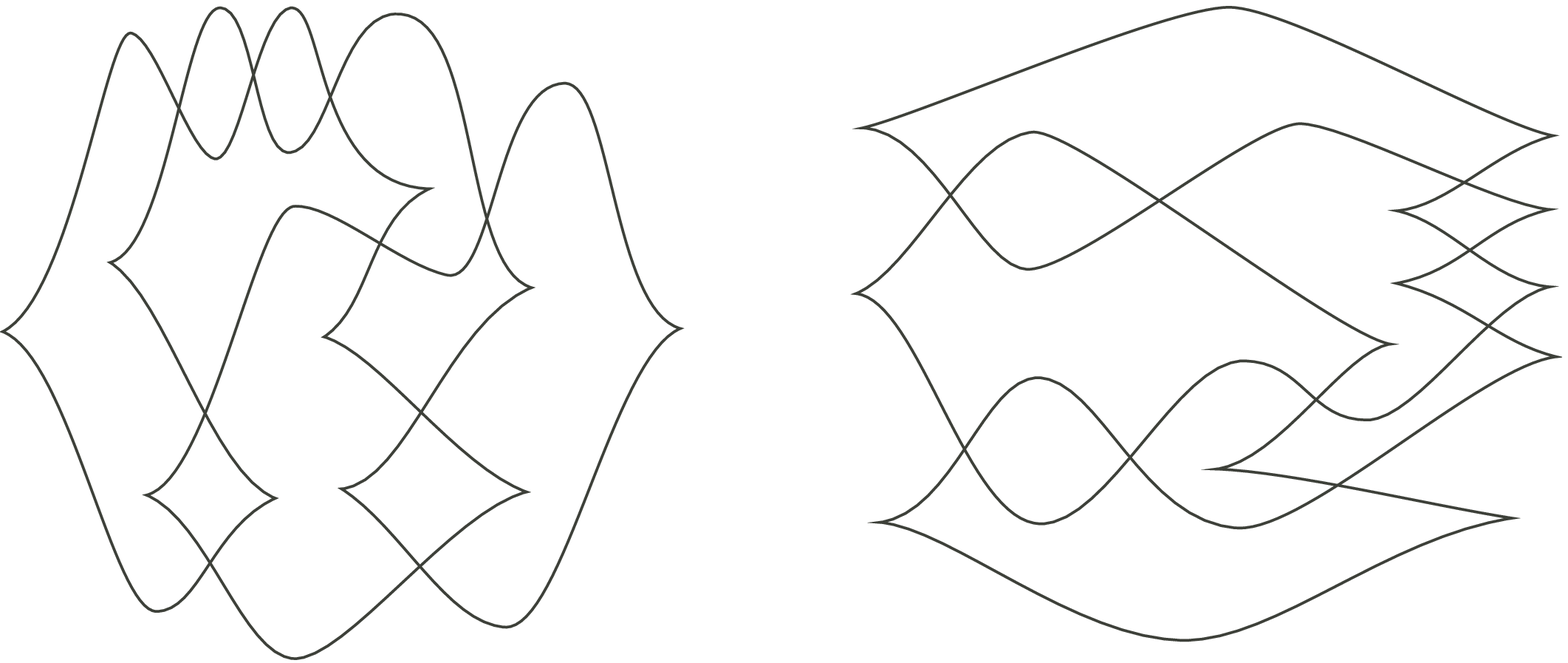}{.14}
&\up{2t}		&\up{(-5,0)}	\\ [1ex]
\up{9_{20}}	&\up{(-12,1)}	&			&\fig{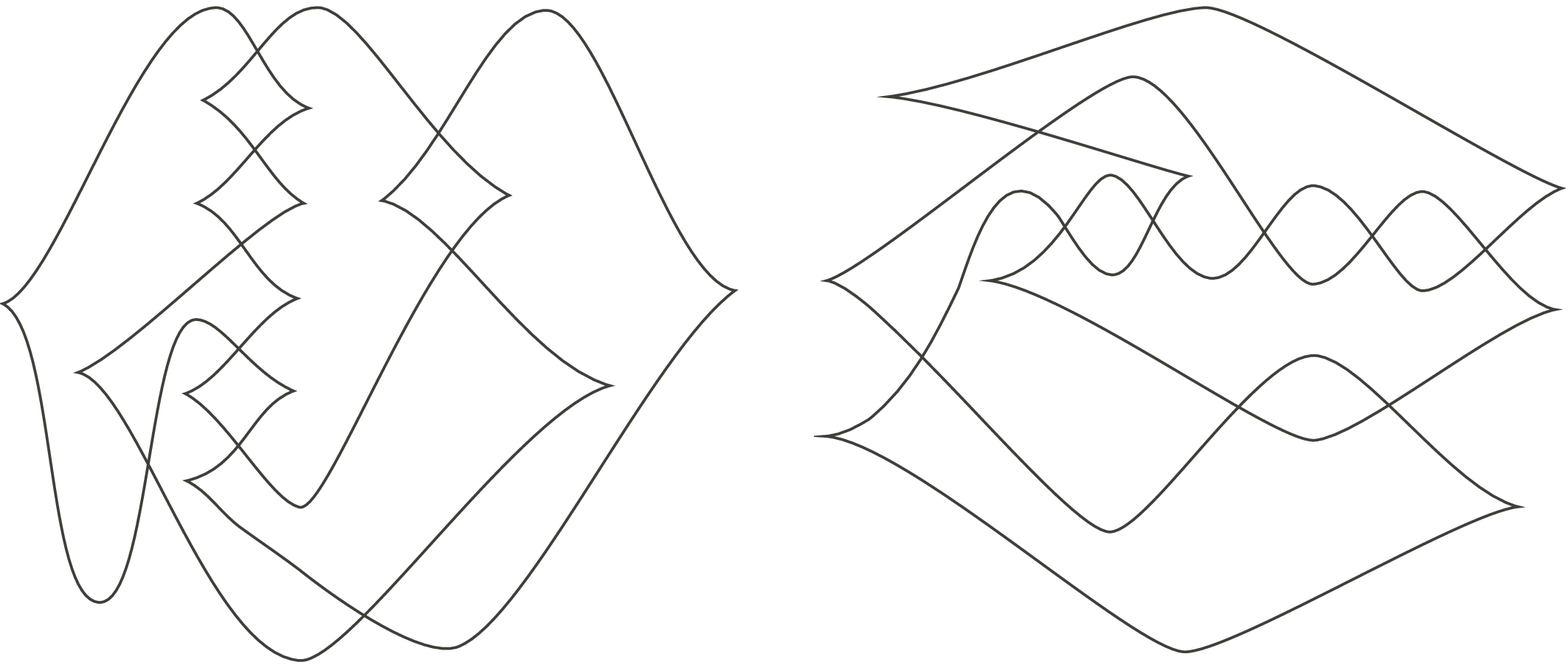}{.14}
&\up{2+t}		&\up{(1,0)}	\\ [1ex]
\up{9_{21}}	&\up{(-1,0)}	&			&\fig{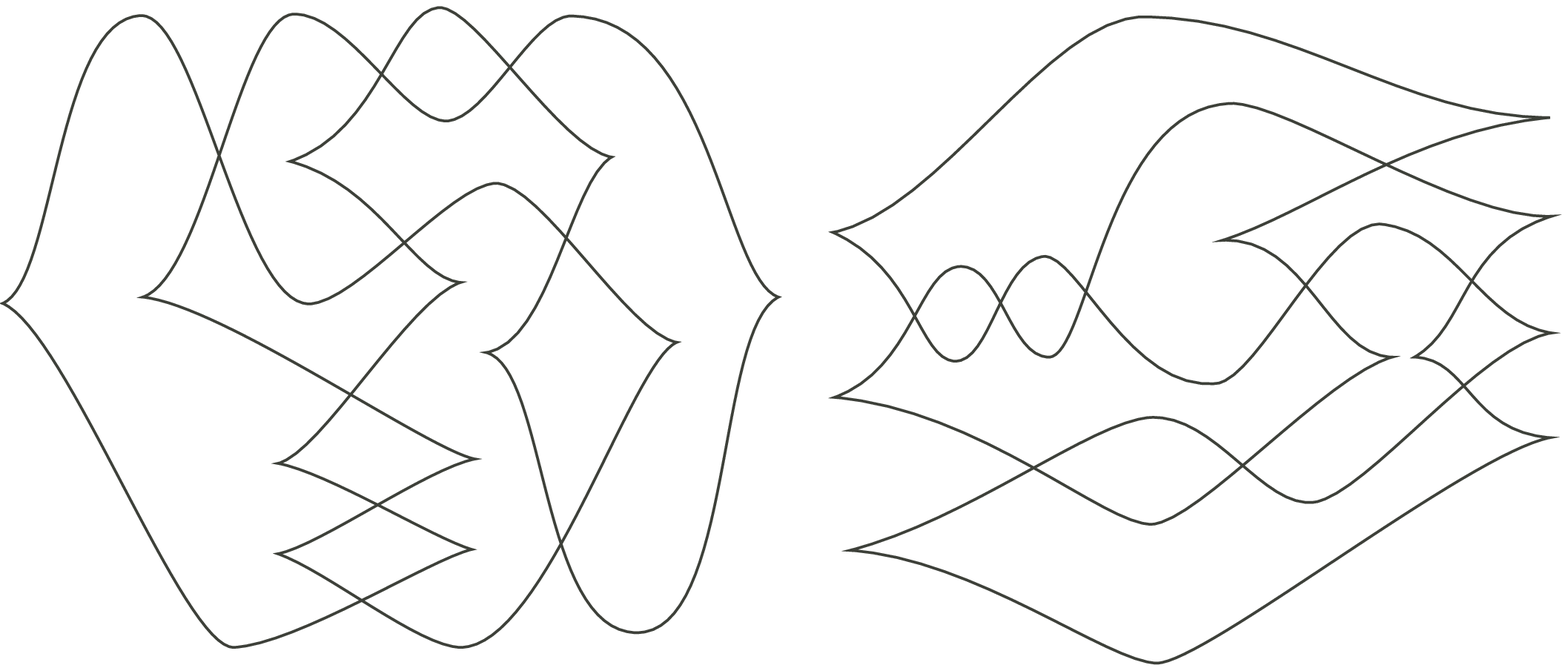}{.14}
&			&\up{(-10,1)}	\\ [1ex]
\up{9_{22}}	&\up{(-3,0)}	&\up{1+2t}	&\fig{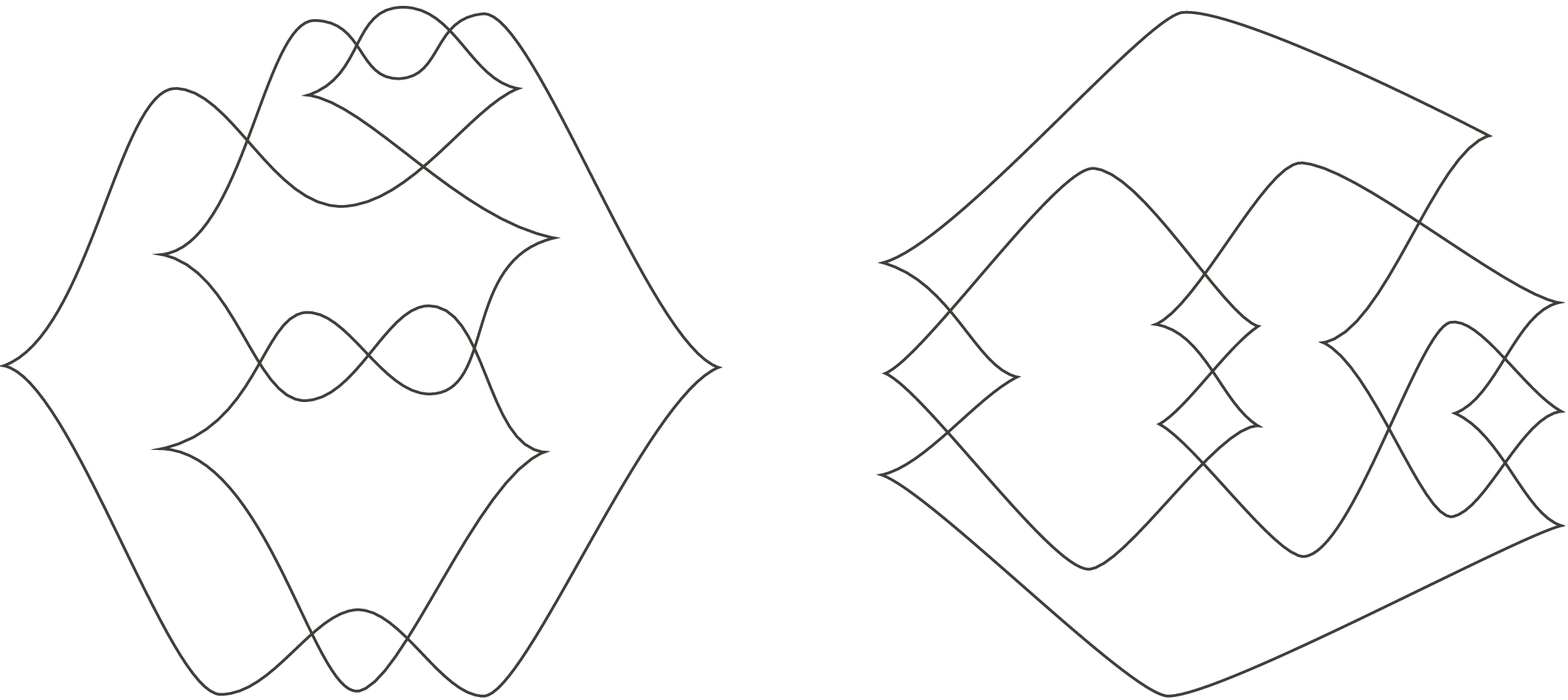}{.14}	&
&\up{(-8,1)}	\\ [1ex]
\up{9_{23}}	&\up{(-14,1)}	&			&\fig{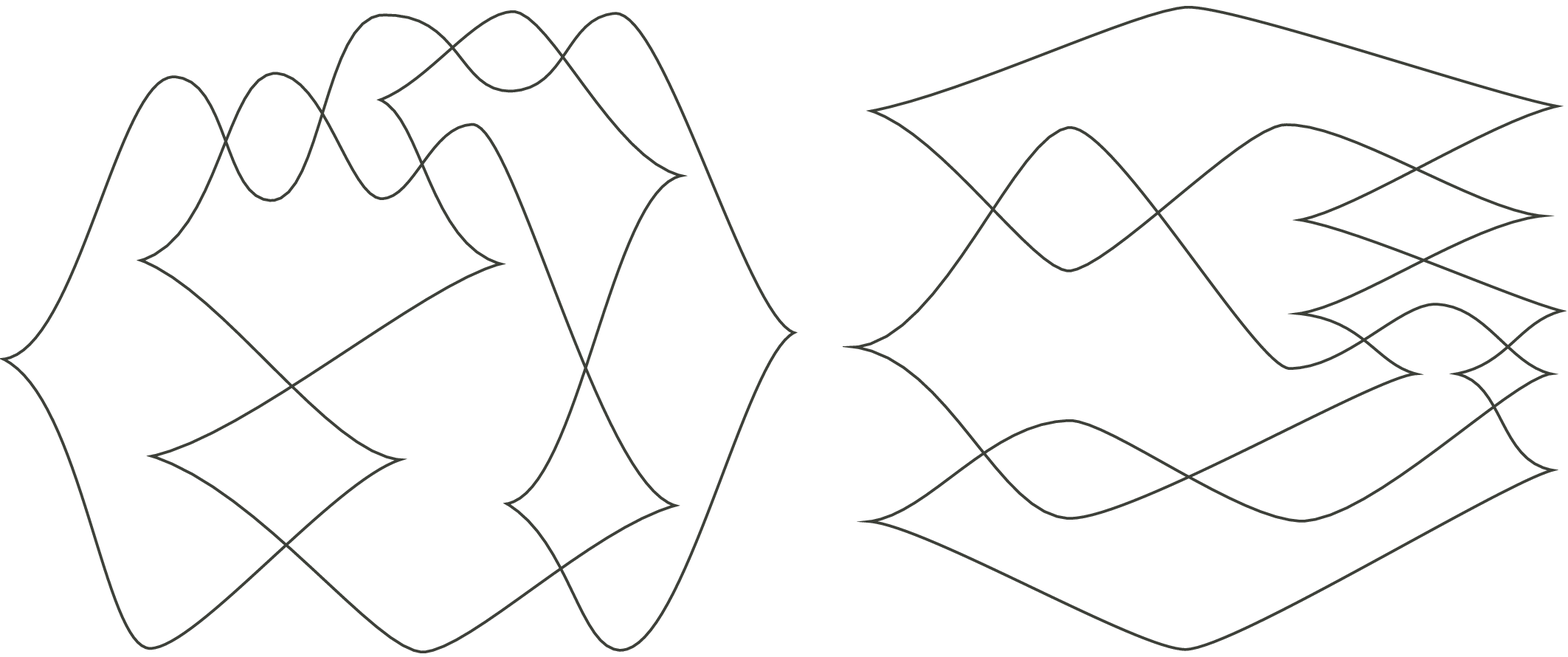}{.14}
&\up{1+t^2}	&\up{(3,0)}	\\ [1ex]
\up{9_{24}}	&\up{(-6,1)}	&			&\fig{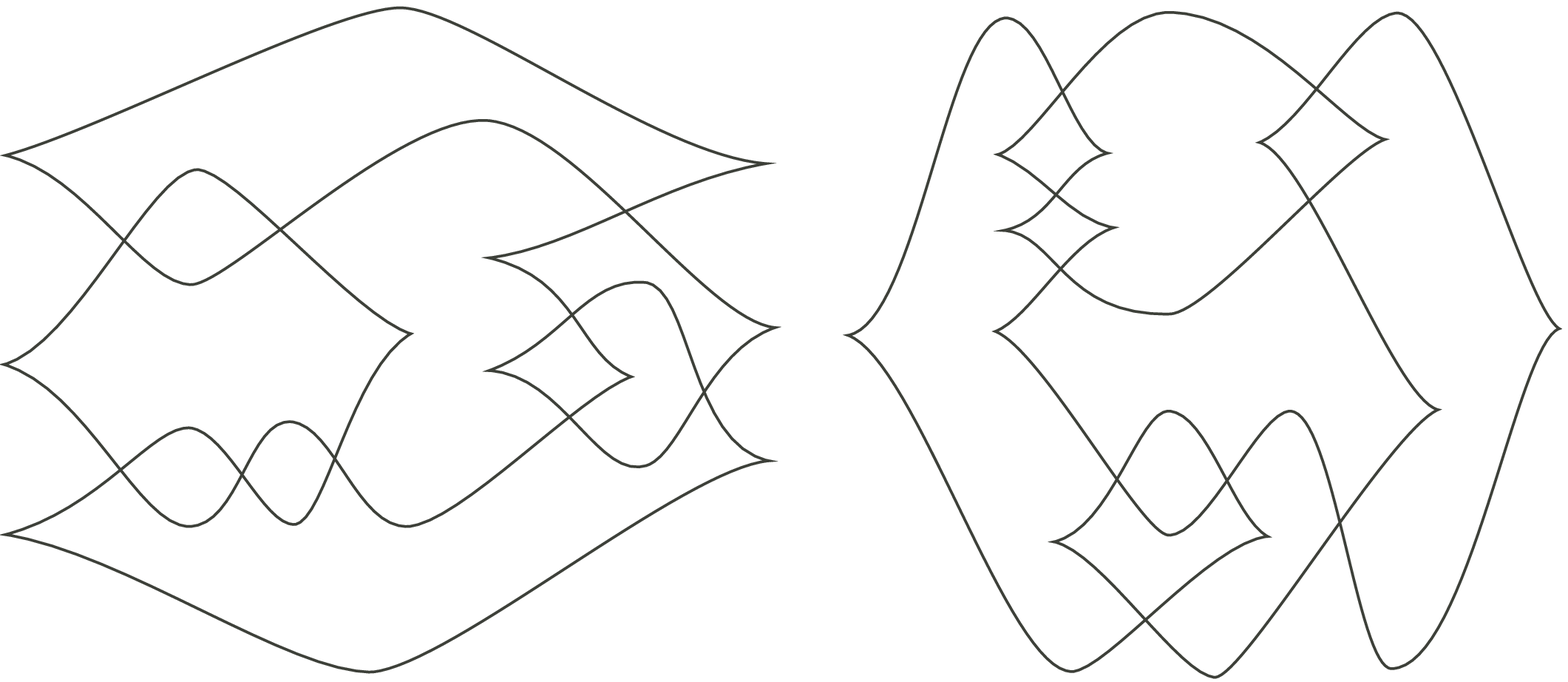}{.14}
&			&\up{(-5,2)}	\\ [1ex]
\hline
\end{tabular}
\end{table}

\clearpage

\begin{table}[!ht]
\tabletop
\up{9_{25}}	&\up{(-10,1)}	&			&\fig{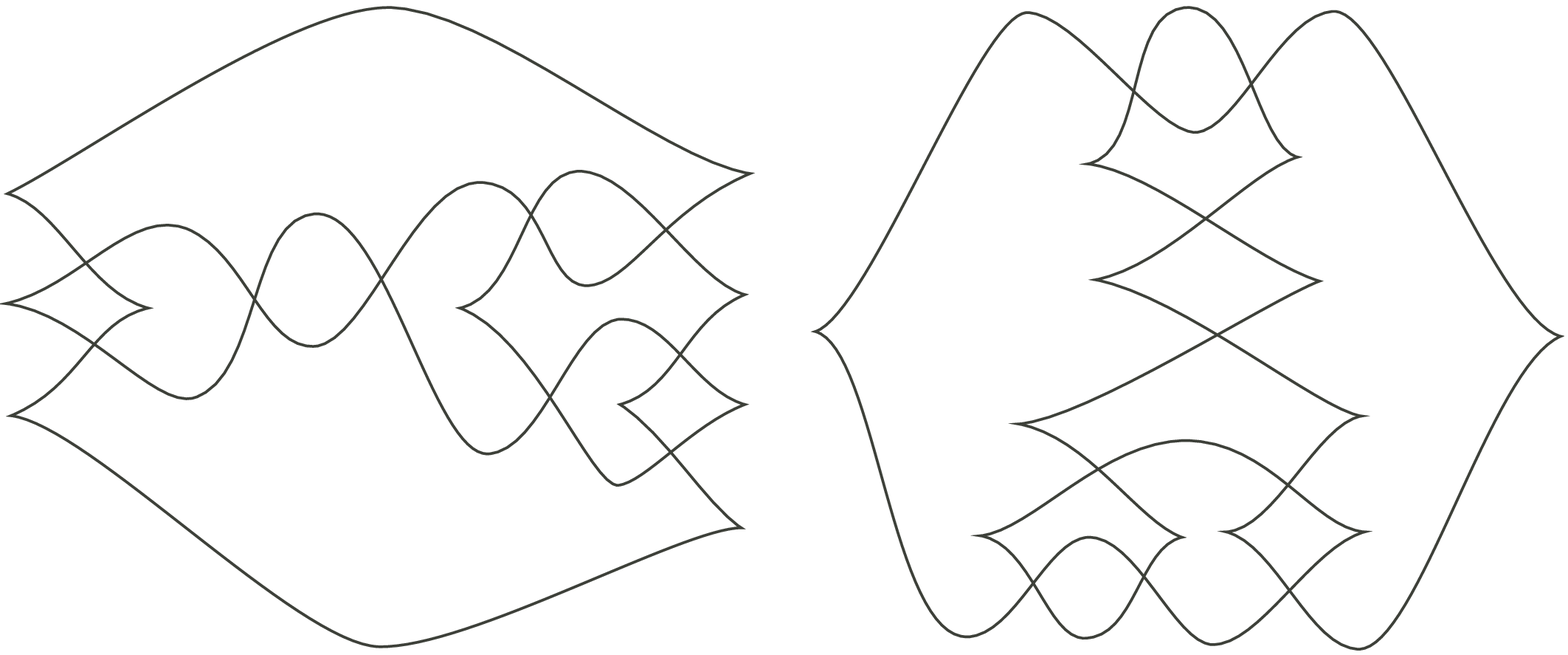}{.14}
&\up{1+t^3}	&\up{(-1,0)}	\\ [1ex]
\up{9_{26}}	&\up{(-2,1)}	&			&\fig{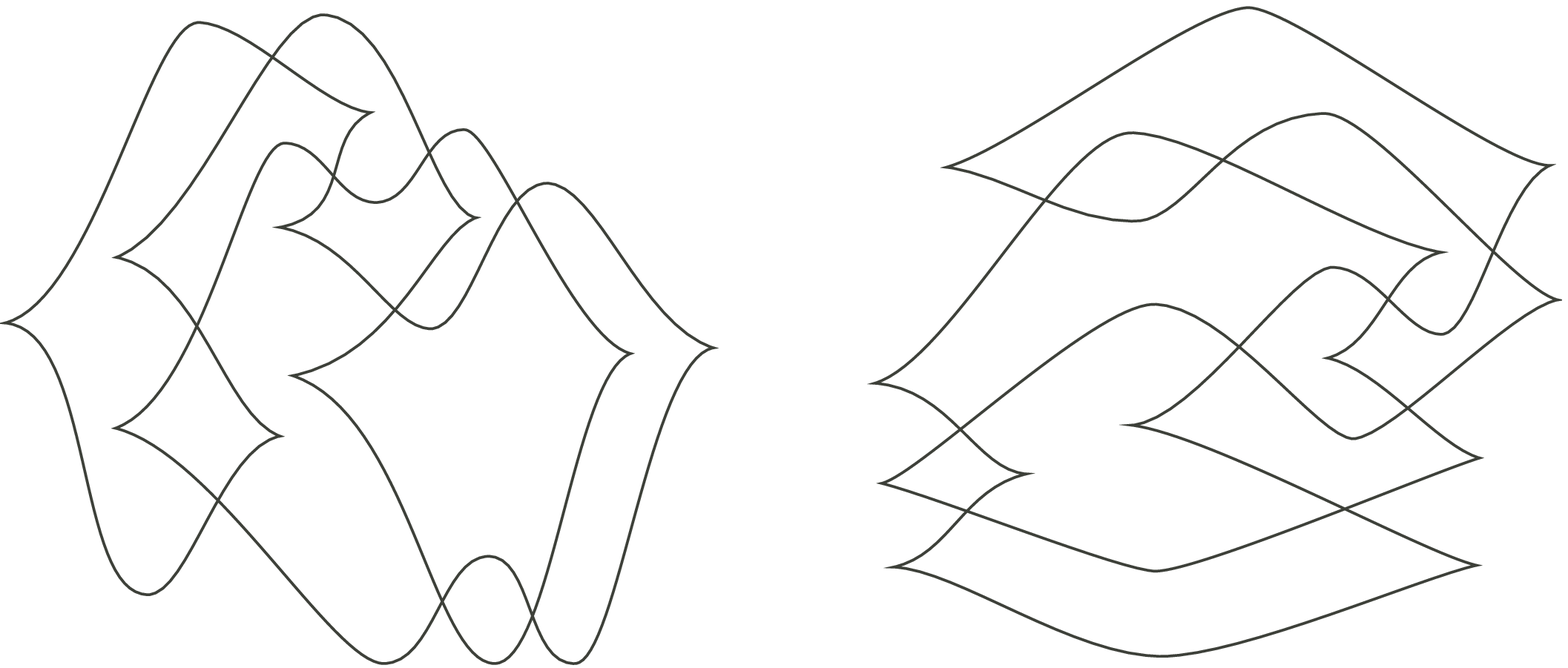}{.14}
&			&\up{(-9,0)}	\\ [1ex]
\up{9_{27}}	&\up{(-6,1)}	&			&\fig{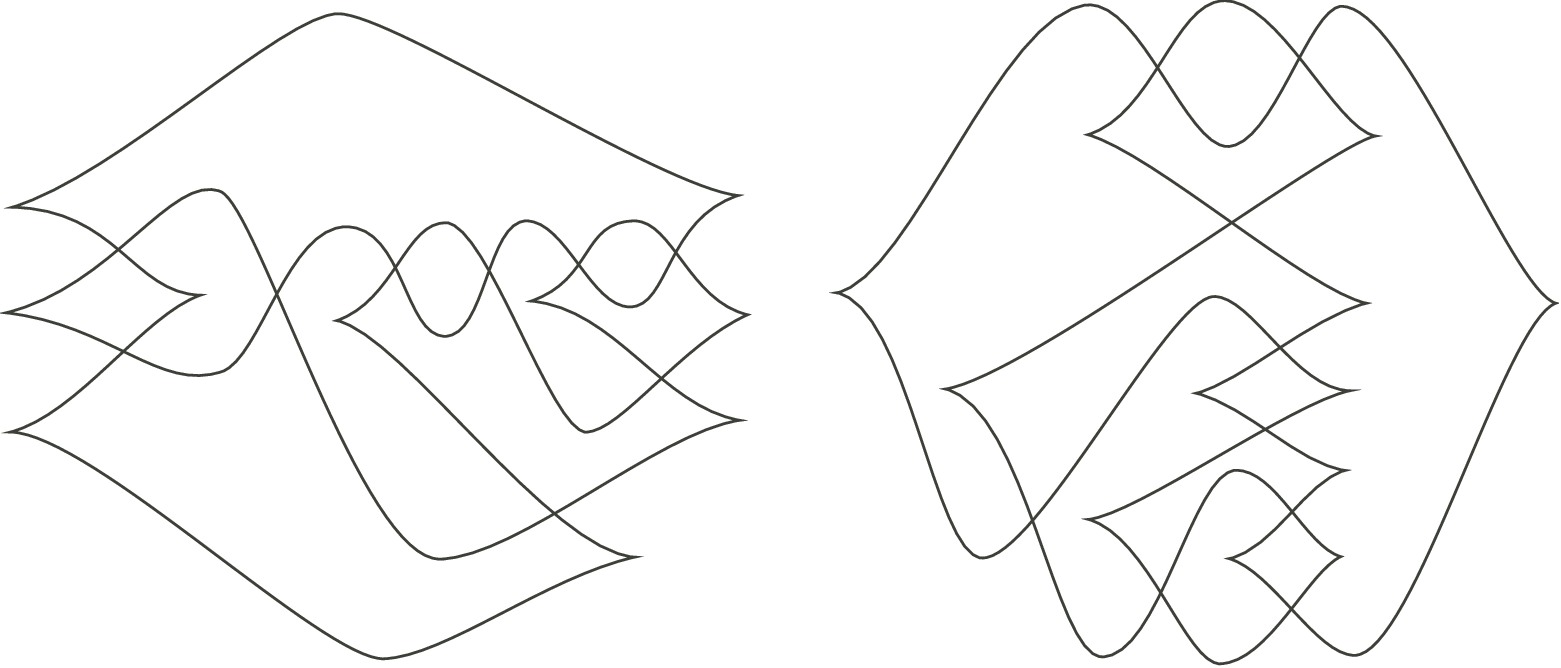}{.14}
&			&\up{(-5,2)}	\\ [1ex]
\up{9_{28}}	&\up{(-9,0)}	&			&\fig{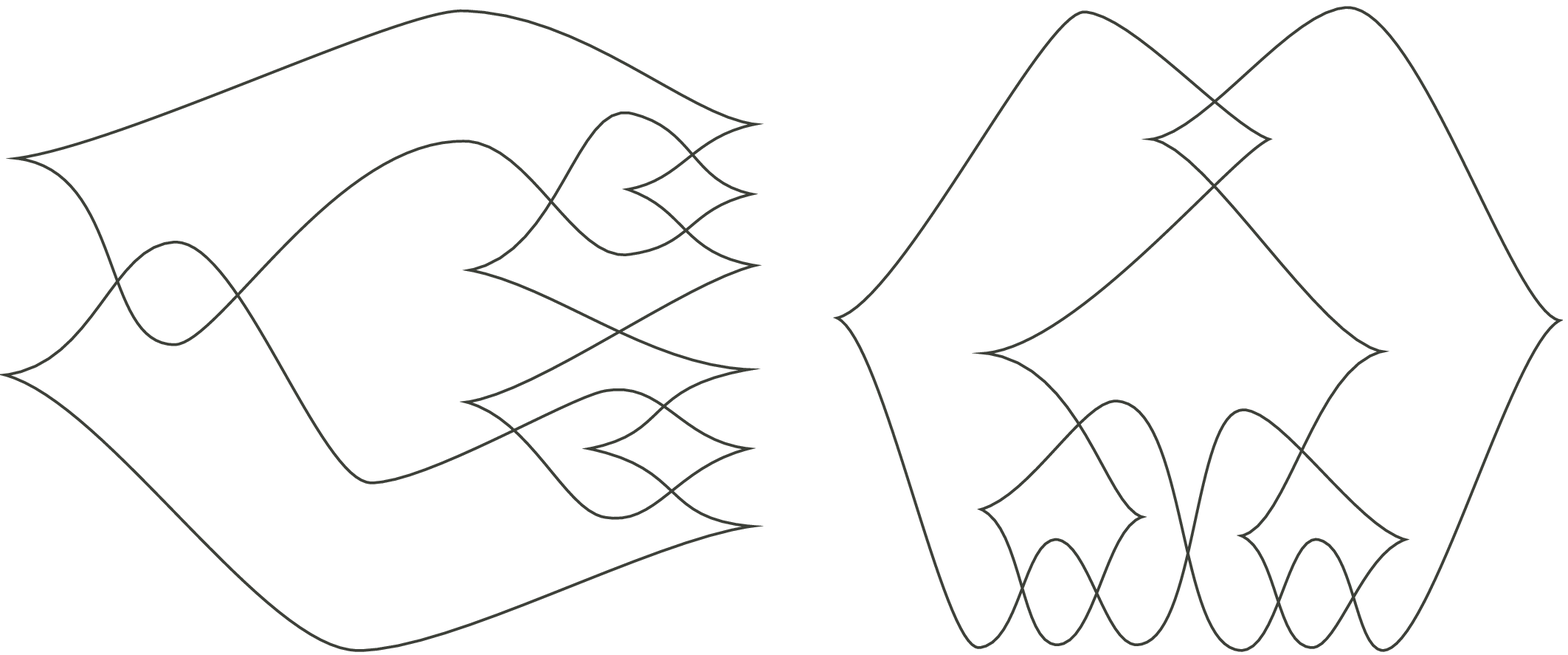}{.14}
&			&\up{(-2,1)}	\\ [1ex]
\up{9_{29}}	&\up{(-8,3)}	&			&\fig{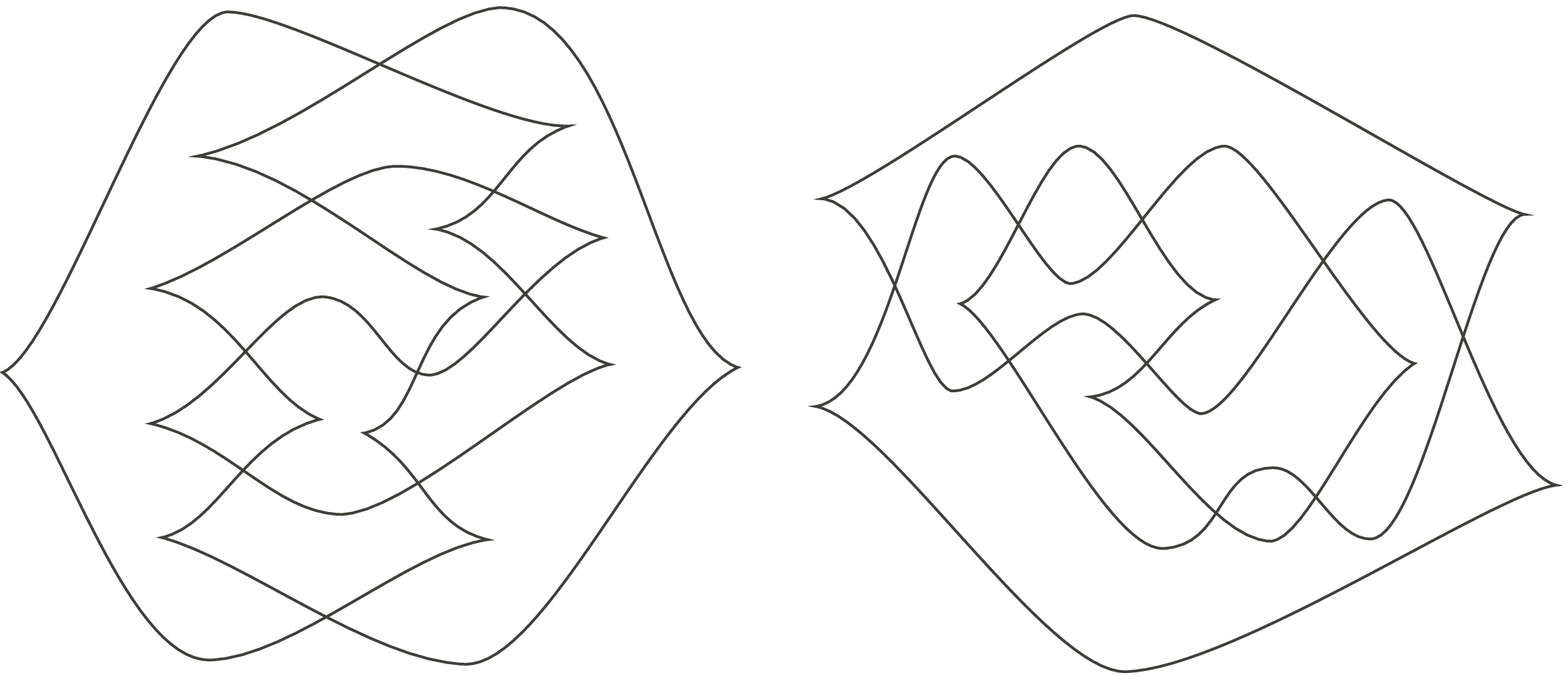}{.14}
&\up{1+2t}	&\up{(-3,0)}	\\ [1ex]
\up{9_{30}}	&\up{(-6,1)}	&			&\fig{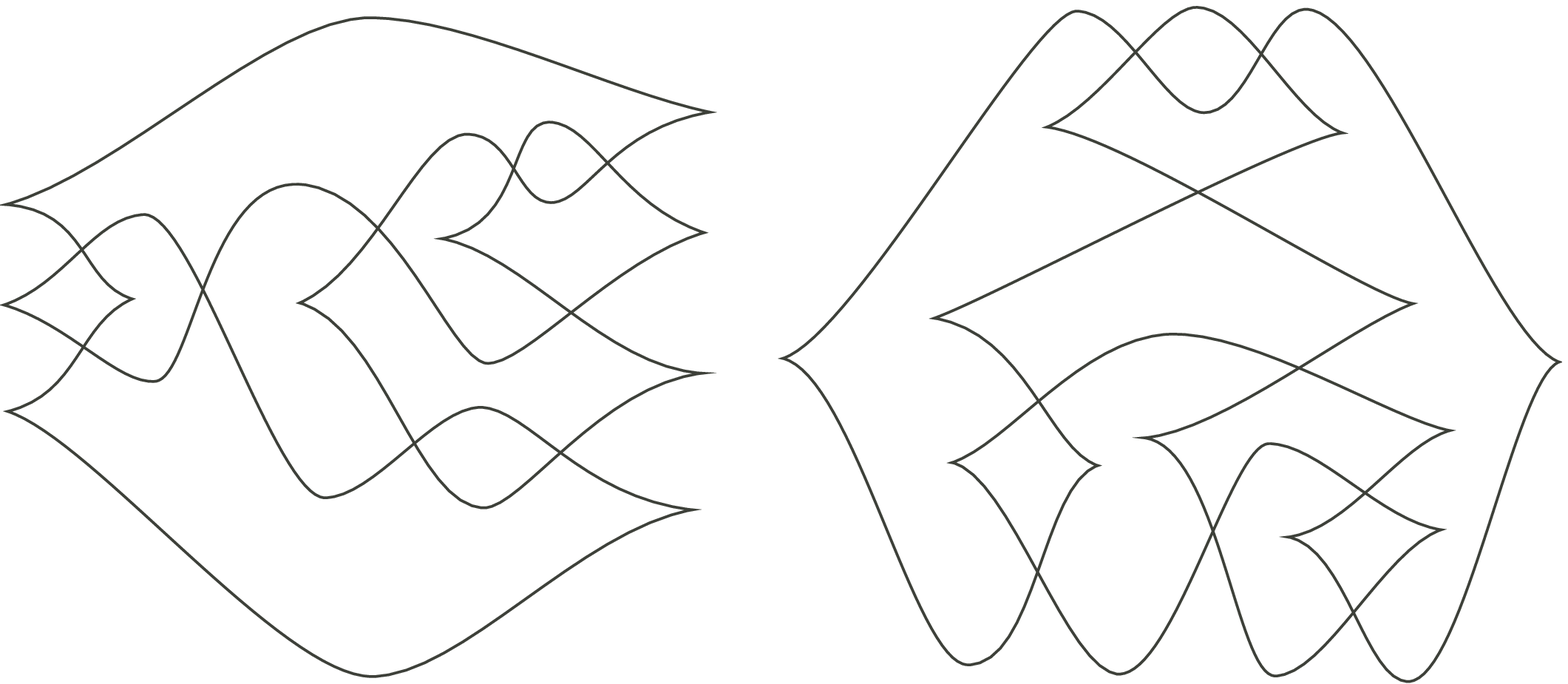}{.14}
&			&\up{(-5,0)}	\\ [1ex]
\up{9_{31}}	&\up{(-9,2)}	&			&\fig{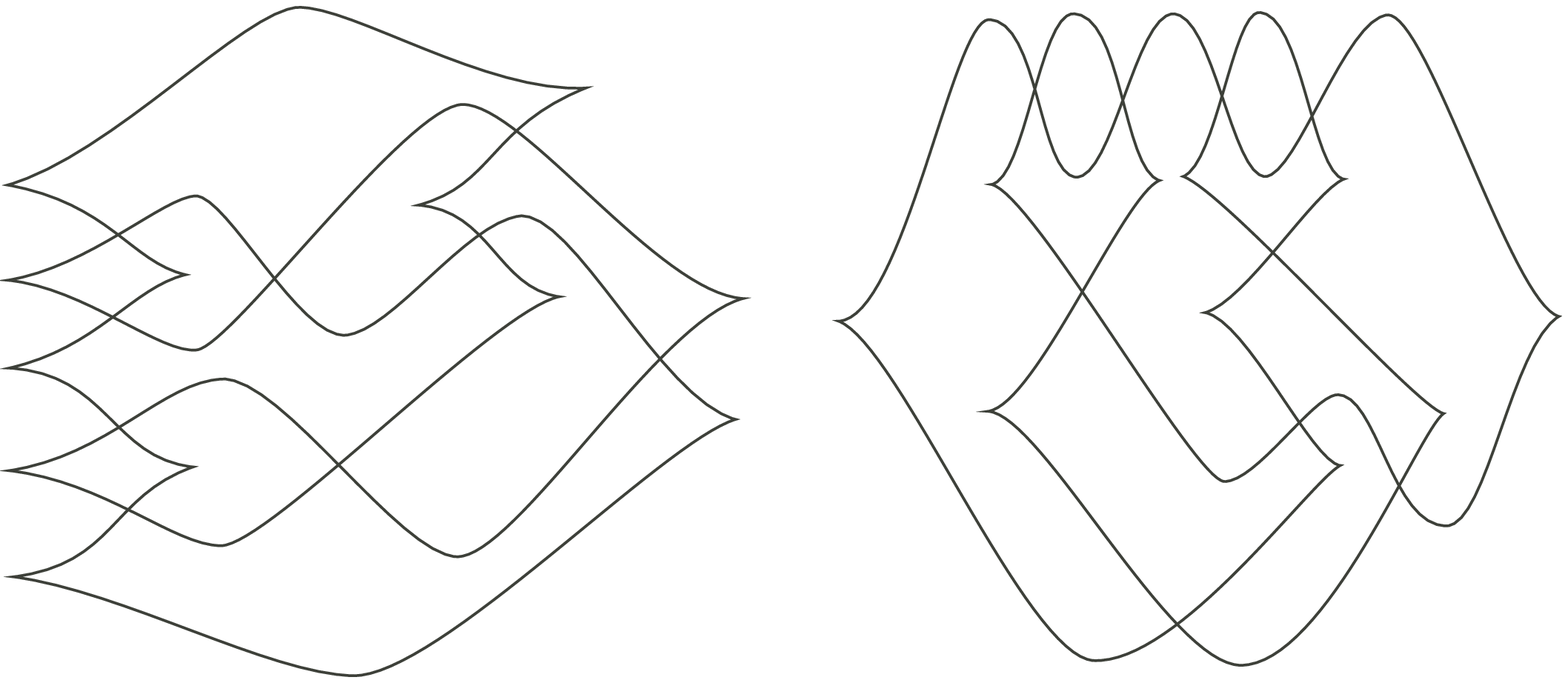}{.14}
&			&\up{(-2,1)}	\\ [1ex]
\up{9_{32}}	&\up{(-2,1)}	&			&\fig{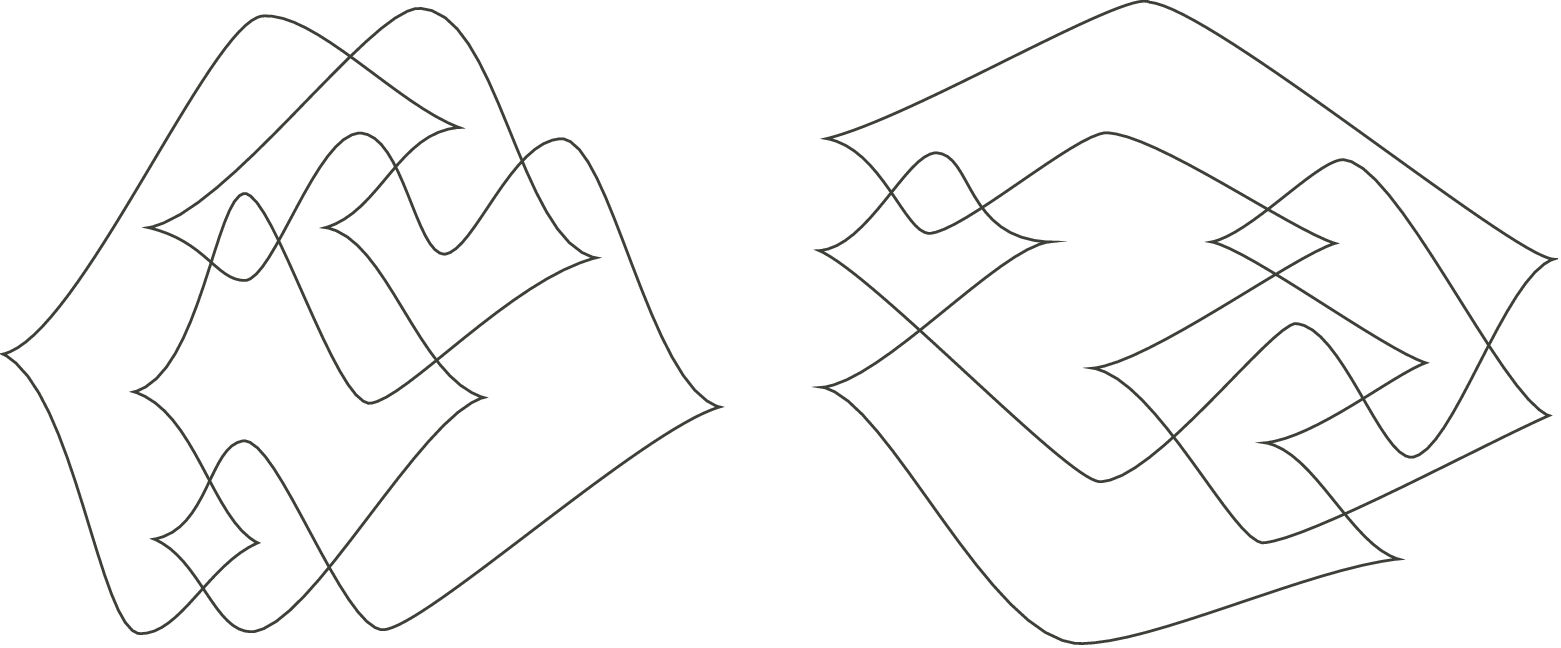}{.14}
&			&\up{(-9,2)}	\\ [1ex]
\up{9_{33}}	&\up{(-6,1)}	&			&\fig{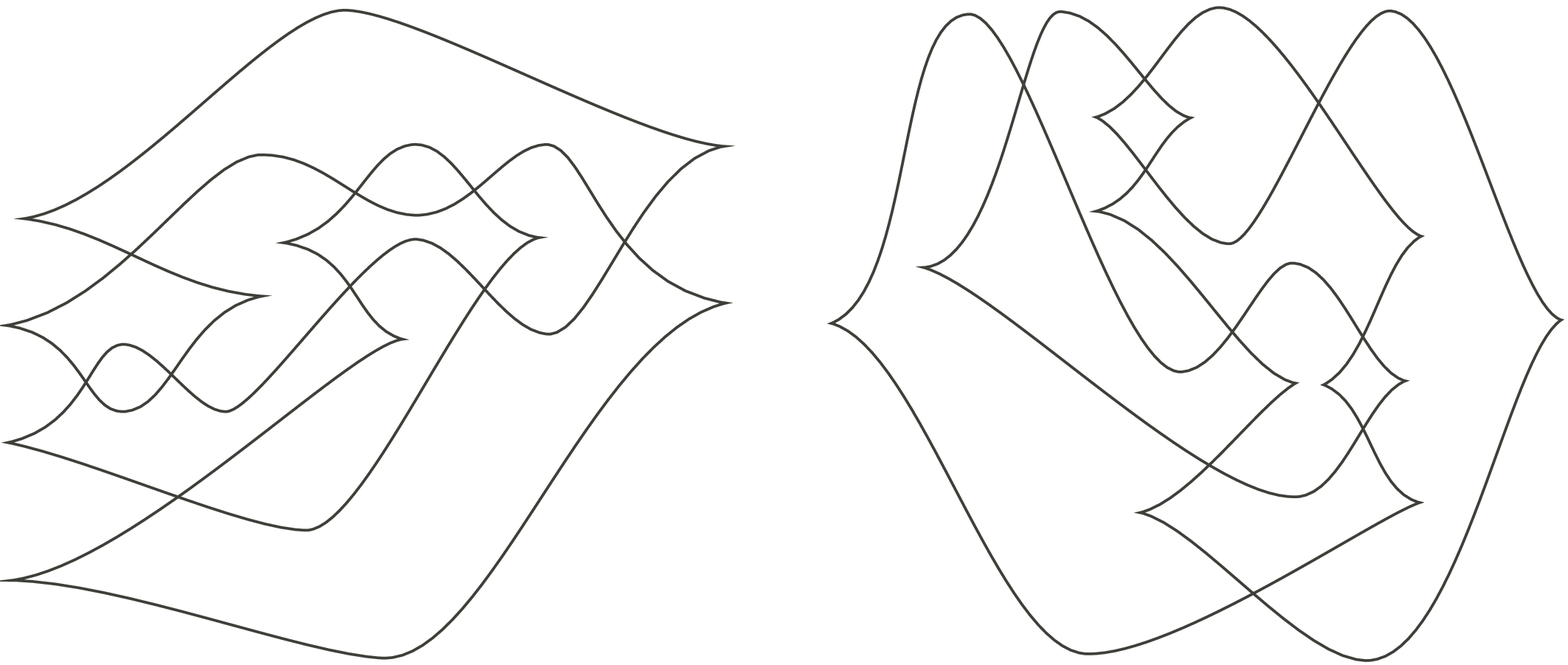}{.14}
&			&\up{(-5,2)}	\\ [1ex]
\up{9_{34}}	&\up{(-6,1)}	&			&\fig{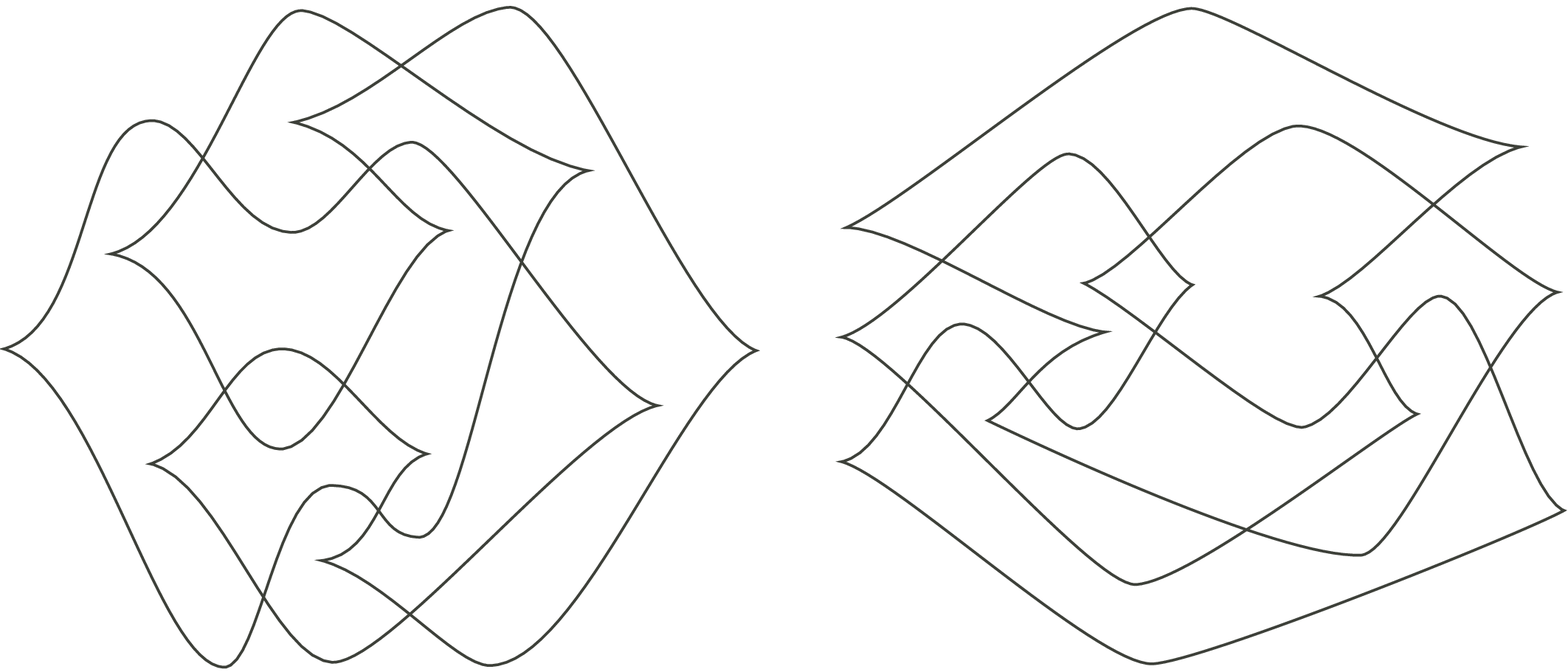}{.14}
&			&\up{(-5,0)}	\\ [1ex]
\hline
\end{tabular}
\end{table}

\clearpage

\begin{table}[!ht]
\tabletop
\up{9_{35}}	&\up{(-12,1)}	&			&\fig{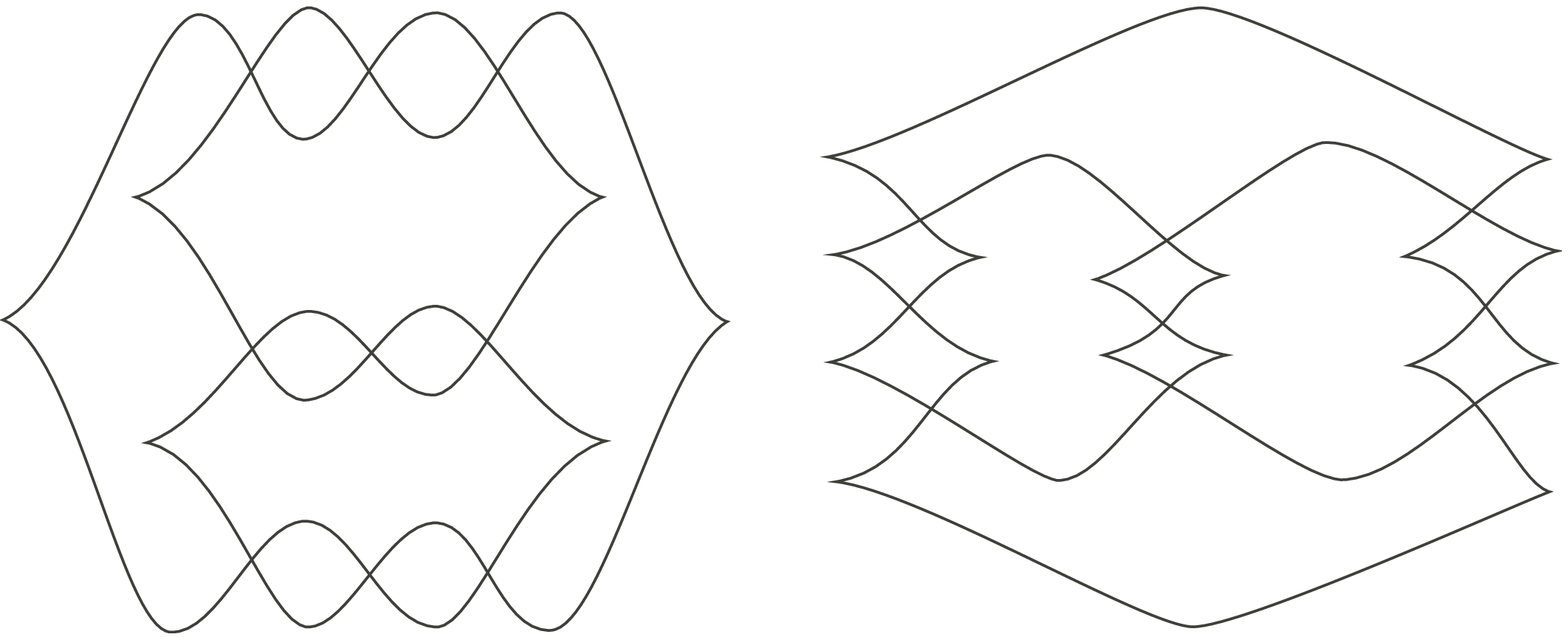}{.14}
&\up{1}		&\up{(1,0)}	\\ [1ex]
\up{9_{36}}	&\up{(1,0)}	&\up{1+t+t^2}	&\fig{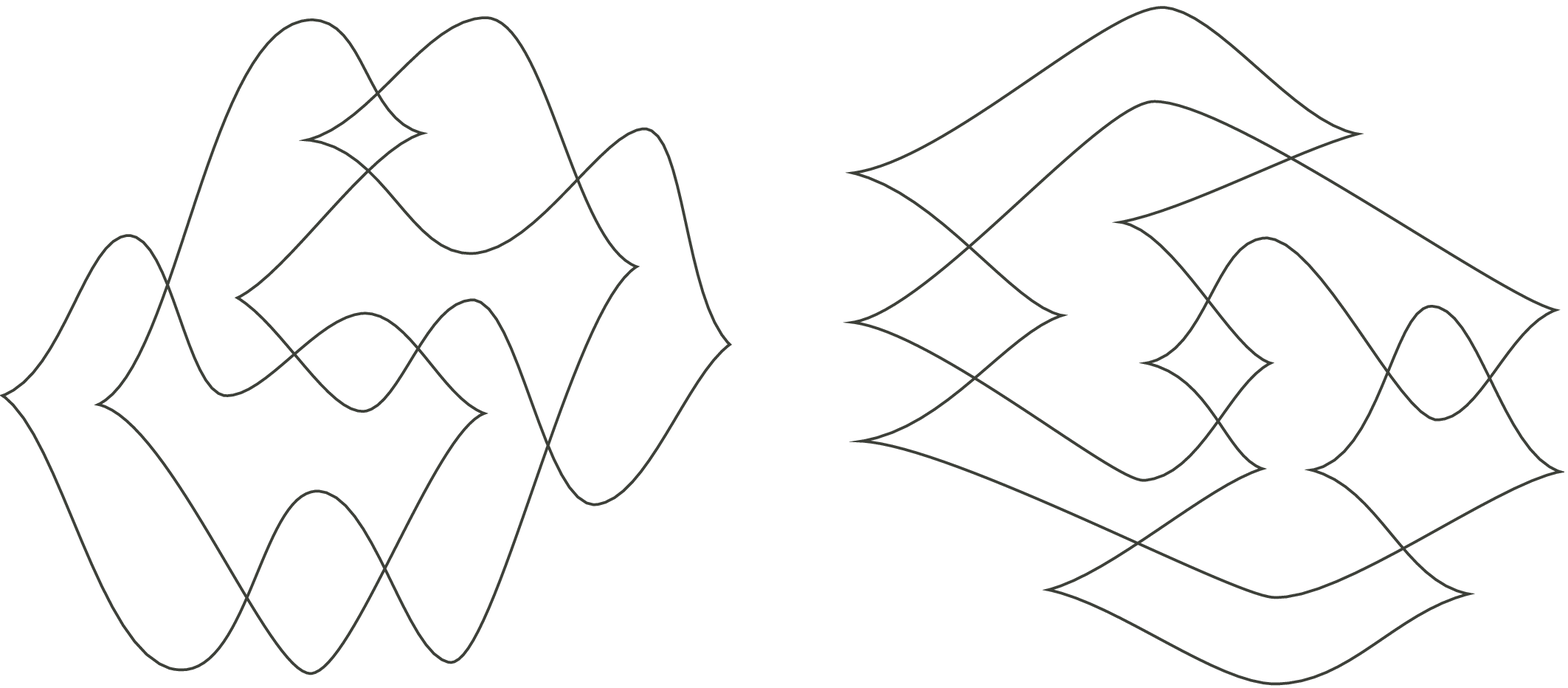}{.14}	&
&\up{(-12,1)}	\\ [1ex]
\up{9_{37}}	&\up{(-6,1)}	&			&\fig{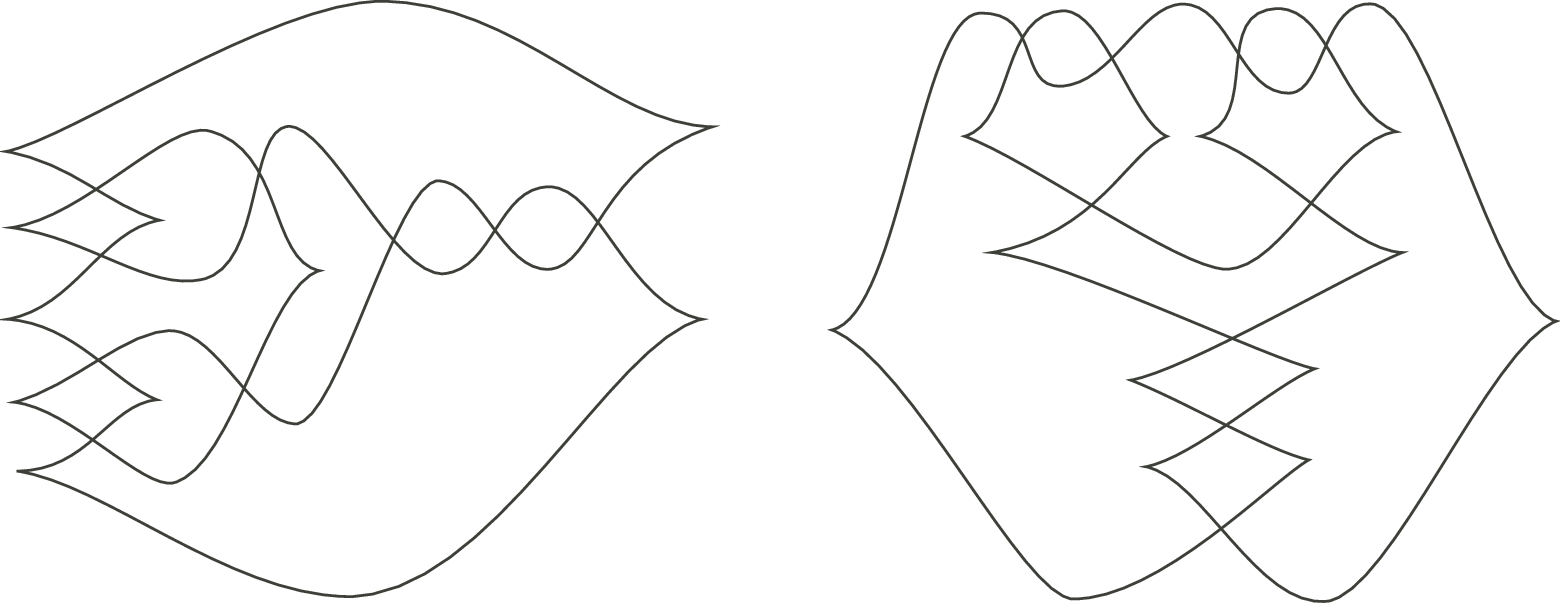}{.14}
&\up{2t^3}		&\up{(-5,0)}	\\ [1ex]
\up{9_{38}}	&\up{(-14,1)}	&			&\fig{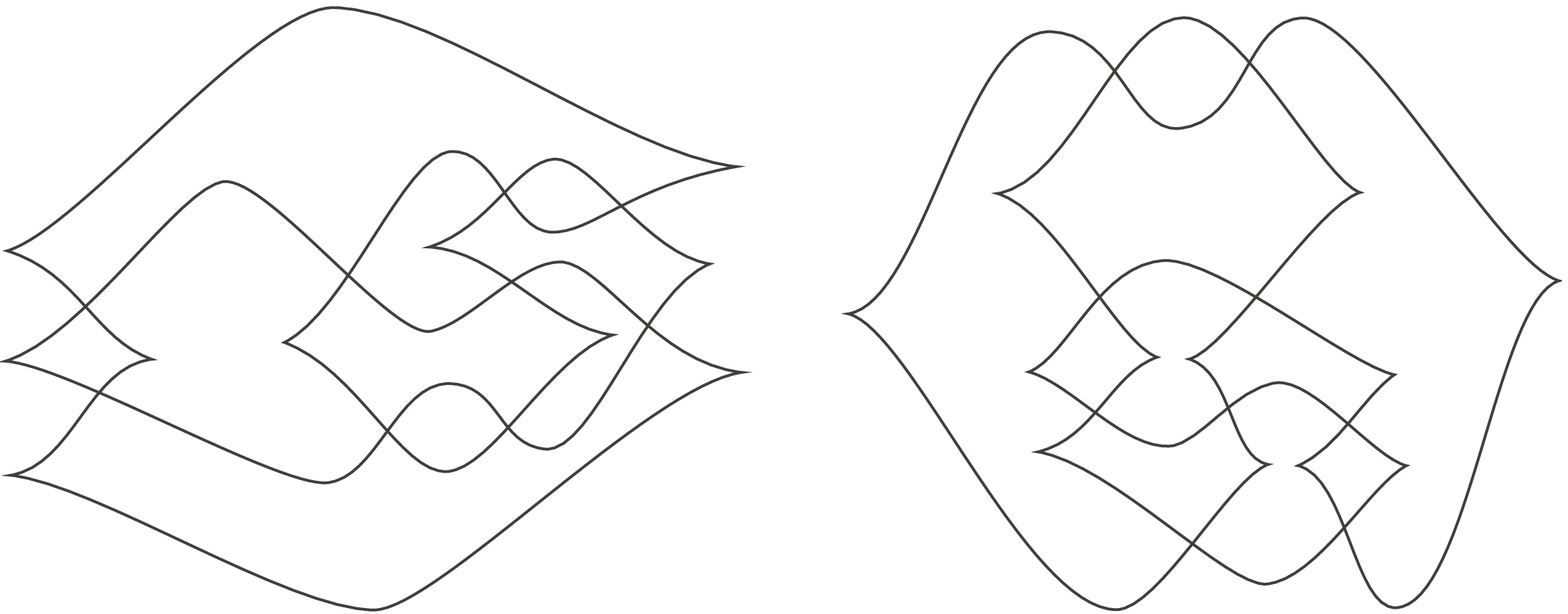}{.14}
&\up{1+t^2}	&\up{(3,0)}	\\ [1ex]
\up{9_{39}}	&\up{(-1,0)}	&\up{1+t}		&\fig{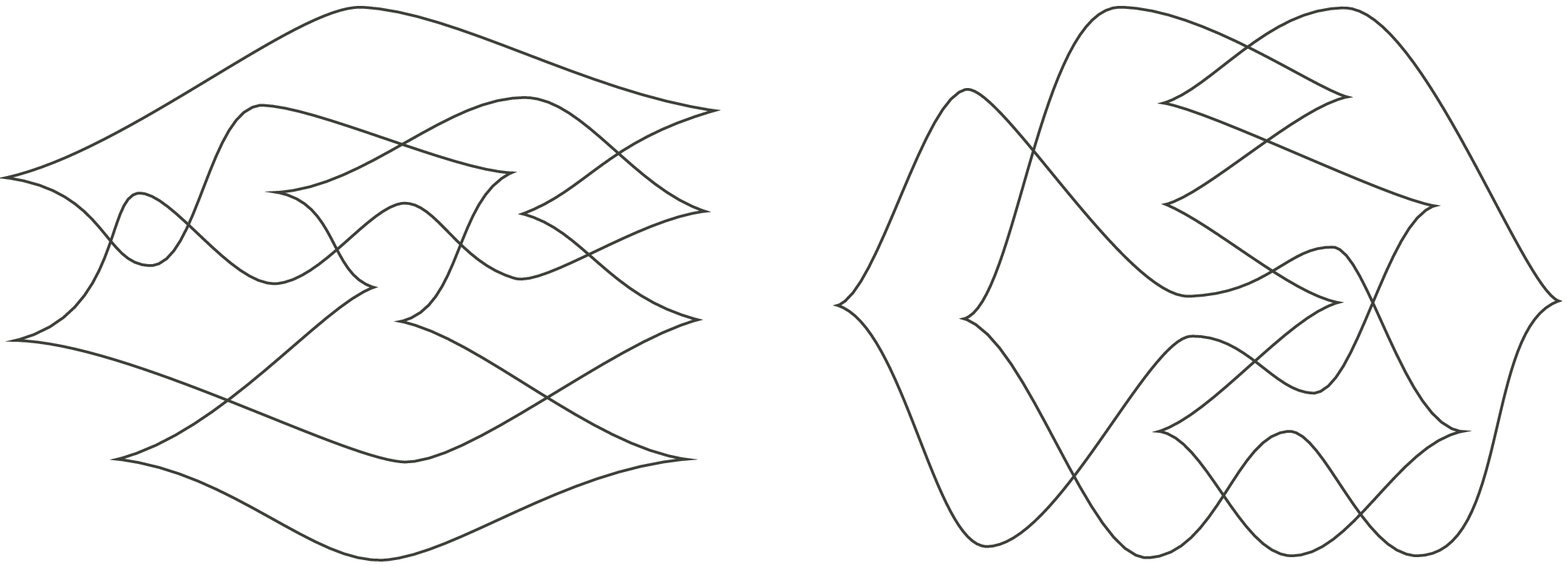}{.14}
&			&\up{(-10,1)}	\\ [1ex]
\up{9_{40}}	&\up{(-9,2)}	&			&\fig{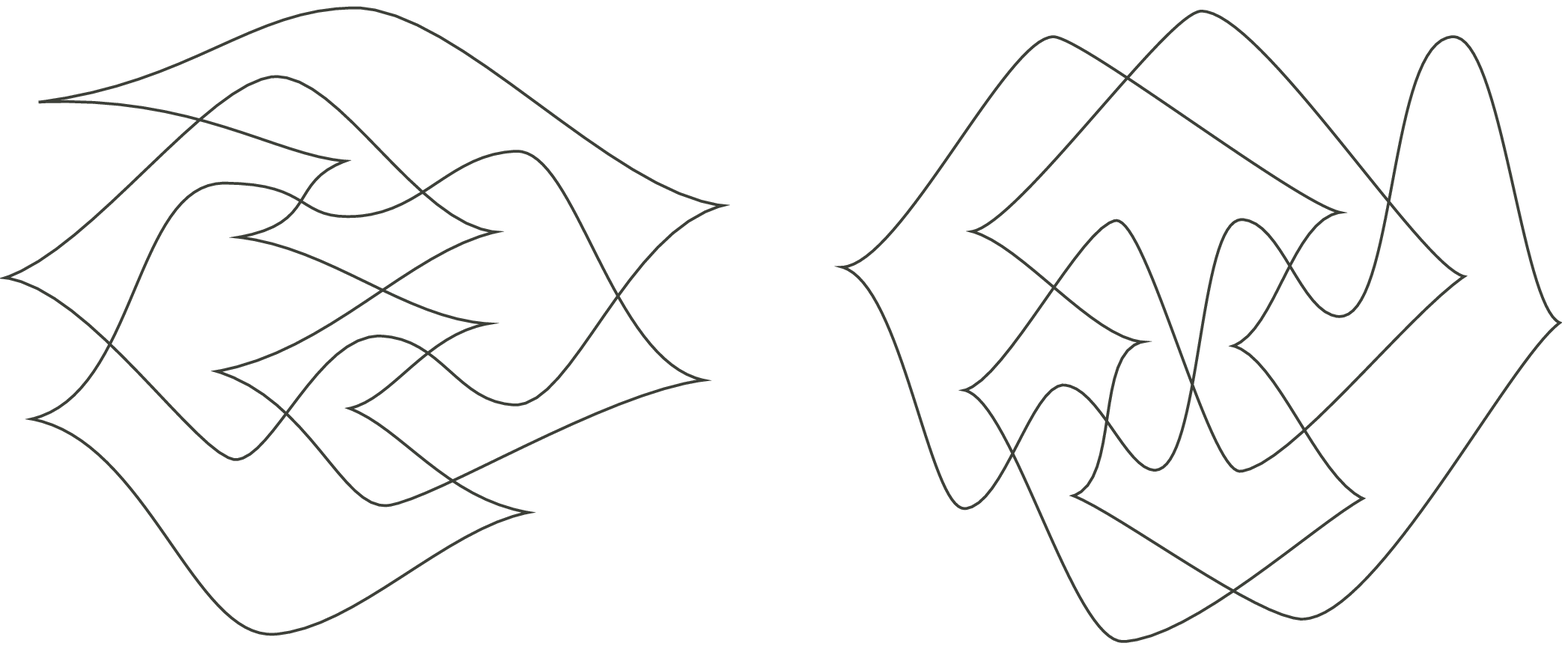}{.14}
&			&\up{(-2,1)}	\\ [1ex]
\up{9_{41}}	&\up{(-7,0)}	&\up{3t}		&\fig{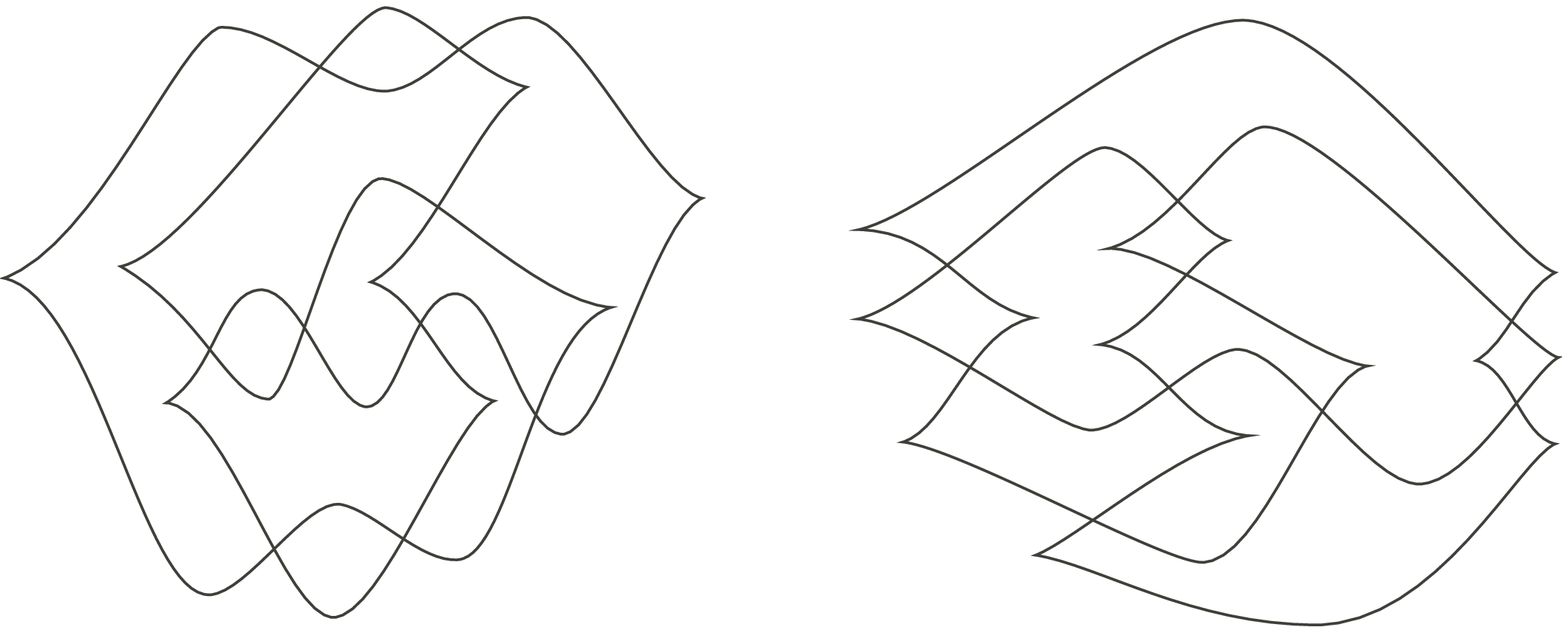}{.14}
&			&\up{(-4,1)}	\\ [1ex]
\up{9_{42}}	&\up{(-3,0)}	&\up{1+2t}	&\fig{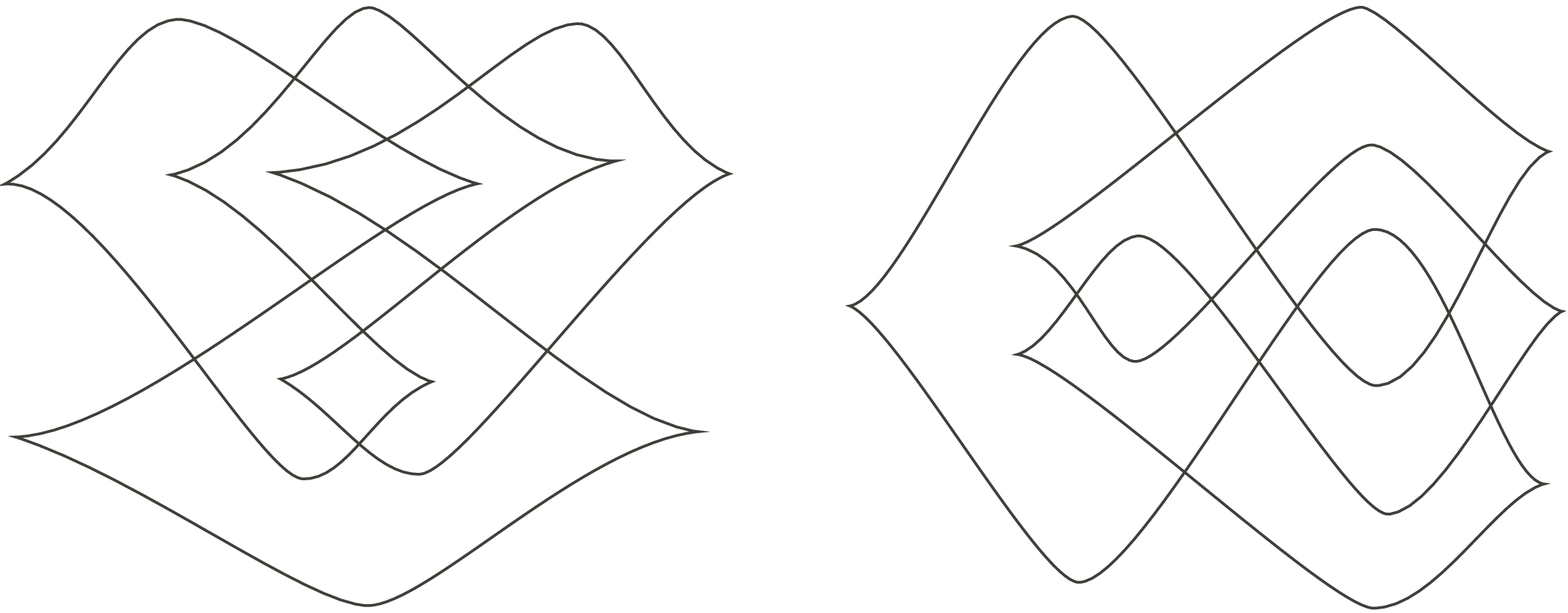}{.14}	&
&\up{(-5,0)}	\\ [1ex]
\up{9_{43}}	&\up{(1,0)}	&\up{2+t}		&\fig{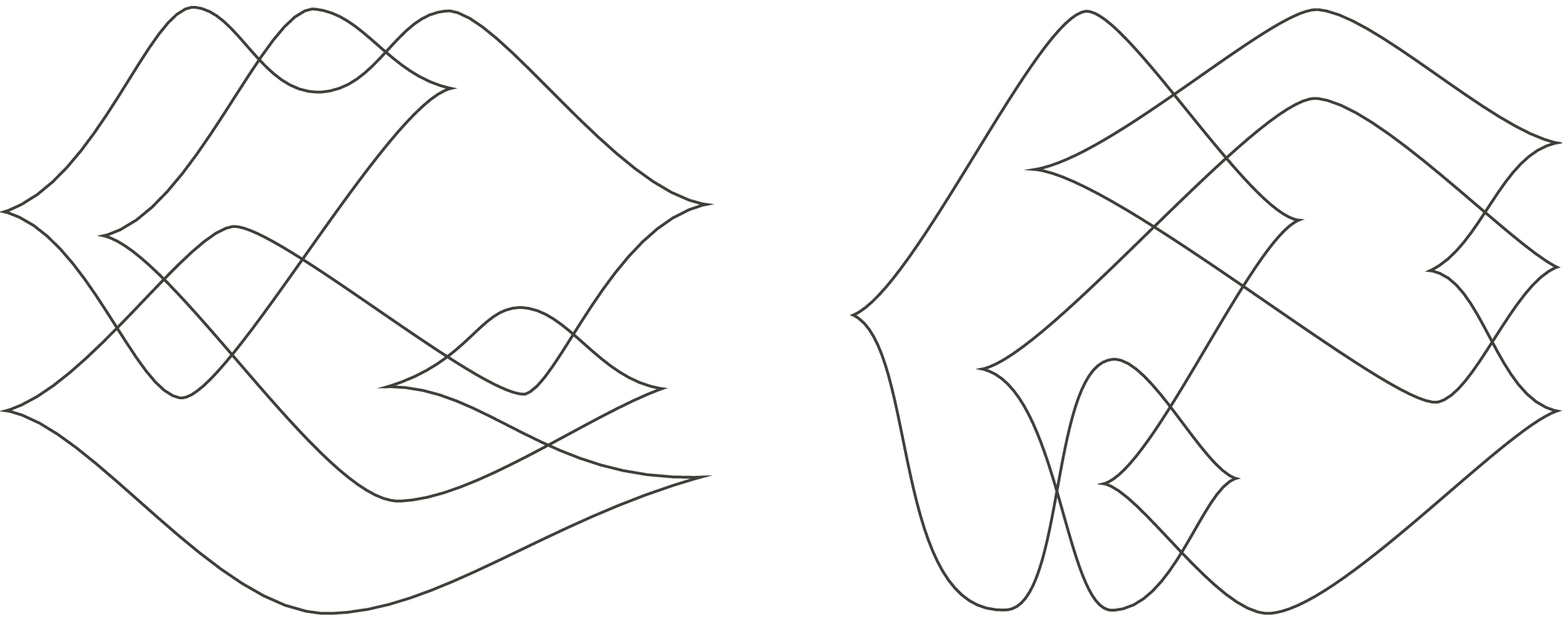}{.14}
&			&\up{(-10,1)}	\\ [1ex]
\up{9_{44}}	&\up{(-6,1)}	&			&\fig{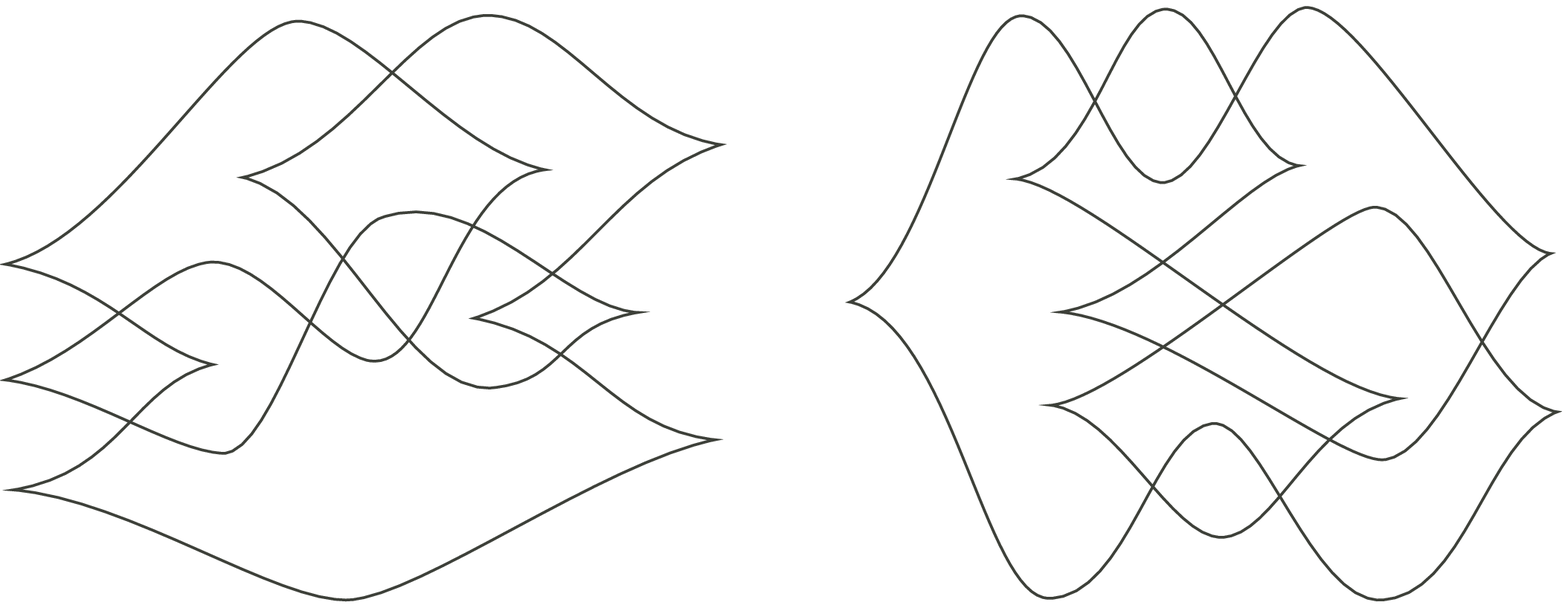}{.14}
&\up{t}		&\up{(-3,0)}	\\ [1ex]
\up{9_{45}} &\up{(-10,1)}	& &\fig{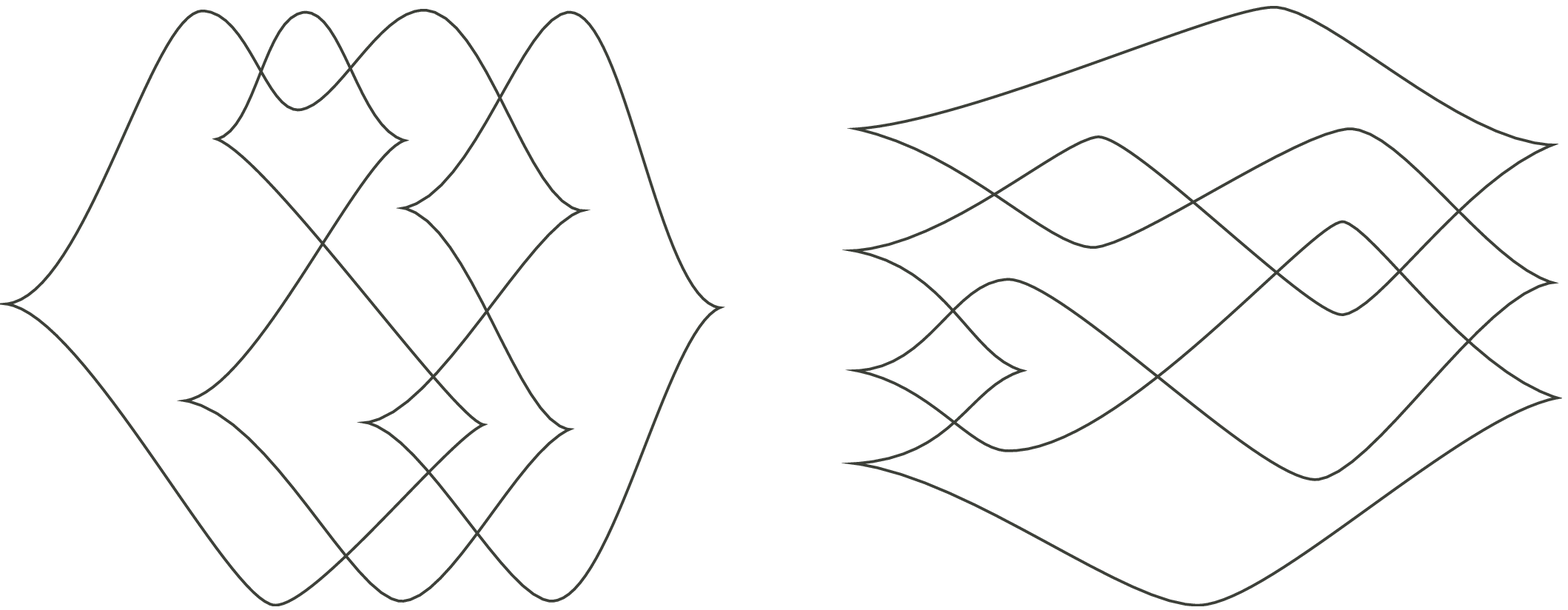}{.14}
&\up{\underset{\textstyle 1+t+t^2}{\overset{\textstyle 1}{}}}
															&\up{(1,0)}
															\\
															[1ex]
\hline
\end{tabular}
\end{table}

\clearpage

\begin{table}[!ht]
\tabletop
\up{9_{46}}	&\up{(-1,0)}	&\up{0}		&\fig{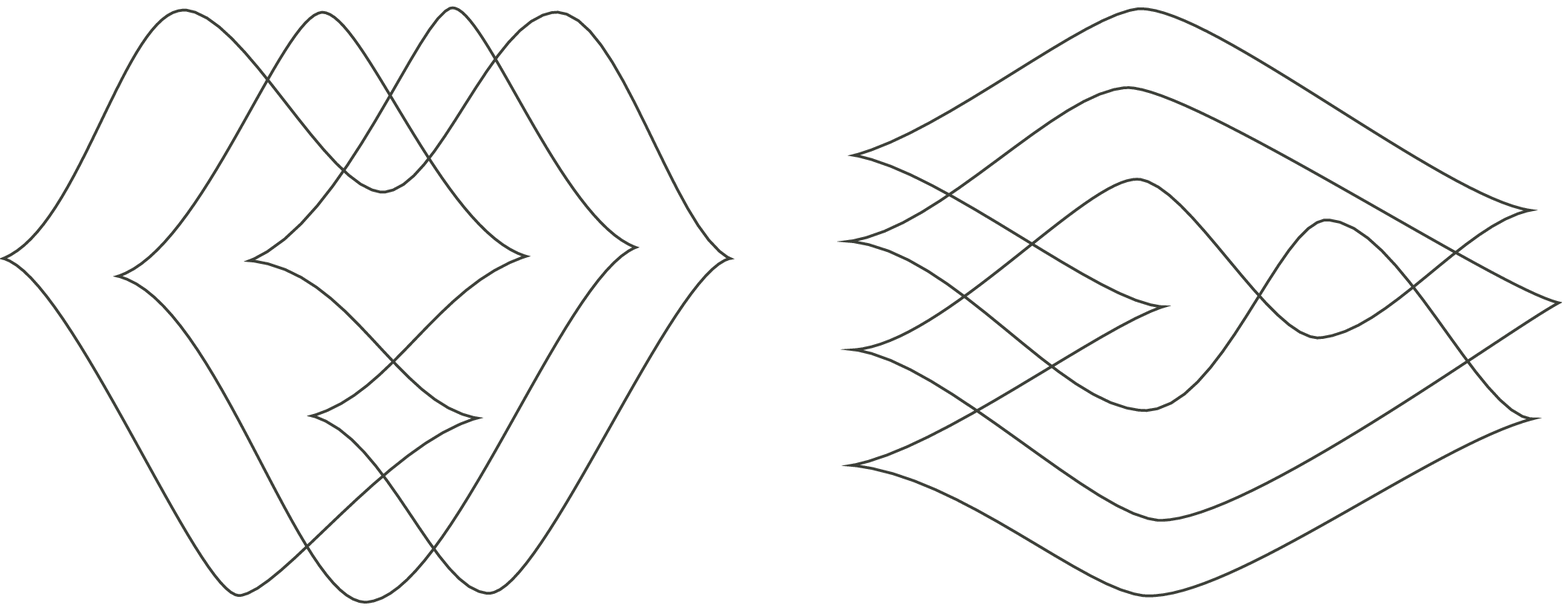}{.14}	&\up{3t}
&\up{(-7,0)}	\\ [1ex]
\up{9_{47}}	&\up{(-2,1)}	&			&\fig{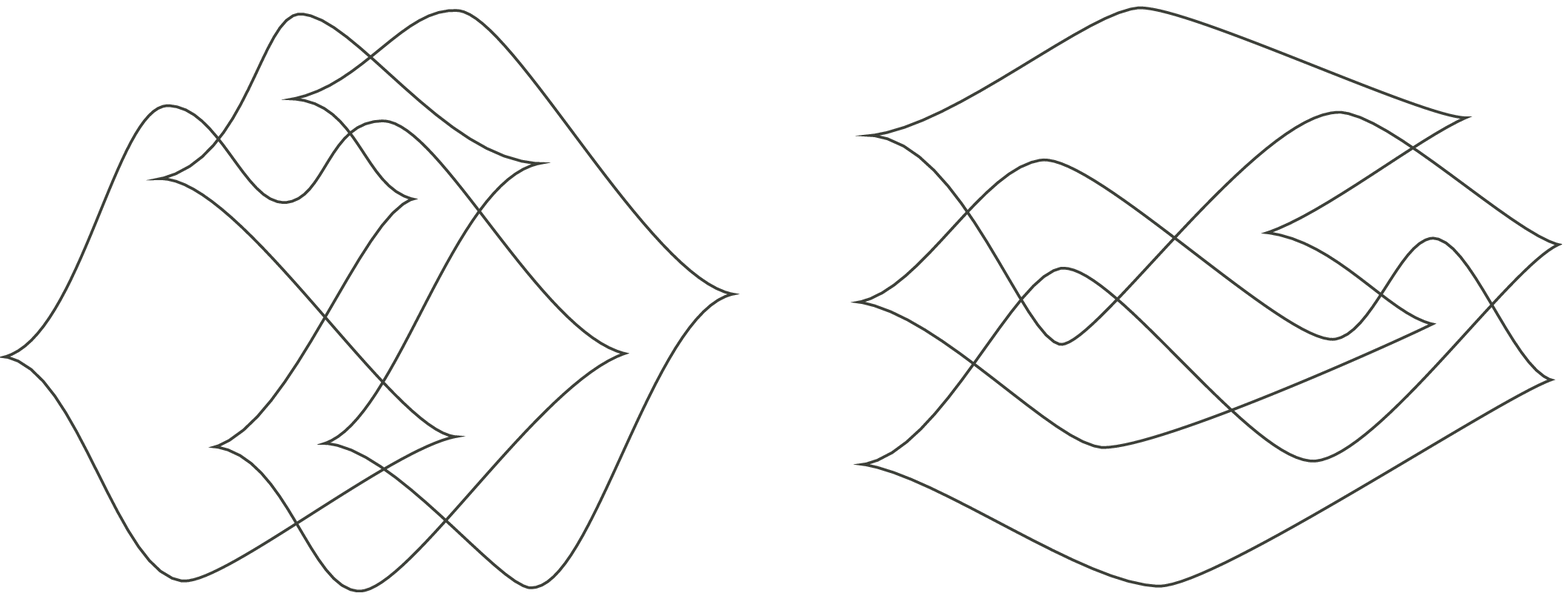}{.14}
&\up{3t}		&\up{(-7,0)}	\\ [1ex]
\up{9_{48}}	&\up{(-1,0)}	&\up{1+t}		&\fig{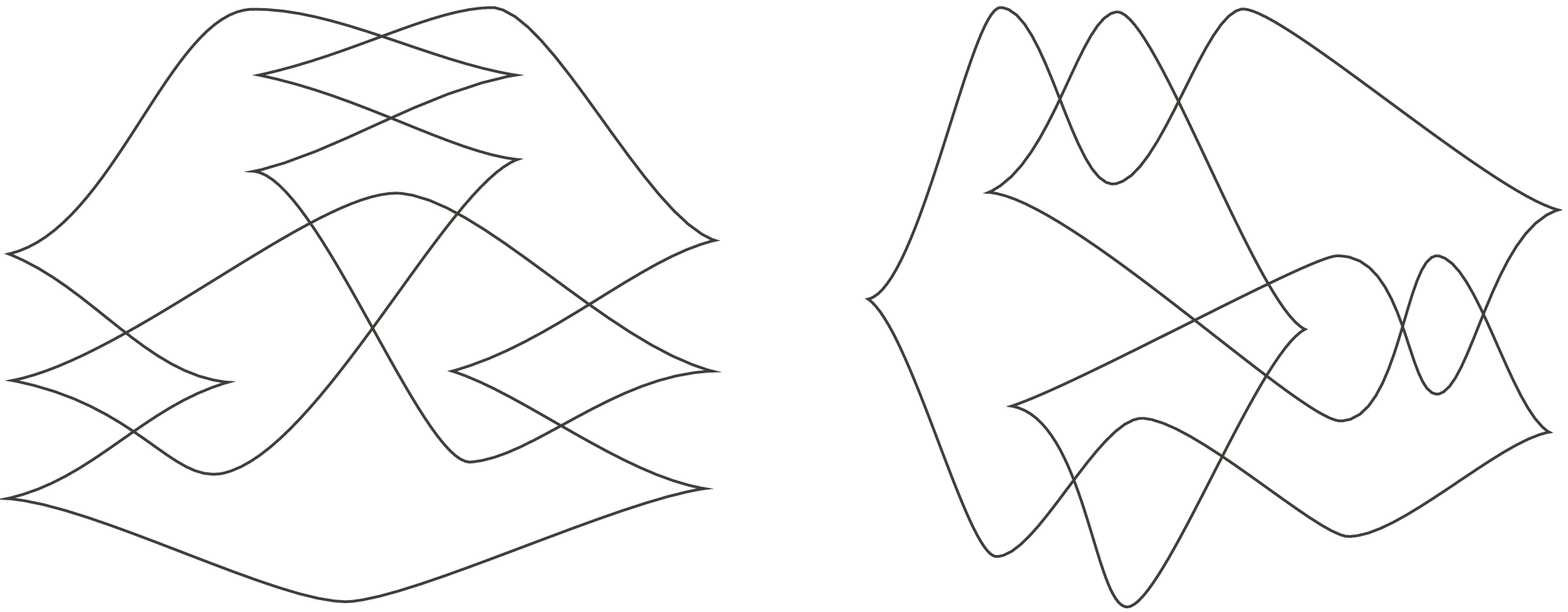}{.14}
&			&\up{(-8,1)}	\\ [1ex]
\up{9_{49}}	&\up{(3,0)}	&\up{2}		&\fig{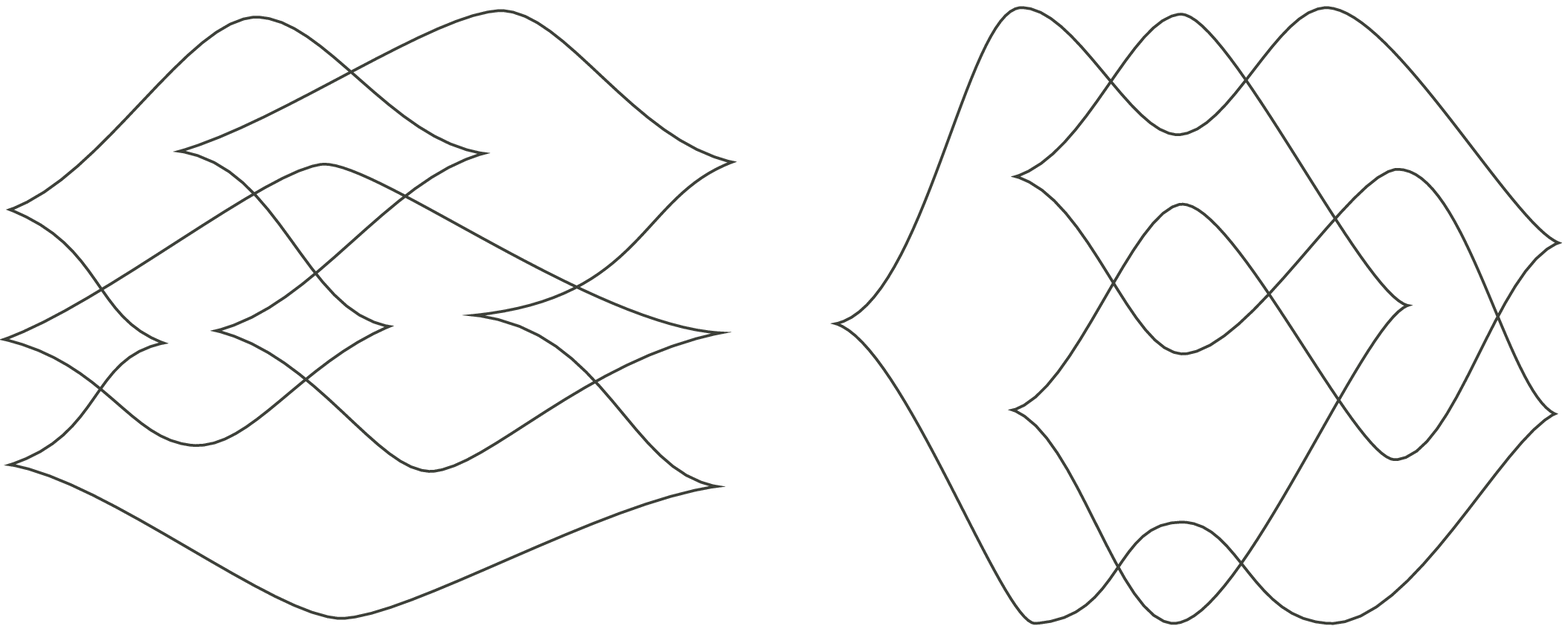}{.14}	&
&\up{(-12,1)}	\\ [1ex]
\hline
\end{tabular}
\end{table}

%%%%%%%%%%%%%
%\end{comment}
%%%%%%%%%%%%%

\end{document}